\documentclass[10pt]{amsart}

\usepackage{amsmath,amsthm,amssymb}
\usepackage[mathscr]{eucal}
\usepackage{graphicx,array,pinlabel,enumerate}
\usepackage[small,bf,margin=24pt]{caption}

\theoremstyle{plain}
\newtheorem{tetel}{Theorem}[section]
\newtheorem{all}[tetel]{Proposition}
\newtheorem{lemma}[tetel]{Lemma}
\newtheorem{kov}[tetel]{Corollary}
\newtheorem{sejt}[tetel]{Conjecture}
\theoremstyle{definition}\newtheorem{Def}[tetel]{Definition}
\theoremstyle{remark}\newtheorem{megj}[tetel]{Remark}
\newtheorem{pelda}[tetel]{Example}

\newcommand*{\N}{\ensuremath{\mathbf N}}
\newcommand*{\R}{\ensuremath{\mathbf R}}
\newcommand*{\Z}{\ensuremath{\mathbf Z}}

\DeclareMathOperator{\bip}{Bip}
\DeclareMathOperator{\conv}{Conv}

\begin{document}

\title{A version of Tutte's polynomial for hypergraphs}

\author{Tam\'as K\'alm\'an}

\address{Tokyo Institute of Technology}
\email{kalman@math.titech.ac.jp}
\urladdr{www.math.titech.ac.jp/\char126 kalman}

\date{}


\thanks{I acknowledge support by a Japan Society for the Promotion of Science (JSPS) Grant-in-Aid for Young Scientists B (no.\ 21740041).}

\begin{abstract}
Tutte's dichromate $T(x,y)$ is a well known graph invariant. Using the original definition in terms of internal and external activities as our point of departure, we generalize the valuations $T(x,1)$ and $T(1,y)$ to hypergraphs. In the definition, we associate activities to hypertrees, which are 
generalizations of the indicator function of the edge set of a spanning tree. We prove that hy\-per\-trees form a lattice polytope which is the set of bases in a polymatroid. In fact, we extend our 
invariants to integer polymatroids as well.
We also examine hypergraphs that can be represented by planar bipartite graphs, write their hypertree polytopes in the form of a determinant, and prove a duality property that leads to an extension of Tutte's Tree Trinity Theorem.
\end{abstract}

\maketitle

\section{Introduction}

In this paper we investigate how much of the theory of the Tutte polynomial $T_G(x,y)$ \cite{tutte2}, which is a large and important branch of graph and matroid theory, can be generalized to hypergraphs and polymatroids. We will see that some of the basic definitions carry over in a new and interesting way. From among the many important properties of the Tutte polynomial, we find an elegant generalization of the equation relating planar dual graphs. The deletion--contraction rule also has an extension to our context but, interestingly, such formulas play a lesser role here than in the classical case.

The central object of our theory is the lattice polytope $Q_\mathscr H$, associated to a hypergraph $\mathscr H$, which generalizes the well known spanning tree polytope of ordinary graphs. The lattice points in $Q_\mathscr H$ are called hypertrees. A hypertree, see Definition \ref{def:hiperfa}, is essentially the valence distribution, taken at one side of the bipartition, of a spanning tree in the bipartite graph naturally associated to $\mathscr H$. It is the use of this concept that sets apart the current paper from earlier work.

From $Q_\mathscr H$, we will read off the one-variable polynomial invariants $I_\mathscr H$ and $X_\mathscr H$ of $\mathscr H$ which generalize $T_G(x,1)$ and $T_G(1,y)$, respectively. They are called the interior and exterior polynomials. Both have positive integer coefficients, the sum of which is the number of hypertrees. Their definitions involve a direct, yet non-obvious generalization of Tutte's idea of ordering the (hyper)edges and using the order to define their internal and external activities with respect to (hyper)trees.

The hypertree polytope is 
naturally embedded as
the set of bases in a certain integer polymatroid. The latter can be viewed as a generalization to hypergraphs of the cycle matroid of a graph. The interior and exterior polynomials can in fact be defined for all integer (extended) polymatroids, in other words for all integer-valued submodular set functions regardless of whether they are non-decreasing. In this paper we will merely mention this possibility as we intend to keep the focus on bipartite graphs and  hypergraphs.

As the counterpart of the planar duality relation $T_{G^*}(x,y)=T_G(y,x)$ we obtain the formula $I_{\mathscr H^*}=X_\mathscr H$, where $\mathscr H^*$ is a naturally associated dual to $\mathscr H$ in cases when the hypergraph $\mathscr H$ can be represented with a plane bipartite graph drawing.

Our invariants have another fundamental property which does not manifest itself in the classical case. Namely, any hypergraph $\mathscr H$ has an abstract dual $\overline{\mathscr H}$ that results from interchanging the roles of its vertices and hyperedges. For such a pair, we conjecture the identity $I_\mathscr H=I_{\overline{\mathscr H}}$. 

A weaker statement is that $\mathscr H$ and $\overline{\mathscr H}$ have the same number of hypertrees. This was recently proven by Alexander Postnikov. Indeed the discovery of $Q_\mathscr H$ (and its dual relationship with $Q_{\overline{\mathscr H}}$) should be attributed to him. In his beautiful paper \cite{post} he puts these polytopes in several important contexts.
The polynomials $I_\mathscr H$ and $X_\mathscr H$ however seem to appear for the first time in this article. I read Postnikov's work after most results presented here have been obtained. Prior to that I was only able to prove that $Q_\mathscr H$ and $Q_{\overline{\mathscr H}}$ had the same number of lattice points when $\mathscr H$ was planar. This brings us to the second half of the paper.

We will revisit a delicate 
picture that William Tutte \cite{tutte1} introduced before his discovery of the dichromate. It consists of three plane bipartite graphs that together triangulate the sphere $S^2$. Their planar duals are canonically directed so that a so-called arborescence number can be associated to each. Tutte's Tree Trinity Theorem states that in such a triple, the three arborescence numbers coincide. We will see that this value is also the number of hypertrees in any of the six hypergraphs found in Tutte's picture. The six hypertree polytopes, for which we will derive a concise determinant formula based on work of Kenneth Berman, form three centrally symmetric pairs. The total number of interior and exterior polynomials associated to the six hypergraphs is reduced from $12$ to $6$ by our planar duality result, and the conjecture above on abstract duality implies that the number is in fact only $3$.

The notions that we are about to introduce also have a strong connection to low-dimensional topology. Indeed, the results contained in this paper are by-products of the author's investigations in knot theory, especially of research done on the Homfly polynomial.
There is a well known method, called the median construction, that one may apply to a plane graph to get an alternating link diagram. When the graph is bipartite, the link has a natural orientation. The interior polynomial was originally developed to describe a pattern found in the Homfly polynomials of these associated links. Then in joint work with Andr\'as Juh\'asz and Jacob Rasmussen, we found a manifestation of the hypertree polytope in Heegaard Floer theory. But as the former of these two results remains a conjecture at the time of this writing, it feels prudent to save the exploration of topological aspects to future papers.

In another forthcoming paper, joint with Alex Bene, we extend the theory of the hypertree polytope and interior and exterior polynomials to topological bipartite graphs and hypergraphs, in the same way that the Bollob\'as--Riordan polynomial generalizes the Tutte polynomial for fatgraphs.

The paper is organized as follows. Section \ref{sec:prelim} fixes terminology and recalls Tutte's definition of the dichromate. 
In Section \ref{sec:hyper}, we introduce (that is to say, recall from \cite{post}) the hypertree polytope and in Section \ref{sec:geometry} we examine its basic geometry. After all this preparation, in section \ref{sec:poly} we define interior and exterior polynomials. The material on polymatroids is found in subsections \ref{ssec:polymat} and \ref{ssec:polypoly}. Section \ref{sec:prop} establishes a variety of properties of the new invariants.
In Section \ref{sec:sejtes} we state the main unsolved problem of the paper, the abstract duality conjecture, along with some supporting evidence. We discuss planar duality in Section \ref{sec:dual}. Finally, Section \ref{sec:trinity} introduces tri\-ni\-ties of plane bipartite graphs and some earlier work on them, and in Section \ref{sec:moretrinity} we lay out our new results on that subject.

{\bf Acknowledgements.} Part of the research reported here was completed while the author worked at the University of Tokyo. My understanding of issues pertaining to the current paper has greatly benefited from discussions with my collaborators Alex Bene, Andr\'as Juh\'asz, and Jake Rasmussen. I also want to express my gratitude to L\'aszl\'o Feh\'er, L\'aszl\'o Lov\'asz, Hitoshi Murakami, Rich\'ard Rim\'anyi, Oliver Riordan, G\'abor Tardos, and Vera V\'ertesi for stimulating discussions. Last but not least, I thank P\'eter Juh\'asz for checking many examples by computer.


\section{Preliminaries}\label{sec:prelim}

We will use the standard notion of a \emph{graph} as an ordered pair $G=(V,E)$ with finite vertex set $V$ and finite edge (multi)set $E$. Loop edges and multiple edges are allowed.
By $G'=(V',E')$ being a \emph{subgraph} of $G$
we will mean 
$V'=V$ and $E'\subset E$ (thus subgraphs of a given graph are equivalent to their edge sets). For example, if a subgraph contains no edge adjacent to some vertex $v$, then $\{v\}$ (more precisely, $(\{v\},\varnothing)$) is a connected component of the subgraph.

We will write $G\setminus\{e\}$ for the graph that results from removing (one copy of) the edge $e$ from $E$. (In \eqref{eq:delcontr} below and in subsection \ref{ssec:delcontr} this will also appear as $G-e$.) The symbol $G-\{v\}$ will stand for the graph that remains after removing the vertex $v$ and all its adjacent edges from $G$.

Graphs can be viewed as one-dimensional cell complexes and in that sense the \emph{nullity} of a graph, $n(G)$, is the rank of its first homology group. (Thus we may also refer to the nullity as the \emph{first Betti number}.)

\subsection{Review of the Tutte polynomial}\label{oldtutte}
Tutte's original construction \cite{tutte2} of the two-variable polynomial $T_G(x,y)$, associated to the graph $G=(V,E)$, proceeds as follows. Order $E$ arbitrarily.
Consider the set $\mathscr T$ of \emph{spanning trees}, that is connected and cycle-free subgraphs, for $G$. (In order for $\mathscr T$ to be non-empty, $G$ needs to be connected. This will almost always be assumed.) 

\begin{Def}\label{def:oldactive}
Given a spanning tree $\Gamma\in\mathscr T$, call an edge $e\in\Gamma$ \emph{internally active} with respect to $\Gamma$ if, after removing $e$ from $\Gamma$, connectedness of the subgraph cannot be restored by adding an edge to $\Gamma\setminus\{e\}$ that is smaller than $e$. 

An edge $e\not\in\Gamma$ is \emph{externally active} if, after adding $e$ to $\Gamma$, cycle-freeness cannot be restored by removing an edge from $\Gamma\cup\{e\}$ that is smaller than $e$. 
\end{Def}

In fact Tutte put `larger' instead of `smaller' in both cases above. But as we will see in Theorem \ref{thm:tutte} below, reversing the order does not affect our main notion.

\begin{Def}\label{def:tuttepoly}
Let $\iota(\Gamma)$ and $\varepsilon(\Gamma)$ denote the number of internally and externally active edges, respectively, with respect to $\Gamma$ (under the given order). Then, the \emph{Tutte polynomial} or \emph{dichromate} of the 
graph $G$ is the generating function
\[T_G(x,y)=\sum_{\Gamma\in\mathscr T}x^{\iota(\Gamma)}y^{\varepsilon(\Gamma)}.\]
\end{Def}

\begin{tetel}[Tutte \cite{tutte2}]\label{thm:tutte}
The polynomial $T_G(x,y)$ is independent of the chosen edge order. 
\end{tetel}

A multitude of properties of $T_G(x,y)$ is known. One of them is the deletion--contraction relation
\begin{equation}\label{eq:delcontr}
\begin{array}{rcll}
T_G&=&xT_{G/e}&\text{if }e\text{ is a bridge}\\
T_G&=&yT_{G-e}&\text{if }e\text{ is a loop}\\
T_G&=&T_{G-e}+T_{G/e}&\text{if }e\text{ is neither a bridge nor a loop.}
\end{array}
\end{equation}
(A bridge is an edge whose removal disconnects the graph. The contraction $G/e$ will be formally defined in subsection \ref{ssec:delcontr}.) For an excellent survey on why the Tutte polynomial is central in graph theory, see \cite{em}.

\begin{pelda}\label{ex:egy}
We show a graph and its Tutte polynomial.
\unitlength .04in
\[\begin{picture}(36,22)
\put(1,6){\circle*{2}}
\put(1,16){\circle*{2}}
\put(11,1){\circle*{2}}
\put(11,11){\circle*{2}}
\put(11,21){\circle*{2}}
\put(21,6){\circle*{2}}
\put(21,16){\circle*{2}}
\thicklines
\put(1,6){\line(0,1){10}}
\put(11,11){\line(0,1){10}}
\put(21,6){\line(0,1){10}}
\put(1,6){\line(2,1){10}}
\put(1,16){\line(2,1){10}}
\put(11,1){\line(2,1){10}}
\put(1,6){\line(2,-1){10}}
\put(11,11){\line(2,-1){10}}
\put(11,21){\line(2,-1){10}}
\put(25,11){\vector(1,0){10}}
\put(25,10){\line(0,1){2}}
\end{picture}
\begin{array}[b]{l@{\hspace{2pt}}l@{\hspace{2pt}}l@{\hspace{2pt}}l@{\hspace{2pt}}l@{\hspace{2pt}}l@{\hspace{2pt}}l}
+y^3&&&&&&\\
+3y^2&+3xy^2&&&&&\\
+2y&+7xy&+6x^2y&+3x^3y&&&\\
&+2x&+6x^2&+7x^3&+6x^4&+3x^5&+x^6.
\end{array}
\]
\end{pelda}

\subsection{Hypergraphs and bipartite graphs} 
For the purposes of this paper these structures are almost equivalent, as follows.

\begin{Def} A \emph{bipartite graph} is a triple $G=(V_0,V_1,E)$, where $V_0$ and $V_1$ are disjoint finite sets, called \emph{color classes}, and $E$ is a finite set of edges, each connecting an element of $V_0$ to an element of $V_1$ (multiple edges are not allowed). We will treat $(V_0,V_1,E)$, $(V_1,V_0,E)$, and the graph $(V_0\cup V_1,E)$ as the same object.
\end{Def}

\begin{Def} A \emph{hypergraph} is a pair $\mathscr H=(V,E)$, where $V$ is a finite set and $E$ is a finite multiset of non-empty subsets of $V$. Elements of $V$ are called \emph{vertices} and the elements of $E$ are the \emph{hyperedges}. 
\end{Def}

Thus, hyperedges with multiplicity (that is, several copies of the same subset of $V$) are allowed. On the other hand, obviously, each hyperedge contains each of its elements exactly once.

A hypergraph is both a generalization and a special case of a graph. The first point is obvious. Conversely, the sets $V$ and $E$ that constitute the hypergraph $\mathscr H$ may be viewed as the color classes of a bipartite graph $\bip\mathscr H$: we connect $v\in V$ to $e\in E$ with an edge if $v\in e$. We will refer to the result as the \emph{bipartite graph associated to the hypergraph} $\mathscr H$.


The construction of $\bip\mathscr H$ above is reversible if we specify one of the two color classes in the bipartite graph $G=(V_0,V_1,E)$. We will use the notation
\begin{equation}\label{eq:absdual}
\mathscr G_0=(V_1,V_0)\quad\text{and}\quad\mathscr G_1=(V_0,V_1)
\end{equation}
for the resulting pair of hypergraphs. (The index of the hypergraph is chosen to emphasize its hyperedges.) Note that $\bip\mathscr G_0=\bip\mathscr G_1=G$.

\begin{Def}
The bipartite graph $G$ above is said to \emph{induce} the hypergraphs $\mathscr G_0$ and $\mathscr G_1$. Two hypergraphs are called \emph{abstract duals} if they can be obtained in the form \eqref{eq:absdual}. In other words, the abstract dual $\overline{\mathscr H}=(E,V)$ of a hypergraph $\mathscr H=(V,E)$ is defined by interchanging the roles of its vertices and hyperedges.
\end{Def}


\section{Hypertrees}\label{sec:hyper}

To generalize the approach of Subsection \ref{oldtutte} to hypergraphs, first we need a notion corresponding to spanning trees. The rest of the paper is built around this concept, and almost all novelty contained herein stems from its use.

\begin{Def}\label{def:hiperfa}
Let $\mathscr H=(V,E)$ be a hypergraph so that its associated bipartite graph $\bip\mathscr H$ is connected. By a \emph{hypertree} in $\mathscr H$ we mean a function (vector) $f\colon E\to\N=\{\,0,1,\ldots\,\}$ so that a spanning tree of $\bip\mathscr H$ can be found which has degree $f(e)+1$ at each $e\in E$. Such a spanning tree is said to \emph{realize} or to \emph{induce} $f$.
\end{Def}

If $f$ is a hypertree, then $0\le f(e)\le|e|-1$ for each $e\in E$. The condition that $\bip\mathscr H$ be connected is 
not essential. We could talk about `hyperforests' realized by spanning forests but, partly to avoid such
terminology, we will generally assume that the bipartite graphs associated to our hypergraphs are connected.

\begin{megj}\label{rem:altalanosit}
Hypertrees generalize spanning trees of graphs: there, an edge $e$ is in the tree if $f(e)=1$ and not in the tree if $f(e)=0$. In our case, we allow hyperedges to be in hypertrees only ``partially.''
\end{megj}

The number of realizations may vary greatly from hypertree to hypertree. However this phenomenon will not be incorporated to our theory in any way. For example, if $f$ is a hypertree and $f(e)=0$, then any edge of $\bip\mathscr H$ adjacent to the hyperedge $e$ can be chosen to be part of a realization of $f$, irrespective of the rest of the construction. The next lemma claims only slightly more.

\begin{lemma}\label{lem:anchor}
Let $\mathscr H=(V,E)$ be a hypergraph and $f\colon E\to\N$ a hypertree. For any collection of edges $\alpha_e$ of $\bip\mathscr H$ so that $\alpha_e$ is adjacent to $e$ for all $e\in E$, there is a realization of $f$ that contains $\alpha_e$ for all $e$.
\end{lemma}

\begin{proof}
Start with an arbitrary realization $\Gamma$ of $f$. If $e$ is such that $\alpha_e\not\in\Gamma$, then add $\alpha_e$ to $\Gamma$. This creates a unique cycle in the subgraph which goes through $e$. Now if we remove the edge of the cycle that is adjacent to $e$ but different from $\alpha_e$, the result is another realization of $f$. Carry this process out at every hyperedge.
\end{proof}

\subsection{The hypertree polytope}\label{ssec:politop}
It turns out that hypertrees are exactly the integer lattice points in a convex polytope. This fact was first realized by Postnikov \cite{post}. But as the author discovered it independently and our points of view and proofs are different, we shall give a full account here. This will be followed by a sampling of Postnikov's ideas, including a sketch of the proof of his duality theorem.

To describe the polytope, let $E'\subset E$ be a non-empty subset of hyperedges and let $\bigcup E'\subset V$ denote the set of its neighbors in $\bip\mathscr H$. Let $\bip\mathscr H\big|_{E'}$, the \emph{bipartite graph induced by $E'$}, be the graph with color classes $E'$ and $\bigcup E'$ and edges inherited from $\bip\mathscr H$. Let us denote the number of its connected components by $c(E')$. 

\begin{tetel}\label{thm:politop}
Let $\mathscr H=(V,E)$ be a hypergraph so that its associated bipartite graph $\bip\mathscr H$ is connected. The hypertrees $f\colon E\to\N$ of $\mathscr H$ are exactly the integer solutions of the following system of linear inequalities in $\R^E$:
\begin{subequations}\label{eq:hypertree}
\begin{alignat}{2}
f(e)&\ge0&\quad&\text{for all }e\in E\label{eq:hypertree1}\\
\sum_{e\in E'}f(e)&\le|\textstyle\bigcup E'|-c(E')&\quad&\text{for all non-empty }E'\subset E\label{eq:hypertree2}\\
\sum_{e\in E}f(e)&=|V|-1.\label{eq:hypertree3}
\end{alignat}
\end{subequations}
\end{tetel}

\begin{megj}\label{rem:-1}
We get an equivalent system of inequalities if we replace \eqref{eq:hypertree2} with $\sum_{e\in E'}f(e)\le|\bigcup E'|-1$, still required for all non-empty $E'\subset E$. This is because if $\bip\mathscr H\big|_{E'}$ is not connected, then its connected components are also of the form $\bip\mathscr H\big|_{E''}$ and we get the stronger version of the inequality by summing the inequalities associated to these smaller subsets $E''$.
\end{megj}

\begin{megj}\label{rem:nemnegkijon}
The conditions \eqref{eq:hypertree1} follow from \eqref{eq:hypertree2} and \eqref{eq:hypertree3}: this is obvious if $E=\{e\}$ and otherwise, with $E'=E\setminus\{e\}$, we have
\[f(e)=|V|-1-\sum_{e'\in E'}f(e')\ge(|V|-|\textstyle\bigcup E'|)+(c(E')-1),\]
where the right hand side is the sum of two non-negative quantities.
\end{megj}

\begin{proof}[Proof of Theorem \ref{thm:politop}]
The conditions \eqref{eq:hypertree} are necessary because of the well known facts that 
\begin{enumerate}[(i)]
\item\label{a} the number of edges in a spanning forest (maximal cycle-free subgraph) of a graph is the number of vertices minus the number of connected components of the graph and 
\item\label{b} any cycle-free subgraph is part of a spanning forest. 
\end{enumerate}
Indeed, for a given hypertree $f$ and non-empty subset $E'\subset E$, any spanning tree of $\bip\mathscr H$ that realizes $f$ has an intersection with $\bip\mathscr H\big|_{E'}$ which is a cycle-free subgraph of the latter. As such, it may have at most $|E'|+|\bigcup E'|-c(E')$ edges. Since each of those has exactly one of its ends at an element of $E'$, we have
\[\sum_{e\in E'}(f(e)+1)\le|E'|+|\bigcup E'|-c(E'),\] 
which is just another form of \eqref{eq:hypertree2}. The claim \eqref{eq:hypertree1} is obvious and \eqref{eq:hypertree3} is immediate from \eqref{a} above.

To see why \eqref{eq:hypertree} is also sufficient, let us formulate a lemma that is marginally stronger than what we currently need, but the form in which we state it will be useful later. Recall that 
a graph has nullity zero if and only if it is cycle-free.

\begin{lemma}\label{lem:kormentes}
Suppose that the integer vector $f\colon E\to\N$ satisfies conditions \eqref{eq:hypertree1} and \eqref{eq:hypertree2} (including for $E'=E$) but not necessarily  \eqref{eq:hypertree3}. Then there exists a cycle-free subgraph in $\bip\mathscr H$ that has valence $f(e)+1$ at all elements $e\in E$.
\end{lemma}

Let the vector $f$ satisfy the conditions and choose an arbitrary subgraph $P$ of $\bip\mathscr H$ whose valence at each $e\in E$ is $f(e)+1$. This is possible because \eqref{eq:hypertree2} applied to $E'=\{e\}$ says $f(e)\le|e|-1$. If $P$ is cycle-free, we are done. Assume it is not.

It suffices to show that there is another subgraph of $\bip\mathscr H$ that has the same valences at elements of $E$ as $P$ (prescribed by $f$) but whose nullity is one less than that of $P$. The subgraph $P$ has a connected component $C$ which contains a cycle. Let $C$ intersect $E$ in the set $E'$.  Applying \eqref{eq:hypertree2} to $E'$ we see that there is a hyperedge $e\in E'$ which is connected by an edge $\alpha$ of $\bip\mathscr H$ to a vertex which is not in $C$. (Otherwise, $P$ and $\bip\mathscr H\big|_{E'}$ would form a connected and not cycle-free intersection; as such a subgraph contains a spanning tree of $\bip\mathscr H\big|_{E'}$, we would get a contradiction between \eqref{eq:hypertree2} and \eqref{a}.)

Now if $e$ is part of a cycle in $C$, we are done because we may remove an edge of that cycle (adjacent to $e$) and replace it with $\alpha$, whereby valences at hyperedges remain the same but the nullity reduces by one. Otherwise, each edge of $C$ coming out of $e$ leads to a different connected component of $C-\{e\}$. At least one of these (call it $C'$) still contains a cycle. Replace one such edge with $\alpha$. This results in a subgraph $P'$ that has the right valences and the same nullity as $P$, but which has a non-tree component $C'$ containing fewer elements of $E$ than $C$ did. Repeat the procedure of this paragraph to $P'$ and $C'$. After a finite number of iterations, the desired reduction in the nullity will occur. This finishes the proof of the lemma.

If we apply Lemma \ref{lem:kormentes} to a vector $f$ that satisfies \eqref{eq:hypertree}, the resulting subgraph is not just cycle-free but because of \eqref{eq:hypertree3} it is actually a spanning tree.
\end{proof}

\begin{Def}
Let $\mathscr H$ be a hypergraph. The \emph{hypertree polytope} of $\mathscr H$, denoted with $Q_\mathscr H$, is the set of solutions in $\R^E$ of the inequalities \eqref{eq:hypertree}. 
\end{Def}

By Theorem \ref{thm:politop}, the set of hypertrees in $\mathscr H$ can be written as $Q_\mathscr H\cap\Z^E$. Furthermore, $Q_\mathscr H$ is a lattice polytope, that is $Q_\mathscr H=\conv(Q_\mathscr H\cap\Z^E)$, where $\conv$ denotes the usual convex hull.

The hypertree polytope generalizes the usual spanning tree polytope of a graph (cf.\ Remarks \ref{rem:altalanosit} and \ref{rem:-1}). Equation \eqref{eq:hypertree3} means that $Q_\mathscr H$ is part of an affine hyperplane in $\R^E$. Thus, $Q_\mathscr H$ is situated in the lattice cut out from $\Z^E$ by the hyperplane. 

\begin{pelda}\label{ex:ketto}
The graph used in Example \ref{ex:egy} is in fact bipartite. Let us denote it with $G$ and label its color classes with $V_0=\{\,a,b,c\,\}$ and $V_1=\{\,p,q,r,s\,\}$ according to Figure \ref{fig:betuzes}.
\begin{figure}[htbp]
\unitlength 3pt
\begin{picture}(36,30)
\put(8,9){\circle*{2}}
\put(8,19){\circle*{2}}
\put(18,4){\circle*{2}}
\put(18,14){\circle*{2}}
\put(18,24){\circle*{2}}
\put(28,9){\circle*{2}}
\put(28,19){\circle*{2}}
\thicklines
\put(8,9){\line(0,1){10}}
\put(18,14){\line(0,1){10}}
\put(28,9){\line(0,1){10}}
\put(8,9){\line(2,1){10}}
\put(8,19){\line(2,1){10}}
\put(18,4){\line(2,1){10}}
\put(8,9){\line(2,-1){10}}
\put(18,14){\line(2,-1){10}}
\put(18,24){\line(2,-1){10}}
\put(17,10){$q$}
\put(4,19){$p$}
\put(17,0){$r$}
\put(30,19){$s$}
\put(17,26){$c$}
\put(4,8){$a$}
\put(30,8){$b$}
\end{picture}
\caption{A bipartite graph.}
\label{fig:betuzes}
\end{figure}
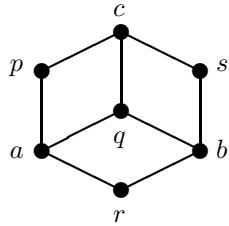
There are $T_G(1,1)=50$ spanning trees in $G$. They give rise to seven hypertrees in $\mathscr G_0$, and the same number in $\mathscr G_1$, cf.\ Theorem \ref{thm:post}. These are shown, with a concrete spanning tree realizing each, in Figure \ref{fig:ujfak}. Because of \eqref{eq:hypertree3}, it is natural to view the polytopes $Q_{\mathscr G_0}$ and $Q_{\mathscr G_1}$ in barycentric coordinate systems with basepoints labeled with
\[\mathbf a(3,0,0),\,\mathbf b(0,3,0),\,\mathbf c(0,0,3)\] 
and 
\[\mathbf p(2,0,0,0),\,\mathbf q(0,2,0,0),\,\mathbf r(0,0,2,0),\,\mathbf s(0,0,0,2),\]
respectively. In Figure \ref{fig:polytopes}, we indicated individual hypertrees by dots and the two hypertree polytopes by thickened edges.

\begin{figure}[htbp]
\labellist
\small
\pinlabel $\mathbf a$ at -10 15
\pinlabel $\mathbf b$ at 340 15
\pinlabel $\mathbf c$ at 160 270
\pinlabel $\mathbf q$ at 575 290
\pinlabel $\mathbf r$ at 455 75
\pinlabel $\mathbf s$ at 710 -10
\pinlabel $\mathbf p$ at 770 150
\endlabellist
   \centering
   \includegraphics[width=4in]{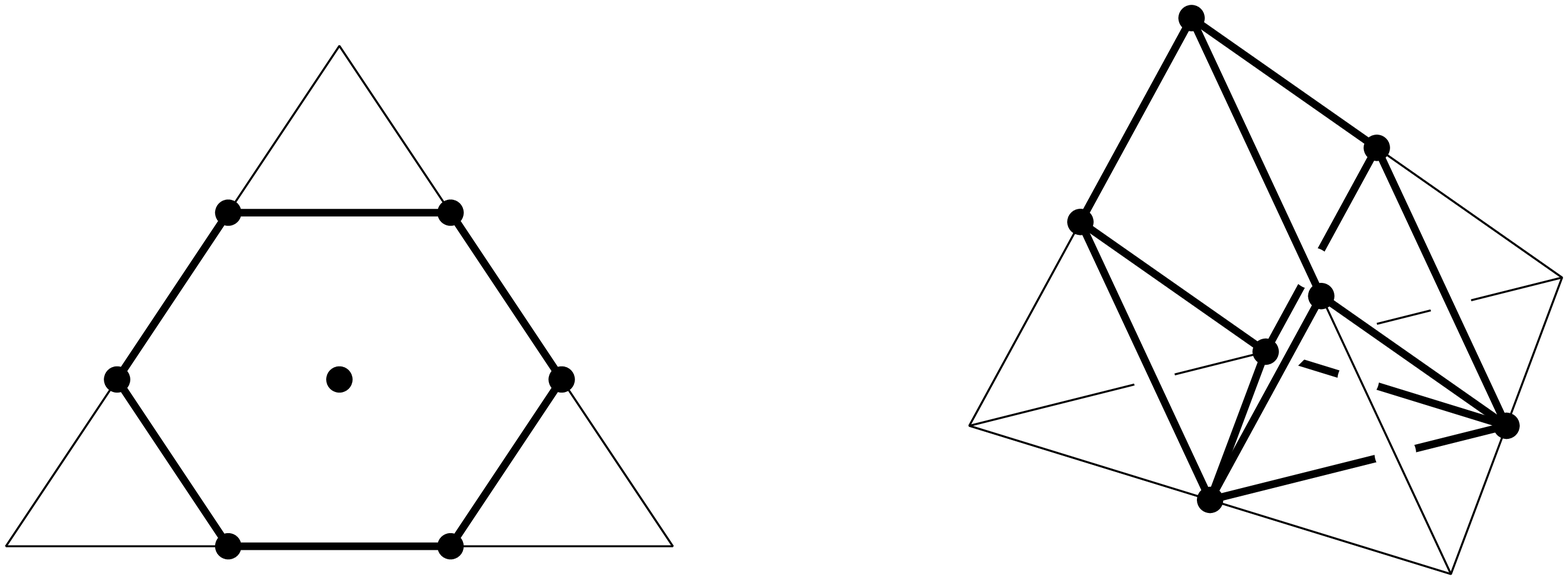} 
   \caption{Hypertree polytopes associated to a pair of abstract dual hypergraphs.}
   \label{fig:polytopes}
\end{figure}
\end{pelda}

The hypertree polytope is of course a polytope, i.e., convex, bounded, and closed. In particular it is closed under convex linear combinations. For its integer lattice points, that is the hypertrees themselves, convexity translates to the following obvious (from Theorem \ref{thm:politop}) but useful property. 

\begin{lemma}\label{lem:discreteconvex}
Let $\{f_i\}$ be some set of hypertrees in the hypergraph $\mathscr H=(V,E)$ and let $g\colon E\to\N$ be a function so that $g(e)\ge0$ for all $e\in E$ and $\sum_{e\in E}g(e)=|V|-1$. If it is also true that for all non-empty $E'\subset E$, there is an index $i$ so that $\sum_{e\in E'}g(e)\le\sum_{e\in E'}f_i(e)$, then $g$ is a hypertree.
\end{lemma}

\subsection{Postnikov's approach}
In the terminology of \cite{post}, a hypertree is a draconian sequence \cite[Definition 9.2]{post} and the hypertree polytope is 
a trimmed generalized permutohedron \cite[Definition 11.2]{post}. Let us now summarize 
Corollary 11.8, Theorem 12.9, and some other surrounding statements from \cite{post}.

\begin{tetel}[Postnikov]\label{thm:post}
Let $G=(V_0,V_1,E)$ be a connected bipartite graph with associated hypergraphs $\mathscr G_0$ and $\mathscr G_1$ as in \eqref{eq:absdual}. Then $|Q_{\mathscr G_0}\cap\Z^{V_0}|=|Q_{\mathscr G_1}\cap\Z^{V_1}|$. In other words, abstract dual hypergraphs have the same number of hypertrees.
\end{tetel}

It turns out that $Q_{\mathscr G_0}$ and $Q_{\mathscr G_1}$ have even more in common. Conjecture \ref{conj:dual} below proposes a generalization of Postnikov's theorem.

Without going into the fine details of an elegant but complicated proof, we will quote some of Postnikov's observations. Of course, an equivalent statement is 
\[|Q_{\mathscr H}\cap\Z^E|=|Q_{\overline{\mathscr H}}\cap\Z^V|,\]
where $\mathscr H=(V,E)$ is an arbitrary hypergraph and $\overline{\mathscr H}=(E,V)$ is its abstract dual. We will use this latter formulation. Denote the standard generators of $\R^E$ with $u_e$, $e\in E$. Recall that the \emph{Minkowski sum} of two subsets $A,B$ of a vector space is the set of vectors one can obtain as $a+b$ with $a\in A$ and $b\in B$. The \emph{Minkowski difference} $A-B$ is defined as the set of vectors $x$ so that $x+B\subset A$. It can also be thought of as the set of all translates of $B$ that are contained in $A$. See \cite[Chapter 3]{bojler} for a thorough introduction to these operations. 

Next, let us quote \cite[Lemma 11.7]{post}.

\begin{all}[Postnikov]\label{pro:felbont}
The hypertree polytope $Q_{\mathscr H}\subset\R^E$ of the hypergraph $\mathscr H=(V,E)$ can be written as 
\[
Q_{\mathscr H}=\Big(\sum_{v\in V}\Delta_v\Big)-\Delta_E,
\]
where $\Sigma$ denotes the Minkowski sum of the simplices $\Delta_v=\conv\{\,u_e\mid v\in e\,\}$, $\Delta_E=\conv\{\,u_e\mid e\in E\,\}$ is the unit simplex and the right hand side is a Minkowski difference.
\end{all}


Of course, if we denote the standard basis for $\R^V$ with $u_v$, $v\in V$ and define the simplices $\Delta_e=\conv\{\,u_v\mid v\in e\,\}\subset\R^V$ for all $e\in E$, as well as $\Delta_V=\conv\{\,u_v\mid v\in V\,\}$, then
\[Q_{\overline{\mathscr H}}=\Big(\sum_{e\in E}\Delta_e\Big)-\Delta_V.\]

The Minkowski sums (in Postnikov's terminology, the (untrimmed) generalized permutohedra)
\[Q_{\mathscr H}^+=\sum_{v\in V}\Delta_v\quad\text{and}\quad Q_{\overline{\mathscr H}}^+=\sum_{e\in E}\Delta_e\] are related through the so-called root polytope 
\[Q=\conv\{\,u_e+u_v\mid v\in e\,\}\subset\R^E\oplus\R^V\]
(an $(|E|+|V|-2)$-dimensional polytope whose vertices are indexed by the edges of $\bip\mathscr H$) as described below.

\begin{all}\label{pro:csel}
If we define the affine subspaces $S_E$ and $S_V$ of $\R^E\oplus\R^V$ as
\[S_E=\pi_V^{-1}\left(\frac1{|V|}\sum_{v\in V}u_v\right)\quad\text{and}\quad S_V=\pi_E^{-1}\left(\frac1{|E|}\sum_{e\in E}u_e\right),\]
where $\pi_V\colon\R^E\oplus\R^V\to\R^V$ and $\pi_E\colon\R^E\oplus\R^V\to\R^E$ are the projections, then, up to translation,
\[Q_{\mathscr H}^+=|V|\left(Q\cap S_E\right)\quad\text{and}\quad Q_{\overline{\mathscr H}}^+=|E|\left(Q\cap S_V\right).\]
\end{all}

Studying triangulations of $Q$, Postnikov observes that a set of vertices of $Q$ is affinely independent if and only if the corresponding edges form a forest in $\bip\mathscr H$. The simplices $T_\Gamma$ that arise this way from spanning trees $\Gamma$ have the same volume. 

Most importantly, if we form the intersections $T_\Gamma\cap S_E$ and $T_\Gamma\cap S_V$ and identify them with subsets of $Q_{\mathscr H}^+$ and $Q_{\overline{\mathscr H}}^+$, respectively, as described in Proposition \ref{pro:csel}, then those subsets contain unique translates of $\Delta_E$ and $\Delta_V$, respectively, by integer vectors. (Also, they are essentially disjoint from other integer translates of those unit simplices.) In fact, those integer vectors are the two hypertrees induced by $\Gamma$. 

Thus Postnikov finds that any triangulation of $Q$ gives a bijection between the integer lattice points of the Minkowski differences which are the hypertree polytopes. We also see that the volume of the root polytope associated to a bipartite graph is proportional to the number of hypertrees in each of its induced hypergraphs.

\section{The geometry of the hypertree polytope}\label{sec:geometry}

\subsection{Polymatroids}\label{ssec:polymat}
It turns out that for our central concepts, namely the interior and exterior polynomials which will be introduced in the next section, the right level of generality is that of a polymatroid. Moreover, basic submodular function techniques will be useful to simplify our arguments, most notably the proof that our polynomials are well defined, even in the hypergraph case. We will recall some elements of this theory here, using \cite[Volume B, Chapter 44]{sch} as basic reference.

\begin{Def}
Let $S$ be a finite ground set and $\mu\colon\mathscr P(S)\to\R$ a set function, i.e., a function defined on all subsets of $S$. We say that $\mu$ is \emph{submodular} if 
\[\mu(U)+\mu(V)\ge\mu(U\cap V)+\mu(U\cup V)\]
holds for all subsets $U,V$ of $S$. The set function $\mu$ is called \emph{non-decreasing} if $U\subset V\subset E$ implies $\mu(U)\le\mu(V)$.
\end{Def}

For an arbitrary set function $\nu\colon\mathscr P(S)\to\R$, we define the polyhedra
\begin{equation*}\begin{split}
P_\nu=\{\,\mathbf{x}\in\R^S\mid\mathbf x\ge\mathbf0,\;\mathbf x\cdot\mathbf i_U\le\nu(U)\text{ for all }U\subset S\,\}\\\text{and }
EP_\nu=\{\,\mathbf{x}\in\R^S\mid\mathbf x\cdot\mathbf i_U\le\nu(U)\text{ for all }U\subset S\,\}.
\end{split}\end{equation*}
Here $\mathbf x\ge\mathbf0$ means that all entries in $\mathbf x$ are non-negative and $\mathbf i_U$ is the indicator function of the subset $U$ so that the dot product $\mathbf x\cdot\mathbf i_U$ becomes the sum of entries in $\mathbf x$ corresponding to elements of $U$.

\begin{Def}
If $\mu\colon\mathscr P(S)\to\R$ is a submodular set function, then $P_\mu$ and $EP_\mu$ are called the \emph{polymatroid} and \emph{extended polymatroid}, respectively, of $\mu$. A vector $\mathbf x\in EP_\mu$ is called a \emph{base} if $\mathbf x\cdot\mathbf i_S=\mu(S)$. We denote the set of bases with $B_\mu$ and call it the \emph{base polytope}.
\end{Def}

The notions of submodular set function and extended polymatroid are essentially equivalent as described in \cite[Section 44.4]{sch}. By the same token, non-decreasing submodular set functions are equivalent to polymatroids. Integer-valued submodular functions correspond to (extended) polymatroids that are \emph{integer}, meaning that they are the convex hulls of their integer lattice points. As to bases, we note 

\begin{lemma}\label{lem:nincstobbbazis}
If $\mu$ is submodular and non-decreasing, then we have $B_\mu\subset P_\mu$.
\end{lemma}

\begin{proof}
Let $\mathbf x$ be a base, $s\in S$, and $U=S\setminus\{s\}$. Then the $s$-component of $\mathbf x$ is $\mu(S)-\mathbf x\cdot\mathbf i_U\ge\mu(S)-\mu(U)\ge\mu(S)-\mu(S)=0$.
\end{proof}

The following definition expresses a basic property of polymatroids  that will have an important role in our treatment.

\begin{Def}\label{def:tight}
Let $\nu\colon\mathscr P(S)\to\R$ be a set function and $\mathbf x\in EP_\nu$. We say that the subset $U\subset S$ is \emph{tight} at $\mathbf x$ if $\mathbf x\cdot\mathbf i_U=\nu(U)$.

The set function $\nu$ is called \emph{tight} if for all $\mathbf x\in EP_\nu$, the set of subsets that are tight at $\mathbf x$ is closed under taking unions and intersections.
\end{Def}

\begin{all}[Theorem 44.2 in \cite{sch}]\label{pro:szoros}
Submodular set functions are tight. 
\end{all}

Our next claim is on certain $3$-dimensional intersections of the base polytope $B_\mu$.

\begin{lemma}\label{lem:teglalap}
Let $\mu\colon\mathscr P(S)\to\R$ be a submodular set function and let $p$, $q$, $r$, and $s$ be distinct elements of $S$. 
Other than these four, give some fixed values
to all components and denote the resulting subset of $B_\mu$ with $B'$. Then, assuming $B'\ne\varnothing$, the face of $B'$ along which the sum of the $k$ and $l$ components takes its maximum is a (possibly degenerate) rectangle with sides parallel to the vectors $\mathbf i_{\{p\}}-\mathbf i_{\{q\}}$ and $\mathbf i_{\{r\}}-\mathbf i_{\{s\}}$.
\end{lemma}

\begin{proof}
The set $B'$ is cut out of the appropriate $3$-dimensional affine subspace $A$ as the intersection of fourteen half-spaces corresponding to the non-trivial subsets of $\{\,p,q,r,s\,\}$. Six of these have the normal vectors $\mathbf i_{\{p\}}+\mathbf i_{\{q\}}-\mathbf i_{\{r\}}-\mathbf i_{\{s\}}$ and so on so that their intersection is a (possibly degenerate) cuboid $C$; see  Figure \ref{fig:muszaki}. The other eight half-spaces have normal vectors such as $\pm(\mathbf i_{\{p\}}+\mathbf i_{\{q\}}+\mathbf i_{\{r\}}-3\mathbf i_{\{s\}})$ and they cut off pieces of $C$ near its vertices. Our task is to show that in $B'$ there can not remain a segment of positive length from any edge of $C$. 

\begin{figure}[htbp]
\labellist
\small
\pinlabel $\mathbf p$ at 530 1130
\pinlabel $\mathbf q$ at -20 980
\pinlabel $\mathbf r$ at 290 660
\pinlabel $\mathbf s$ at 255 810
\pinlabel $\mathbf x$ at 1130 990
\endlabellist
   \centering
   \includegraphics[width=4in]{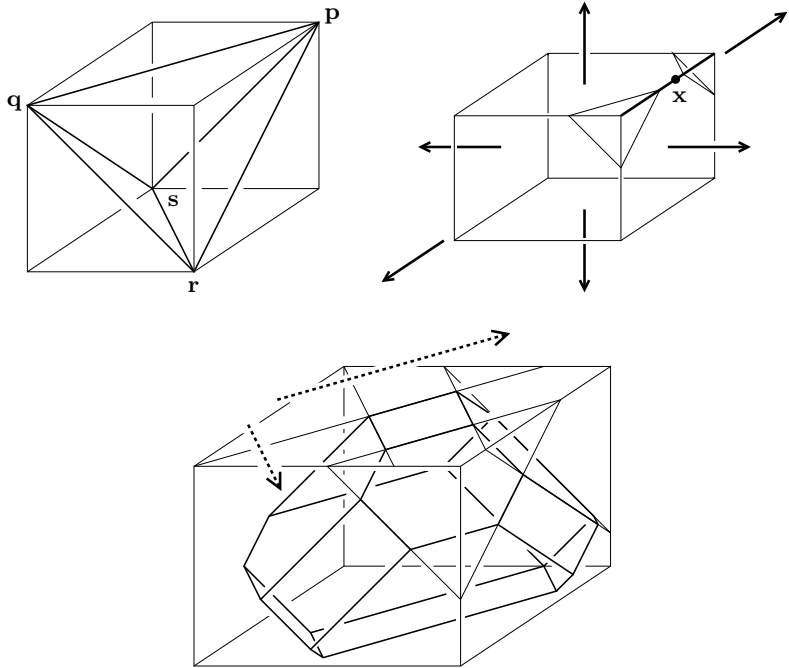} 
   \caption{Upper left: the portion of the simplex $\mu(S)\Delta_S$ that lies in the affine subspace $A$. If $\sigma$ denotes $\mu(S)$ minus the sum of the fixed components then, writing them in the order $p,q,r,s$, the other four components of the points shown are $\mathbf p(\sigma,0,0,0)$ and so on. Upper right: The cuboid $C$ with the normal vectors of its faces. Bottom: the polytope $B'=B_\mu\cap A$ and the direction vectors of two of its edges.}
   \label{fig:muszaki}
\end{figure}

If the opposite was the case and, for instance, such an interior point $\mathbf x$ remained on, say, the edge of $C$ that we thickened in Figure \ref{fig:muszaki}, then $\mathbf x$ would be a base in $B_\mu$ at which the sets $\{\,p,q\,\}$ and $\{\,r,s\,\}$ are tight but neither their union $\{\,p,q,r\,\}$ nor their intersection $\{p\}$ is tight. As this is impossible, 
we see that the correct form of $B'$ is as shown in the lower panel of Figure \ref{fig:muszaki}.
\end{proof}

Submodular functions are relevant in this paper due to the following fact, pointed out to the author by L\'aszl\'o Lov\'asz.

\begin{all}\label{pro:submodular}
Let $\mathscr H=(V,E)$ be a hypergraph. The function 
\begin{equation}\label{eq:subm}
\mu(E')=|\textstyle\bigcup E'|-c(E')
\end{equation}
of Theorem \ref{thm:politop}, extended to the empty set as $\mu(\varnothing)=0$, is non-decreasing and submodular on the set $E$ of hyperedges.
\end{all}

\begin{proof}
We leave the proof of the first assertion to the reader. As to the second, according to \cite[Theorem 44.1]{sch}, it suffices to prove
\begin{equation}\label{eq:laci}
\big(\mu(E'\cup\{e_1\})-\mu(E')\big)+\big(\mu(E'\cup\{e_2\})-\mu(E')\big)\ge\mu(E'\cup\{\,e_1,e_2\,\})-\mu(E')
\end{equation}
for all $E'\subset E$ and distinct $e_1,e_2\in E\setminus E'$. Assume that from the connected components of $\bip\mathscr H\big|_{E'}$, the number of those that are connected to both $e_1$ and $e_2$ is $m_{12}$. Aside from these, let $e_i$, $i=1,2$ be connected to another $m_i$ components of $\bip\mathscr H\big|_{E'}$. Furthermore, let us denote the number of such common elements of $e_1$ and $e_2$ that are not in $\bigcup E'$ with $n_{12}$ and assume that other than these, $e_i$ contains $n_i$ non-elements of $\bigcup E'$. Then we have
\begin{equation*}\begin{split}
\mu(E'\cup\{e_1\})-\mu(E')=n_1+n_{12}+m_1+m_{12}-1\text{ and}\\
\mu(E'\cup\{e_2\})-\mu(E')=n_2+n_{12}+m_2+m_{12}-1,
\end{split}\end{equation*}
whereas the right hand side of \eqref{eq:laci} is
\[\mu(E'\cup\{\,e_1,e_2\,\})-\mu(E')=\left\{\begin{array}{l>{\hspace{-12pt}}r}
n_1+n_2+m_1+m_2-2&\text{if }n_{12}=m_{12}=0\\
n_1+n_2+n_{12}+m_1+m_2+m_{12}-1&\text{otherwise.}
\end{array}\right.\]
The required inequality holds in all cases (with equality if and only if $n_{12}=m_{12}=0$ or one of the two is $1$ and the other $0$).
\end{proof}

\begin{Def}\label{def:PH}
The \emph{polymatroid of the hypergraph} $\mathscr H$ is the polymatroid $P_\mathscr H=P_\mu$ associated to the non-decreasing submodular set function $\mu$ of equation \eqref{eq:subm}.
\end{Def}

The next Proposition is just a reformulation of Theorem \ref{thm:politop} and Lemma \ref{lem:kormentes} which was used in its proof.

\begin{all}
The integer lattice points in the polymatroid $P_\mathscr H$ of the hypergraph $\mathscr H$ are exactly the functions $f\colon E\to\N$ so that $\bip\mathscr H$ has a cycle-free subgraph with valence $f(e)+1$ at each $e\in E$. The hypertree polytope $Q_\mathscr H$ coincides with the base polytope of $P_\mathscr H$.
\end{all}


Each graph has a so-called cycle matroid. This is a classical construction that goes back to the very beginning of matroid theory. The rank function of any matroid is submodular and therefore defines a polymatroid. (In fact, one way to define a matroid is as a polymatroid associated to a non-negative integer-valued, non-decreasing submodular set function that assigns $0$ to $\varnothing$ and $0$ or $1$ to singleton sets.) Thus every graph has a natural polymatroid associated to it, and Definition \ref{def:PH} generalizes this association to hypergraphs.




\subsection{Hyperedges in sequence}
We are going to need the following elementary observations later. 

\begin{lemma}\label{lem:erdo}
Let the spanning tree $\Gamma\subset\bip\mathscr H$ induce the hypertree $f\colon E\to\N$. The inequality \eqref{eq:hypertree2} is sharp for $E'\subset E$ if and only if $\Gamma\cap(\bip\mathscr H\big|_{E'})$ is a spanning forest in the latter. Given such a $\Gamma$ and $E'$, if we remove $\Gamma\cap(\bip\mathscr H\big|_{E'})$ from $\Gamma$ and replace it with another spanning forest $F$ of $\bip\mathscr H\big|_{E'}$, then the result $\tilde\Gamma$ is another spanning tree of $\bip\mathscr H$. It realizes a hypertree which agrees with $f$ on the set $E\setminus E'$.
\end{lemma}

\begin{proof}
The first assertion is trivial from the proof of Theorem \ref{thm:politop}. For the second, notice that paths in $\Gamma$ and $\tilde\Gamma$ are in a one-to-one correspondence, as follows. Any path $\varphi\subset\Gamma$ has maximal segments that fall within $\bip\mathscr H\big|_{E'}$ and obviously each such segment is in one connected component of $\bip\mathscr H\big|_{E'}$. Hence, we get a new path $\tilde\varphi\subset\tilde\Gamma$ by replacing each segment with the unique connection that exists between its endpoints in $F$. This correspondence of paths has an inverse constructed in the analogous way. As $\varphi$ and $\tilde\varphi$ share the same endpoints, we see that $\tilde\Gamma$ is connected and cycle-free. The final claim of the lemma is of course just stating the obvious.
\end{proof}

Now suppose that the set $E$ of hyperedges in $\mathscr H=(E,V)$ has been ordered and labeled so that
\[e_1<e_2<\cdots<e_{|E|}.\] 
Think of $\mathscr H$ (and of $\bip\mathscr H$) as being built step-by-step by adding one hyperedge at a time in the prescribed order. We will use the notation
\[G_k=\bip\mathscr H\big|_{\{e_1,\ldots,e_k\}}\]
for the bipartite graphs of the intermediate stages. In each step, the nullity of the graph may go up. We record this by introducing the \emph{nullity-jump function} of the chosen order,
\[nj(e_k)=n(G_k)-n(G_{k-1}),\]
where $k=1,2,\ldots,|E|$. Obviously, $nj(e)\ge0$ for all $e$ and $\sum_{e\in E}nj(e)=n(\bip\mathscr H)$.

\begin{lemma}\label{lem:moho}
For any order of the hyperedges, the vector $g(e)=|e|-1-nj(e)$ is a hypertree. It is the unique hypertree that makes the inequality \eqref{eq:hypertree2} true with an equality sign for all subsets $\{\,e_1,\ldots,e_k\,\}$, $k=1,2,\ldots,|E|$.
\end{lemma}

\begin{proof}
Using the order, it is easy to construct a spanning tree that realizes $g$. The $k$'th stage of the construction will be a spanning forest $F_k$ of the bipartite graph $G_k$. All edges of $\bip\mathscr H$ that are adjacent to $e_1$ will be part of $F_1$ (note that $nj(e_1)=0$). Suppose that a forest $F_{k-1}\subset G_{k-1}$, with the valences $g(e_i)+1$ at $e_1,\ldots,e_{k-1}$, respectively, has been defined. 

The graph $G_k-\{e_k\}$ consists of $G_{k-1}$ and some isolated points. Suppose that $e_k$ is joined by an edge to $c$ of its connected components. It is easy to see that $nj(e_k)=|e_k|-c$, since the `first edge' of $G_k$ to connect $e_k$ to one of the $c$ components does not increase the nullity, whereas all others after the first increase it by $1$. We define $F_k$ by adding to $F_{k-1}$ a collection of $c=|e_k|-nj(e_k)=g(e_k)+1$ edges, all adjacent to $e_k$ and leading to different components of $G_k-\{e_k\}$. This is a spanning forest of $G_k$. After the last stage, $F_{|E|}$ is a spanning tree of $G_{|E|}=\bip\mathscr H$ which realizes $g$.

From the fact that $F_k$ is a spanning forest of $G_k$ and the first claim in Lemma \ref{lem:erdo} it is immediate that \eqref{eq:hypertree2} is sharp for $g$ and each set $E'=\{\,e_1,\ldots,e_k\,\}$. No other hypertree can have that property for the simple reason that
\[g(e_k)=\mu(\{\,e_1,\ldots,e_k\,\})-\mu(\{\,e_1,\ldots,e_{k-1}\,\})\]
is uniquely determined as the difference of two consecutive right-hand sides.
\end{proof}

\begin{megj}
The greedy algorithm offers another way of defining the hypertree $g$ of the previous lemma: if for all $k=1,\ldots,|E|$, we choose $f(e_k)$ to be the maximal value so that the valences $f(e_i)+1$, $i=1,\ldots,k$, can be realized by a cycle-free subgraph of $G_k$, then it can be shown that $f=g$.
\end{megj}

We will also refer to the hypertree of Lemma \ref{lem:moho} as the \emph{greedy hypertree} of the given order. The next definition will be fundamental in the next section.

\begin{Def}\label{def:transfer}
Let $\mathscr H=(V,E)$ be a hypergraph, $f\colon E\to\N$ a hypertree and $a,b\in E$ hyperedges. We say that \emph{$f$ is such that a transfer of valence is possible} from $a$ to $b$ if the function obtained from $f$ by reducing $f(a)$ by $1$ and increasing $f(b)$ by $1$ is also a hypertree.
\end{Def}


\begin{lemma}\label{cor:nemeles}
For any non-empty collection $E'\subset E$ of hyperedges, there exists a hypertree $f\colon E\to\N$ so that \eqref{eq:hypertree2} is sharp for $f$ and $E'$. If $E'$ and $f$ are so that \eqref{eq:hypertree2} is not sharp, then $f$ is such that it is possible to transfer valence from some element of $E\setminus E'$ to some element of $E'$.
\end{lemma}

\begin{proof}
Choose any order in which the elements of $E'$ are the smallest and construct its greedy hypertree $g\colon E\to\N$. By Lemma \ref{lem:moho}, $g$ has the required property.

To prove the second claim, we construct an indirect argument using Propositions \ref{pro:submodular} and \ref{pro:szoros}. Let $f$ be a hypertree and $E'\subset E$ a subset so that no transfer of valence is possible from any element of $E\setminus E'$ to any element of $E'$. By Theorem \ref{thm:politop}, this implies that for any $a\in E'$ and $b\in E\setminus E'$, there exists a subset $U_{a,b}\subset E$ of hyperedges so that with the set function $\mu$ of \eqref{eq:subm}, we have
\[a\in U_{a,b},\quad b\not\in U_{a,b},\quad\text{and}\quad\sum_{x\in U_{a,b}}f(x)=\mu(U_{a,b}).\]
In other words, $U_{a,b}$ is tight at $f$. 
Then so is $E'=\bigcup_{a\in E'}\left(\bigcap_{b\in E\setminus E'}U_{a,b}\right)$, i.e., \eqref{eq:hypertree2} is sharp for $E'$ and $f$.
\end{proof}

\begin{lemma}\label{lem:mohofelett}
Let $\mathscr H=(V,E)$ be a hypergraph with its hyperedges ordered as above. Let $f\colon E\to\N$ be a hypertree and $e\in E$ a hyperedge so that $f(e)>g(e)$, where $g$ is the greedy hypertree of the order. Then there exists a hyperedge $x<e$ so that $f$ is such that a transfer of valence is possible from $e$ to $x$.
\end{lemma}

\begin{proof}
Assume the contrary. Then for all $x<e$, there exists a set $U_{e,x}$ of hyperedges that is tight at $f$, contains $x$, and does not contain $e$. Let $U=\bigcup_{x<e}U_{e,x}$. Then we have $e\not\in U$ and $\{\,x\mid x<e\,\}\subset U$ while, by Proposition \ref{pro:szoros}, $U$ is also tight at $f$.

Now define a new order $<'$ on $E$ by letting the elements of $U$ be smallest while keeping the original order among them. (It does not matter how the rest of $E$ gets ordered.) Let $g'$ be the greedy hypertree of $<'$. By Lemma \ref{lem:moho}, it is such that both $\{\,x\mid x<e\,\}$ and $U$ are tight at $g'$. Fix realizations $\Gamma'$ of $g'$ and $\Gamma$ of $f$. By Lemma \ref{lem:erdo}, these are such that both $\Gamma\cap(\bip\mathscr H\big|_U)$ and $\Gamma'\cap(\bip\mathscr H\big|_U)$ are spanning forests in $\bip\mathscr H\big|_U$. If in $\Gamma$ we replace the former with the latter, the result, again by Lemma \ref{lem:erdo}, is a spanning tree in $\bip\mathscr H$. Now for the hypertree $f'$ realized by this spanning tree we have
\begin{multline*}
\sum_{x\le e}f'(x)=f'(e)+\sum_{x<e}f'(x)=f(e)+\sum_{x<e}g'(x)=f(e)+\sum_{x<e}g(x)\\
>g(e)+\sum_{x<e}g(x)=\sum_{x\le e}g(x)=\mu(\{\,x\mid x\le e\,\}),
\end{multline*}
which contradicts Theorem \ref{thm:politop}.
\end{proof}

\section{Internal and external polynomials}\label{sec:poly}

\subsection{Activities}
Let $\mathscr H=(V,E)$ be a hypergraph so that $\bip\mathscr H$ is connected. Just as in Subsection \ref{oldtutte}, we order 
the set $E$ arbitrarily. With regard to a fixed hypertree $f$, we make the following definitions. 

\begin{Def}\label{def:activity}
A hyperedge $e\in E$ is \emph{internally active} with respect to the hypertree $f$ if it is not possible to decrease $f(e)$ by $1$ and increase $f$ of a hyperedge smaller than $e$ by $1$ so that another hypertree results. 

A hyperedge $e\in E$ is \emph{externally active} with respect to $f$ if it is not possible to increase $f(e)$ by $1$ and to decrease $f$ of a smaller hyperedge by $1$ so that another hypertree results. 
\end{Def}

Recall that in Definition \ref{def:transfer} the operation on hypertrees used above was termed a transfer of valence.  For example, a hyperedge is internally active with respect to a hypertree if it cannot transfer valence to smaller hyperedges, and it is externally active if valence cannot be transferred to it from smaller hyperedges. Regarding transfers of valence, we will need the following lemma later.

\begin{lemma}\label{lem:rombusz}
Let $a$, $b$, and $c$ be hyperedges in $\mathscr H=(V,E)$ and $f$ such a hypertree that $a$ can transfer valence to $b$ and $b$ can transfer valence to $c$. Then $f$ is also such that $a$ can transfer valence to $c$. In other words, regarding the rhombus of Figure \ref{fig:rombusz} (explained below), if the three lattice points indicated by full dots are hypertrees then so is the one represented by the hollow dot.
\end{lemma}

\begin{proof}
This is immediate from Lemma \ref{lem:discreteconvex} but we think it worthwile to examine the picture. In $\R^E$, the three-dimensional affine subspace through $f$ spanned by $\{\,u_a,u_{b},u_{c}\}$ forms a two-dimensional intersection $Q_0$ with $Q_\mathscr H$. The triangle we labelled with $abc$ in Figure \ref{fig:rombusz} is in fact spanned by the vectors whose $u_a$, $u_{b}$, and $u_{c}$-components are equal to $(f(a)+f(b)+f(c),0,0)$, $(0,f(a)+f(b)+f(c),0)$, and $(0,0,f(a)+f(b)+f(c))$, respectively, whereas their other components are identical to those of $f$.

\begin{figure}[htbp]
\unitlength .8pt
\begin{picture}(180,160)
\put(10,10){\line(1,0){160}}
\put(10,10){\line(3,5){80}}
\put(170,10){\line(-3,5){80}}
\thicklines
\put(80,30){\circle*{6}}
\put(62,60){\circle*{6}}
\put(116,30){\circle*{6}}
\put(98,60){\circle{6}}
\put(80,30){\line(1,0){36}}
\put(80,30){\line(-3,5){18}}
\put(80,30){\line(3,5){16}}
\put(2,5){$a$}
\put(173,5){$b$}
\put(87,148){$c$}
\put(66,26){$f$}
\put(121,26){$f_1$}
\put(46,60){$f_2$}
\put(103,60){$\hat f$}
\end{picture}
\caption{A rhombus of hypertrees.}
\label{fig:rombusz}
\end{figure}
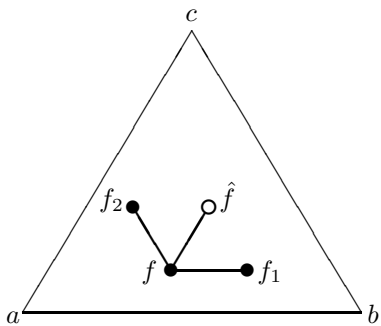

The conditions in the lemma are that $f$, $f_1$, and $f_2$ are in $Q_\mathscr H$ and the conclusion is that $\hat f$ is, too. Now this is obvious from the fact (Theorem \ref{thm:politop}) that $Q_0$ is cut out from the plane by lines parallel to the sides of the triangle.
\end{proof}

It is easy to see that  the notions of Definition \ref{def:activity} generalize the ones found in Definition \ref{def:oldactive}. It will be slightly more convenient for us to work with their negations in what follows.

\begin{Def}
Let $\mathscr H=(V,E)$ be a hypergraph so that $\bip\mathscr H$ is connected. With respect to a given hypertree $f\colon E\to\N$ and some order of the hyperedges, let the number of internally inactive hyperedges of $\mathscr H$ be denoted with $\bar\iota(f)$ and the number of externally inactive hyperedges be denoted with $\bar\varepsilon(f)$. Then, let the \emph{internal polynomial} and the \emph{external polynomial} of $\mathscr H=(V,E)$ be defined as
\begin{equation}\label{polinomok}
I_\mathscr H(\xi)=\sum_f \xi^{\bar\iota(f)}\quad\text{and}\quad X_\mathscr H(\eta)=\sum_f \eta^{\bar\varepsilon(f)},
\end{equation}
respectively, where both summations are over all hypertrees $f$ in $\mathscr H$. 
\end{Def}

These notions generalize the valuations 
\begin{equation}\label{eq:regitutte}
\xi^{|V|-1}T_G(1/\xi,1)\quad\text{and}\quad\eta^{|E|-|V|+1}T_G(1,1/\eta),
\end{equation} 
respectively, of the classical Tutte polynomial $T_G$ of the graph $G=(V,E)$. 

There are two ways of treating hypergraphs whose associated bipartite graphs are disconnected. One is to assign to them the product of the polynomials associated to their connected components. The other is to extend the definition verbatim, which results in the value $0$ for both polynomials because a disconnected graph has no spanning trees and therefore the set of hypertrees is empty. In this paper we take the latter approach. This will hardly matter since we almost always assume $\bip\mathscr H$ to be connected.

Our first order of business, however, is to address well definedness.

\begin{tetel}\label{thm:independent}
The formulas \eqref{polinomok} for the internal and external polynomials do not depend on the chosen order of the hyperedges.
\end{tetel}

We will need to following statement.



\begin{lemma}\label{lem:nyil}
Let $e_1\ne e_2$ be hyperedges in the hypergraph $\mathscr H=(V,E)$ and let $E_0=E\setminus\{\,e_1,e_2\,\}$. Fix a function $f_0\colon E_0\to\N$ and consider its extensions to $E$ that are hypertrees in $\mathscr H$. Among these, let $f_1$ and $f_2$ be such that $f_1(e_1)>f_2(e_1)$ and let also $x\in E_0$ be a hyperedge.
\begin{enumerate}
\item If $f_1$ is such that valence can be transferred from $e_1$ to $x$ then $f_2$ is such that valence can be transferred from $e_2$ to $x$.
\item If $f_1$ is such that valence can be transferred from $x$ to $e_2$ then $f_2$ is such that valence can be transferred from $x$ to $e_1$.
\end{enumerate}
In other words, regarding both `staples' of Figure \ref{fig:nyil}, if the three lattice points indicated by full dots are hypertrees then so is the one represented by the hollow dot.
\end{lemma}

\begin{proof}
The argument is very similar to the proof of Lemma \ref{lem:rombusz}. Examining Figure \ref{fig:nyil}, it is easy to see that if $f_1$, $f_2$, and $g_1$ satisfy all constraints of Theorem \ref{thm:politop}, then so does $\hat g_1$. Likewise, if $f_1$, $f_2$, and $g_2$ are in $Q_\mathscr H$ then so is $\hat g_2$.
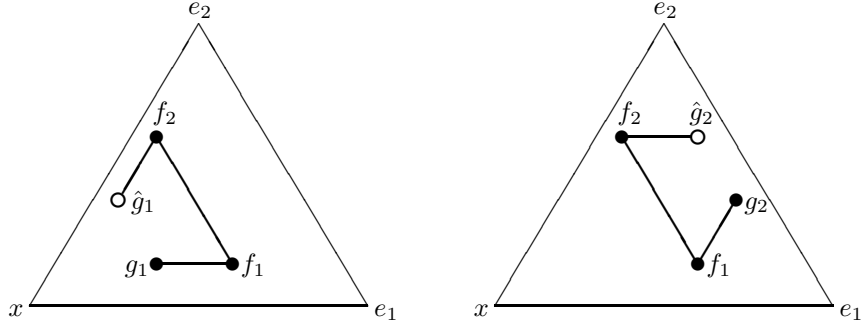
\begin{figure}[htbp]
\unitlength .8pt
\begin{picture}(393,160)
\put(10,10){\line(1,0){160}}
\put(10,10){\line(3,5){80}}
\put(170,10){\line(-3,5){80}}
\put(0,5){$x$}
\put(173,5){$e_1$}
\put(85,148){$e_2$}
\put(110,26){$f_1$}
\put(68,98){$f_2$}
\put(55,26){$g_1$}
\put(58,56){$\hat g_1$}
\put(230,10){\line(1,0){160}}
\put(230,10){\line(3,5){80}}
\put(390,10){\line(-3,5){80}}
\put(220,5){$x$}
\put(393,5){$e_1$}
\put(305,148){$e_2$}
\put(330,26){$f_1$}
\put(288,98){$f_2$}
\put(348,56){$g_2$}
\put(322,98){$\hat g_2$}
\thicklines
\put(70,90){\line(3,-5){36}}
\put(70,90){\line(-3,-5){16}}
\put(106,30){\line(-1,0){36}}
\put(70,90){\circle*{6}}
\put(106,30){\circle*{6}}
\put(70,30){\circle*{6}}
\put(52,60){\circle{6}}
\put(290,90){\line(3,-5){36}}
\put(290,90){\line(1,0){33}}
\put(326,30){\line(3,5){18}}
\put(290,90){\circle*{6}}
\put(326,30){\circle*{6}}
\put(344,60){\circle*{6}}
\put(326,90){\circle{6}}
\end{picture}
\caption{Two `staples' formed by hypertrees.}
\label{fig:nyil}
\end{figure}
\end{proof}

\begin{proof}[Proof of Theorem \ref{thm:independent}]
We will emulate Tutte's original proof, i.e., analyze the effect of changing the relative position in the order of two adjacent hyperedges. Let $a,b\in E$ and assume that $o_1$ is an order on $E$ so that $a<_1b$ with no other hyperedge in between, whereas the order $o_2$ only differs from $o_1$ in that $b<_2a$. The rest of the hyperedges are partitioned in two sets so that in both orders, elements of $E_-$ are smaller than $a$ and $b$ while elements of $E_+$ are larger than both.

Let $f\colon E\to\N$ be a hypertree. Our goal is to compare the values $\bar\iota_1(f)$ and $\bar\iota_2(f)$ (as well as $\bar\varepsilon_1(f)$ and $\bar\varepsilon_2(f)$) resulting from the two orders. If a hyperedge differs from both $a$ and $b$, then it is easy to check that (in the internal as well as in the external sense) it is active with respect to $f$ in $o_1$ if and only if the same holds in $o_2$. Let us now separate three cases.

\begin{enumerate}[I.]
\item\label{case1} If $f$ is such that no transfer of valence is possible between $a$ and $b$, then the activity statuses of $a$ and $b$ are also unaffected by the change in the order. Hence in such cases $\bar\iota_1(f)=\bar\iota_2(f)$ and $\bar\varepsilon_1(f)=\bar\varepsilon_2(f)$.

\item\label{case2} Next, assume that $f$ is such that valence can be transferred between $a$ and $b$ in both directions. Then in the order $o_1$ the hyperedge $b$ is not active whereas with respect to $o_2$, the hyperedge $a$ is not active. (This holds in both the internal and external senses.) Now according to Lemma \ref{lem:rombusz}, if $a$ is not active in $o_1$ (i.e., there is some $x\in E_-$ so that valence can be transferred from $a$ to $x$ (internal case) or from $x$ to $a$ (external case)), then $b$ is not active in $o_2$ and vice versa. Therefore we again have $\bar\iota_1(f)=\bar\iota_2(f)$ and $\bar\varepsilon_1(f)=\bar\varepsilon_2(f)$.

\item\label{case3} Lastly, let $f$ be such that valence can be transferred between $a$ and $b$ but only in one direction. Before analyzing activities, let us establish that these kinds of hypertrees are partitioned into pairs. Indeed, the line 
\[\left\{\,g\colon E\to\R\biggm| g\big|_{E\setminus\{a,b\}}=f\big|_{E\setminus\{a,b\}}, \sum_{e\in E}g(e)=|V|-1\,\right\}\subset\R^E\]
intersects the polytope $Q_\mathscr H$ in a segment\footnote{The endpoints of this segment are $f$ and another lattice point $f^*$. This can be shown using a property of polymatroids called total dual integrality. However we will not need this presently.}. Among the lattice points of the segment, $f$ represents one extreme; let the other extreme be denoted with $f^*$. This is necessarily different from $f$ by our assumption and if a transfer of valence is possible only from $a$ to $b$ (and not in the other direction) at one of $f$ and $f^*$, then the opposite transfer (and only that) is possible at the other one. It is also clear that $(f^*)^*=f$, so we get the desired pairs. (The lattice points between $f$ and $f^*$ were discussed in case \ref{case2}, whereas case \ref{case1} is when the segment degenerates to a point.)

Without loss of generality we may assume that $f$ is such that $a$ can transfer valence to $b$ and the opposite transfer is possible at $f^*$. Let us first examine the hyperedges $a$ and $b$ themselves. We will start with their internal activities.

\begin{enumerate}
\item Assume that $\bar\iota_1(f)=\bar\iota_2(f)$. We will show that in this case $\bar\iota_1(f^*)=\bar\iota_2(f^*)$ holds, too. As the activity status of $b$ with respect to $f$ is the same in $o_1$ as in $o_2$, a similar property has to hold for $a$. Because $a$ is not active with respect to $f$ in $o_2$, this is only possible if there exists a hyperedge $x\in E_-$ so that $f$ is such that valence can be transferred from $a$ to $x$. By Lemma \ref{lem:nyil}, this implies that $b$ is inactive with respect to $f^*$ in both $o_1$ and $o_2$. Whether $a$ is active with respect to $f^*$ depends on the same thing in both orders: namely, on whether $f^*$ is such that valence can be transferred from $a$ to some element $y\in E_-$.
\item\label{bbb} The only way $\bar\iota_1(f)$ and $\bar\iota_2(f)$ can be different is if $f$ is such that $a$ cannot transfer valence to any element of $E_-$. In such cases $b$ has the same property by Lemma \ref{lem:rombusz} so that $\bar\iota_1(f)=\bar\iota_2(f)-1$ ($b$ is active with respect to $f$ in both orders and $a$ is active only in $o_1$). Examining $f^*$ now, by Lemma \ref{lem:nyil} we see that it has to be such that $b$  cannot transfer valence to any element of $E_-$. This implies, using Lemma \ref{lem:rombusz}, that $a$ has the same property. Therefore we have $\bar\iota_1(f^*)=\bar\iota_2(f^*)+1$ as $a$ is active with respect to $f^*$ in both orders and $b$ is active in $o_2$ only. 
\end{enumerate}
External activities can be handled in almost the same way.
\begin{enumerate}[(a')]
\item If $\bar\varepsilon_1(f)=\bar\varepsilon_2(f)$, then (since the activity status of $a$ with respect to $f$ is independent of order) there must be some $x\in E_-$ so that valence can be transferred from $x$ to $b$. Then by Lemma \ref{lem:nyil}, $f^*$ is such that valence can be transferred from $x$ to $a$, making $a$ inactive with respect to $f^*$ in both orders. Since the activity status of $b$ with respect to $f^*$ is the same in both orders, we obtain $\bar\varepsilon_1(f^*)=\bar\varepsilon_2(f^*)$.
\item\label{bbbb} If $\bar\varepsilon_1(f)\ne\bar\varepsilon_2(f)$, that is if $f$ is such that no $x\in E_-$ can transfer valence to $b$, then the same holds for transfers of valence from $x$ to $a$ by Lemma \ref{lem:rombusz}. Therefore $a$ is active with respect to $f$ in both orders, whereas $b$ is active in $o_2$ and inactive in $o_1$. By the preceding paragraph, a similar analysis has to apply to $f^*$, i.e., $b$ is active with respect to $f^*$ in both orders, whereas $a$ is active in $o_1$ and inactive in $o_2$. So we have obtained $\bar\varepsilon_1(f)=\bar\varepsilon_2(f)+1$ and $\bar\varepsilon_1(f^*)=\bar\varepsilon_2(f^*)-1$. 
\end{enumerate}
\end{enumerate}

The only hypertrees whose $\bar\iota$ or $\bar\varepsilon$ values actually changed as a result of switching from $o_1$ to $o_2$ were described in the cases (b) and (b') above. We saw that they occur in pairs of the form $\{f,f^*\}$. Now to complete our proof it suffices to show that for such a pair and any hyperedge $y$ different from $a$ and $b$, the activity status 
of $y$ with respect to $f$ is the same as with respect to $f^*$. (We saw that the choice between $o_1$ and $o_2$ does not matter for $y$, only for $a$ and $b$. Switching to a new hypertree, on the other hand, could in principle make a big difference.) Indeed, that will imply $\bar\iota_1(f)=\bar\iota_2(f^*)$ and $\bar\iota_2(f)=\bar\iota_1(f^*)$ in case (b), as well as $\bar\varepsilon_1(f)=\bar\varepsilon_2(f^*)$ and $\bar\varepsilon_2(f)=\bar\varepsilon_1(f^*)$ in case (b') so that, regardless of order, the generating functions $I_\mathscr H$ and $X_\mathscr H$ always encode the same information.

As we are about to relate activities with respect to different hypertrees, this last part of the proof is where we will rely most heavily on our assumptions. What we need is essentially Lemma \ref{lem:teglalap} but we chose to spell the argument out without an explicit reference to it. It will be convenient to branch out into several cases again.

\begin{enumerate}[(1)]
\item We will deal with internal activities first. Recall that we are under the assumptions of the case \ref{case3}(b), in particular $f$ and $f^*$ are such that neither $a$ nor $b$ can transfer valence to any element of $E_-$. Let $y$ be a hyperedge which is internally inactive with respect to $f$. Our goal is to show that $y$ is also internally inactive with respect to $f^*$. 
\begin{enumerate}[i.]
\item Assume $y\in E_+$. If $f$ is such that $y$ can transfer valence to $a$ or $b$, then Lemma \ref{lem:rombusz} implies that it can definitely transfer to $b$ and then Lemma \ref{lem:nyil} says that $f^*$ is such that $y$ can transfer valence to $a$. Therefore $y$ is inactive with respect to $f^*$.

If $y\in E_+$ but $f$ is such that $y$ cannot transfer valence to $a$ or $b$, then $f^*$ has to have the same property by the usual combination of Lemmas \ref{lem:rombusz} and \ref{lem:nyil}. Because $y$ is inactive with respect to $f$, the hypertree $f$ must be such that $y$ can transfer valence to some $x<y$ with $a\ne x\ne b$. We will show that $f^*$ is also such that $y$ can transfer valence to $x$. Suppose the opposite is true. Then there are sets of hyperedges $U_a\ni a$, $U_b\ni b$, and $U_x\ni x$ so that all three are tight at $f^*$ and neither contains $y$. (Here tightness is in the sense of Definition \ref{def:tight} with regards the set function $\mu$ of Proposition \ref{pro:submodular}.) By Proposition \ref{pro:szoros}, the union of the three subsets is also tight at $f^*$. Since the union contains both $a$ and $b$, the sum of the $f^*$-values over it agrees with the sum of the $f$-values. Therefore $U_a\cup U_b\cup U_x$ is also tight at $f$, which contradicts our assumption that $f$ is such that $y$ (which is not in the union) can transfer valence to $x$ (which is).
\item Let now $y\in E_-$. Then of course there is a hyperedge $x\in E_-$ so that $f$ is such that $y$ can transfer valence to $x$. We will show that $f^*$ is also such that $y$ can transfer valence to $x$. Assuming the contrary, we find sets of hyperedges $U_y\not\ni y$, $U_a\not\ni a$, and $U_b\not\ni b$ which are tight at $f^*$ and all contain $x$. Their intersection is also tight at $f^*$ by Proposition \ref{pro:szoros}. As $U_a\cap U_b\cap U_y$ contains neither $a$ nor $b$, it is also tight at $f$. But that is a contradiction with the transfer of valence from $y\not\in U_a\cap U_b\cap U_y$ to $x\in U_a\cap U_b\cap U_y$ which is possible at $f$.
\end{enumerate}
\item In the external case, the assumptions of \ref{case3}(b') were that $f$ and $f^*$ are both such that no transfer of valence is possible from elements of $E_-$ to $a$ or to $b$. In such a situation, let $y$ be externally inactive with respect to $f$. We wish to prove that $y$ is also externally inactive with respect to $f^*$. 
\begin{enumerate}[i'.]
\item If $y\in E_+$, then it may be that $f$ is such that $a$ or $b$ can transfer valence to $y$. If $b$ can, then so can $a$ by Lemma \ref{lem:rombusz}. Thus in either case we can apply Lemma \ref{lem:nyil} to conclude that $f^*$ is such that $b$ can transfer valence to $y$.

Suppose now that $y\in E_+$ but $f$ (and consequently $f^*$) is such that neither $a$ nor $b$ can transfer valence to $y$. Then some other hyperedge $x<y$ can and we will prove that the same transfer is also possible at $f^*$. If not, then there are sets of hyperedges $U_a\not\ni a$, $U_b\not\ni b$, and $U_x\not\ni x$, all containing $y$, that are tight at $f^*$. Then their intersection has the same property by Proposition \ref{pro:szoros} and furthermore, as it contains neither $a$ nor $b$, it is tight at $f$ as well. That contradicts the assumption that at $f$, a transfer of valence is possible from outside of the set (from $x$) to one of its elements (namely $y$).
\item Assuming $y\in E_-$, it being externally inactive implies the existence of another hyperedge $x\in E_-$ so that $f$ is such that $x$ can transfer valence to $y$. If the same transfer is not possible at $f^*$, that means that there are sets of hyperedges $U_y\ni y$, $U_a\ni a$, and $U_b\ni b$ that are tight at $f^*$ and none of them contains $x$. Then their union has the same properties and it is also tight at $f$, which yields the usual contradiction.
\end{enumerate}
\end{enumerate}
This completes the proof.
\end{proof}

\begin{pelda}\label{ex:harom}
We revisit the hypertrees that were discussed in Example \ref{ex:ketto}.
They are shown in Figure \ref{fig:ujfak} with an actual spanning tree realization for each. For the particular orders $a<b<c$ and $p<q<r<s$, respectively, we also indicated the internal and external activity status of each hyperedge with respect to each hypertree, as well as the resulting $(\bar\iota,\bar\varepsilon)$ values. Thus, we find that 
\begin{equation}\label{eq:egyik}
I_{\mathscr G_0}(\xi)=1+3\xi+3\xi^2\quad\text{and}\quad X_{\mathscr G_0}(\eta)=1+3\eta+3\eta^2,
\end{equation} 
whereas 
\begin{equation}\label{eq:masik}
I_{\mathscr G_1}(\xi)=1+3\xi+3\xi^2\quad\text{and}\quad X_{\mathscr G_1}(\eta)=1+2\eta+3\eta^2+\eta^3.
\end{equation}
These values are in line with Propositions \ref{pro:const} and \ref{pro:rang} below, as well as with Conjecture \ref{conj:dual}. 

Up to isomorphism there is only one order on the three-element set $V_0$. The four-element set $V_1$, on the other hand, has four essentially different orders depending on the position of $q$. The reader may check that the other three give rise to the same interior and exterior polynomials as in \eqref{eq:masik}.

Note the lack of similarity between our polynomials and the one shown in Example \ref{ex:egy}. Indeed, the author is not aware of any formula relating the Tutte polynomial $T_G(x,y)$ of a bipartite graph $G$ and the interior and exterior polynomials of its induced hypergraphs.
\begin{figure}[htbp]
\labellist
\pinlabel $(1,1)$ at 170 630
\pinlabel $(2,0)$ at 585 630
\pinlabel $(1,2)$ at 1000 630
\pinlabel $(2,2)$ at 1415 630
\pinlabel $(2,1)$ at 1830 630
\pinlabel $(0,2)$ at 2245 630
\pinlabel $(1,1)$ at 2660 630
\pinlabel $(1,2)$ at 170 80
\pinlabel $(0,3)$ at 585 80
\pinlabel $(2,1)$ at 1000 80
\pinlabel $(2,1)$ at 1415 80
\pinlabel $(1,2)$ at 1830 80
\pinlabel $(1,2)$ at 2245 80
\pinlabel $(2,0)$ at 2660 80
\endlabellist
   \centering
   \includegraphics[width=\linewidth]{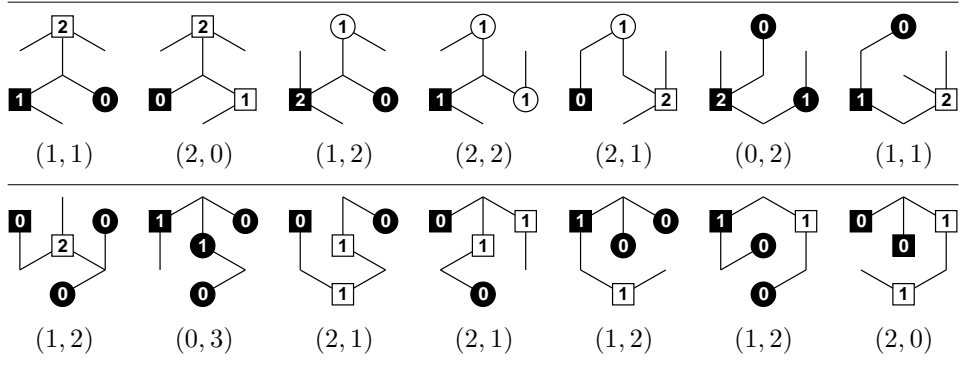} 
   \caption{Hypertrees in a pair of abstract dual hypergraphs. Hollow shapes are internally inactive and full ones are internally active; circles are externally inactive while squares are externally active.}
   \label{fig:ujfak}
\end{figure}
\end{pelda}

\begin{megj}
Unfortunately, for a hypergraph that is not a graph, the two-variable polynomial $\sum_f \xi^{\bar\iota(f)}\eta^{\bar\varepsilon(f)}$ (cf.\ Definition \ref{def:tuttepoly}) does depend on the choice of order. For instance, for the hypergraph $\mathscr G_1$ of the previous example, the order that we used there gives rise to the two-variable generating function $\eta^3+3\xi\eta^2+2\xi^2\eta+\xi^2$ (this is directly read off of the bottom row of Figure \ref{fig:ujfak}), whereas the order $p<r<s<q$ gives $\xi\eta^3+\xi^2\eta^2+\xi\eta^2+\eta^2+2\xi^2\eta+\xi$. (Note how the two polynomials do coincide after setting $\xi=1$ or $\eta=1$.)

\end{megj}


\subsection{Extension to polymatroids}\label{ssec:polypoly}
The interior and exterior polynomials are constructed from the hypertree polytope. The same procedure can be carried out after replacing $Q_\mathscr H$ with the set of bases $B_\mu$ of an integer extended polymatroid. Because of Lemma \ref{lem:nincstobbbazis} and the fact that integer polymatroids can always be defined by integer-valued non-decreasing submodular set functions, we immediately obtain invariants of integer polymatroids as well. In this subsection we outline the minimal modifications that are needed in our arguments to get the generalization.

Let $S$ be a finite set and $\mu\colon\mathscr P(S)\to\R$ an integer-valued submodular set function, with associated base polytope $B_\mu\subset EP_\mu$. We say that the base $\mathbf x\in B_\mu\cap\Z^S$ is such that a \emph{transfer} is possible from $s_1\in S$ to $s_2\in S$ if by decreasing the $s_1$-component of $\mathbf x$ by $1$ and increasing its $s_2$-component by $1$, we get another base.

Now order $S$ arbitrarily. Call an element $s\in S$ \emph{internally active} with respect to the base $\mathbf x\in B_\mu\cap\Z^S$ if $\mathbf x$ is such that a transfer is possible from $s$ to a smaller element of $S$. We say that $s$ is \emph{externally active} with respect to $\mathbf x$ if it is such that a transfer is possible to $s$ from a smaller element of $S$. 

For $\mathbf x\in B_\mu\cap\Z^S$, define $\bar\iota(\mathbf x)$ to be the number of elements of $S$ that are not internally active with respect to $\mathbf x$. Similarly, let $\bar\varepsilon(\mathbf x)$ denote the number of elements of $S$ that are not externally active with respect to $\mathbf x$. 

\begin{Def}\label{def:polym}
Let
\[I_\mu(\xi)=\sum_{\mathbf x\in B_\mu\cap\Z^S}\xi^{\bar\iota(\mathbf x)}\quad\text{and}\quad X_\mu(\eta)=\sum_{\mathbf x\in B_\mu\cap\Z^S}\eta^{\bar\varepsilon(\mathbf x)}\]
and call these quantities the \emph{interior polynomial}, respectively the \emph{exterior polynomial} of the submodular set function $\mu$.
\end{Def}

The fact that these polynomials do not depend on the order that was used to write them can be shown just like in the proof of Theorem \ref{thm:independent}. Indeed, the first half of that proof depended on the elementary lemmas \ref{lem:rombusz} (rhombus lemma) and \ref{lem:nyil} (staple lemma) where we only used the fact that $Q_\mathscr H$ is cut out by placing an upper bound on each partial sum of the components of its elements. We also used the triviality that lines intersect bounded convex sets in segments. Then in the second half of the argument we relied on Proposition \ref{pro:szoros} and the fact that our upper bounds are values of a tight set function.

\begin{figure}[htbp]
\labellist
\small
\pinlabel $\mathbf a$ at 610 660
\pinlabel $\mathbf b$ at -30 10
\pinlabel $\mathbf c$ at 1170 320
\pinlabel $\mathbf d$ at 450 980
\endlabellist
   \centering
   \includegraphics[width=2.5in]{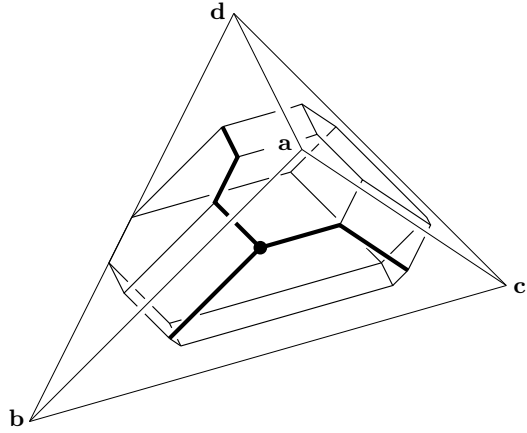} 
   \caption{The geometry of the interior polynomial.}
   \label{fig:szakoca}
\end{figure}

To illustrate the geometry underlying our polynomials, we show in Figure \ref{fig:szakoca} the hypothetical base polytope of a polymatroid $P_\mu$ with a ground set $S=\{\,a,b,c,d\,\}$ of four elements. (We reused the diagram from the proof of Lemma \ref{lem:teglalap} but placed it in a different frame of reference.) It is situated in a tetrahedron with vertices $\mathbf a(\mu(S),0,0,0)$, $\mathbf b(0,\mu(S),0,0)$, $\mathbf c(0,0,\mu(S),0)$, and $\mathbf d(0,0,0,\mu(S))$.  If we set the order $a<b<c<d$, then the marked vertex represents the only hypertree with $\bar\iota=0$; the other lattice points along the six thickened edges are the ones with $\bar\iota=1$; the rest of the lattice points that are visible (along a total of seven faces) in the view we used are those with $\bar\iota=2$; finally the invisible lattice points have $\bar\iota=3$.

A similar picture applies to the exterior polynomial. It would be very interesting to see whether invariants of arbitrary (i.e., not necessarily integer) extended polymatroids can be defined in a way analogous to Definition \ref{def:polym}, but replacing counts of lattice points along certain faces by taking the volumes of those faces.

It may also be worth it to investigate if there is a wider class of set functions (or polytopes) for which interior and exterior polynomials are well defined. Already in the set function case some assumption is 
necessary, as the next example shows.

\begin{pelda}
Consider the tetrahedron
\[Q=\conv\{\,(1,1,0,0),(1,0,1,0),(1,0,0,1),(0,0,1,1)\,\}\subset\R^4_{xyzt}.\]
The four given points are exactly the integer solutions of the equation $x+y+z+t=2$ and the linear inequalities
\[\begin{array}{rcl@{\hspace{.5in}}rcl@{\hspace{.5in}}rcl}
x&\le&1;&x+y&\le&2;&y+z+t&\le&2;\\
y&\le&1;&x+z&\le&2;&x+z+t&\le&2;\\
z&\le&1;&x+t&\le&2;&x+y+t&\le&2;\\
t&\le&1;&y+z&\le&1;&x+y+z&\le&2.\\
&&&y+t&\le&1;&&&\\
&&&z+t&\le&2;&&&
\end{array}\]
Furthermore, each of the fourteen (fifteen) inequalities is sharp for at least one of the four points. The right hand sides 
\emph{do not} give a tight set function. For example, the set of the $y$ and $z$ coordinates, as well as the set of the $y$ and $t$ coordinates is tight at the point $(1,1,0,0)$. However the union of the two sets is not tight at the same point. We also see that Lemma \ref{lem:teglalap} is violated: the face along which the sum of the $y$ and $z$ (or the $y$ and $t$) coordinates takes its maximum is a triangle instead of a rectangle.

From the order $x<y<z<t$, we obtain the interior and exterior polynomials $1+2\xi+\xi^2$ and $1+2\eta+\eta^2$. If we use $y<z<t<x$ instead, the interior polynomial becomes $2+2\xi^2$. For $x<t<z<y$, the exterior polynomial is $2+2\eta^2$. 
\end{pelda}

\section{Properties}\label{sec:prop}

\subsection{Low-order terms}
We first make an observation on the degrees of the interior and exterior polynomials and then turn our attention to some individual coefficients.

\begin{all}
For a hypergraph $\mathscr H=(V,E)$, the degree of its interior polynomial is at most $\min\{\,|E|,|V|\,\}-1$, while the degree of its exterior polynomial is at most $|E|-1$.
\end{all}

\begin{proof}
That $|E|-1$ is an upper bound for both degrees is obvious from the observation that the smallest hyperedge in an order is both internally and externally active with respect to any hypertree.

As to $\deg I_\mathscr H\le|V|-1$, note that in order for a hyperedge $e$ to be internally inactive with respect to a hypertree $f$, we have to have $f(e)\ge1$. This combined with \eqref{eq:hypertree3} gives the result.
\end{proof}

For the constant term in the exterior polynomial and the two lowest-order coefficients in the interior polynomial, the following two results offer a more direct way of showing their order-independence.

\begin{all}\label{pro:const}
Both the internal and external polynomials of the hypergraph $\mathscr H=(E,V)$ have $1$ for constant term.
\end{all}

\begin{proof}
Let us fix an arbitrary order on $E$ and denote the greedy hypertree of Lemma \ref{lem:moho} with $g$. We claim that $g$ is the unique hypertree with respect to which (and the given order) every hyperedge is internally active. Comparing the assertions of Lemma \ref{lem:moho} and Theorem \ref{thm:politop} with the definition of internal activity, it is clear that $g$ is indeed one such hypertree.

The fact that there are no others follows easily either from Lemma \ref{cor:nemeles} or from Lemma \ref{lem:mohofelett}. Yet it may be interesting to see an argument that does not rely on abstract properties of submodular functions.

To this end, let $f\colon E\to\N$ be a hypertree different from $g$. We need to show that there is an internally inactive hyperedge with respect to it. We will make use of notation from subsection \ref{ssec:politop}, in particular the sequence of forests $F_1\subset F_2\subset\cdots\subset F_{|E|}$ that was constructed in the proof of Lemma \ref{lem:moho}. Let $e_k$ be the smallest hyperedge in the order so that $f(e_k)\ne g(e_k)$, which of course implies $f(e_k)<g(e_k)$. By Lemma \ref{lem:erdo} we may choose a realization $\Gamma$ for $f$ so that $\Gamma\cap G_{k-1}=F_{k-1}$.

There is at least one connected component of $G_k-\{e_k\}$ which is not connected to $e_k$ by any edge in $\Gamma$, even though there is an edge $\alpha$ of $\bip\mathscr H$ between $e_k$ and this component (we can take for example the one which was selected into $F_k$). Adding $\alpha$ to $\Gamma$ creates a cycle which has to go through a hyperedge $e_l$ with $l>k$. By removing from $\Gamma$ an edge of the cycle adjacent to $e_l$ and replacing it with $\alpha$, we have created a new spanning tree for $\bip\mathscr H$. It induces a hypertree which only differs from $f$ at $e_l$ (where it is one smaller than $f$) and $e_k$ (where it is one bigger). This shows that the hyperedge $e_l$ is internally inactive with respect to $f$.

In the case of the external polynomial, the unique hypertree without an (externally) inactive hyperedge is again as in Lemma \ref{lem:moho} but constructed using the reverse order. The rest of the proof can be carried out just like above.
\end{proof}

\begin{tetel}\label{pro:rang}
For any hypergraph $\mathscr H=(E,V)$, the coefficient of the linear term in the interior polynomial $I_\mathscr H$ is the nullity (first Betti number) $n(\bip\mathscr H)$ of the bipartite graph $\bip\mathscr H$.
\end{tetel}

\begin{proof}
We will extend the analysis carried out in the proofs of Lemma \ref{lem:moho} and Proposition \ref{pro:const} a little further. After fixing an order on $E$, we are going to construct $n(\bip\mathscr H)$ hypertrees. Namely, if $e$ is a hyperedge, then we will associate to it $nj(e)$ hypertrees as follows.

Recall the greedy hypertree $g(e)=|e|-1-nj(e)$ of Lemma \ref{lem:moho} and fix one of its realizations. 
If $e$ is a hyperedge that has $nj(e)>0$ with respect to the order, then add to the realization one more edge adjacent to $e$. This creates a cycle which goes through another hyperedge $e'$. Remove an edge of this cycle adjacent to $e'$ to get a new spanning tree. Its induced hypertree is the result of a transfer of valence from $e'$ to $e$. Then add another edge and make another transfer and so on until all edges adjacent to $e$ are used up. The hyperedges that we transfer valence from are all smaller than $e$ because of Lemma \ref{lem:moho} and Theorem \ref{thm:politop}.

The argument in the previous paragraph shows that for all $1\le i\le nj(e)$, there is an $i$-element multiset of hyperedges smaller than $e$ so that the result of reducing their associated $g$-values by $1$ and increasing $g(e)$ to $g(e)+i$ is a hypertree in $\mathscr H$. Now define $f_{e,i}\colon E\to\N$ ($e\in E$, $1\le i\le nj(e)$) to be the hypertree among these so that its associated multiset $M_{e,i}$ is largest in reverse lexicographical order, i.e., if $e_1<\cdots<e_k$ are all the hyperedges smaller than $e$, then $M_{e,i}$ has the highest possible multiplicity at $e_k$; then among those the highest possible multiplicity at $e_{k-1}$, and so on.

We claim that the $f_{e,i}$ are exactly those hypertrees in $\mathscr H$ that have a unique internally inactive hyperedge in the given order. 

It is easy to see that the $f_{e,i}$ do have this property. Namely, $e$ is the unique hyperedge that is not internally active with respect to $f_{e,i}$. Indeed, if $e'>e$ then $e'$ is internally active with respect to $f_{e,i}$ because $\sum_{x<e'}f_{e,i}(x)=\sum_{x<e'}f(x)$ is already at the largest value allowed by \eqref{eq:hypertree2}. If $e'<e$, then it is internally active because a downward transfer of valence from $e'$ would contradict the way $M_{e,i}$ was chosen. Finally, $e$ itself is not active because by Lemma \ref{lem:discreteconvex} applied to $f_1=g$ and $f_2=f_{e,i}$, it can transfer valence to any element of $M_{e,i}$.

Let now $f$ be a hypertree in $\mathscr H$ so that it has a unique internally inactive hyperedge $e$. Our goal is to show that $f$ is one of the $f_{e,i}$.

First we claim that \eqref{eq:hypertree2} is sharp for $f$ and the set $E'=\{\,x\in E\mid x\le e'\,\}$ for all $e'\ge e$: indeed if it was not, then by Lemma \ref{cor:nemeles} there had to be hyperedges larger than $e'$ (and hence different from $e$) that are internally inactive with respect to $f$. Because of the similar property of $g$ stated in Lemma \ref{lem:moho}, we see that $f=g$ (and hence $f=f_{e,i}$ for any $i$) for all hyperedges larger than $e$. It also follows that $f(e)\ge g(e)$ and in fact $f(e)>g(e)$, because otherwise $e$ could not be inactive.

Next, we note that if $e'$ is a hyperedge smaller than $e$, then $f(e')\le g(e')$ because if this was not the case then Lemma \ref{lem:mohofelett} would imply that $e'$ is internally inactive. So far we have shown that $f$ is obtained from $g$ by transferring valence to $e$ from a (non-empty) multiset of hyperedges that are smaller than $e$. 

Assume now that $f\ne f_{e,i}$, where we set $i=f(e)-g(e)$. Because of the way $f_{e,i}$ was constructed, the largest hyperedge $e'$ where the two hypertrees differ is so that $e'<e$ and $f(e')>f_{e,i}(e')$. We are going to argue that $e'$ is internally inactive with respect to $f$ by showing that $f$ is such that $e'$ can transfer valence to at least one element of the set $S=\{\,x\in E\mid f(x)<f_{e,i}(x)\,\}$. (Note that $S$ is non-empty because $\sum_{y\in E}f(y)=\sum_{y\in E}f_{e,i}(y)$ and also that its elements are smaller than $e'$.)

Suppose again that the opposite is true. Then by the usual argument based on Proposition \ref{pro:szoros}, there exists a set $U$ of hyperedges that is tight at $f$, contains $S$, and does not contain $e'$. Over $E\setminus U$, the sum of $f$-values is higher than the sum of $f_{e,i}$-values; therefore over $U$, the sum of $f_{e,i}$-values is higher. But this means that $f_{e,i}$ and $U$ contradict the inequality \eqref{eq:hypertree2} because the sum of $f$-values over $U$ is already $\mu(U)$, where $\mu$ is as in equation \eqref{eq:subm}.

Therefore $e'$ is internally inactive, but that contradicts our assumption on $f$. The only possible conclusion, then, is $f=f_{e,i}$.
\end{proof}

\begin{megj}
If we change the definition of $f_{e,i}$ in the proof above to now require that $M_{e,i}$ have the lowest possible multiplicity at $e_1$, then among those choices the lowest possible multiplicity at $e_2$ and so on, then a similar argument reveals that this, too, is a hypertree in which $e$ is the unique internally inactive hyperedge. Therefore the two descriptions define the same hypertree.
\end{megj}

It is not hard to see (cf.\ \cite[Section 44.6.c]{sch}) that Proposition \ref{pro:const} generalizes to arbitrary integer extended polymatroids, i.e., interior and exterior polynomials always have a constant term of $1$. Theorem \ref{pro:rang}, in turn, offers a way of defining the nullity of an integer extended polymatroid via its interior polynomial.

In the classical case, Proposition \ref{pro:const} translates to a well known property of the Tutte polynomial. Proposition \ref{pro:rang}, however, reveals (to the best of the author's knowledge) previously undiscovered information.

\begin{kov}
Let $G=(V,E)$ be a connected graph. By summing the coefficients of its Tutte polynomial $T_G(x,y)$ in front of terms that contain $x^{|V|-2}$, we obtain the nullity of $G$. In other words, under any order, the number of spanning trees that contain exactly one internally inactive edge is the first Betti number of $G$ as a one-dimensional cell complex.
\end{kov}

\subsection{Product formulas}
Let us now discuss situations when the interior and exterior polynomials behave multiplicatively. Our first claim is obvious.


\begin{lemma}\label{lem:farkinca}
From $\mathscr H=(V,E)$, construct another hypergraph $\mathscr H'$ by
\begin{enumerate}[(a)]
\item adding a singleton hyperedge $e'=\{v\}$ to $E$ for some $v\in V$;
\item adding a new vertex $v'$ to $V$ and making it part of exactly one hyperedge $e\in E$.
\end{enumerate}
Then $I_{\mathscr H'}=I_{\mathscr H}$ and $X_{\mathscr H'}=X_{\mathscr H}$.
\end{lemma}

\begin{proof}
In both cases, new and old hypertrees are in an obvious bijection. In the first, new hypertrees are exactly the extensions of old ones by assigning $0$ to $e'$. The new hyperedge $e'$ is both internally and externally active with respect to all hypertrees regardless of the order used.

In the second case, increase $f(e)$ by $1$ in each hypertree $f$ of $\mathscr H$ to get the collection of hypertrees in $\mathscr H'$. With respect to any order, this correspondence preserves $\bar\iota$ and $\bar\varepsilon$ values.
\end{proof}

The main result of this subsection examines the picture when two disjoint bipartite graphs are joined by identifying two of their edges.

\begin{tetel}\label{thm:edgesum}
Let $\mathscr H_1=(V_1,E_1)$ and $\mathscr H_2=(V_2,E_2)$ be hypergraphs so that $V_1\cap V_2=\{v\}$. Suppose that $e_1\in E_1$ and $e_2\in E_2$ both contain $v$ and form the hypergraph 
\[\mathscr H=(V,E)=\big(V_1\cup V_2,(E_1\setminus\{e_1\})\cup(E_2\setminus\{e_2\})\cup\{e\}\big)\] 
by merging $e_1$ and $e_2$
into a single hyperedge $e=e_1\cup e_2$. For hypertrees $f_1$ in $\mathscr H_1$ and $f_2$ in $\mathscr H_2$, define the function $f_1\#f_2\colon E\to\N$ by
\[(f_1\#f_2)(e)=f_1(e_1)+f_2(e_2)\]
and by otherwise letting $(f_1\#f_2)\big|_{E_i\setminus\{e_i\}}=f_i$, $i=1,2$. Then $\#$ defines a bijection
\[(Q_{\mathscr H_1}\cap\Z^{E_1})\times(Q_{\mathscr H_2}\cap\Z^{E_2})\cong Q_{\mathscr H}\cap\Z^{E}.\] Consequently, 
\[I_{\mathscr H}=I_{\mathscr H_1}I_{\mathscr H_2}\quad\text{and}\quad X_{\mathscr H}=X_{\mathscr H_1}X_{\mathscr H_2}.\]
\end{tetel}

\begin{proof}
For the claim on the sets of hypertrees, choose arbitrary elements $f_i\in Q_{\mathscr H_i}\cap\Z^{E_i}$ and realize them with spanning trees which contain the edges connecting $e_i$ and $v$. This can be done by Lemma \ref{lem:anchor}. After merging $e_1$ and $e_2$, the union of the two trees becomes a spanning tree of $\bip\mathscr H$ realizing $f_1\#f_2$, which is therefore a hypertree.

Conversely, if $f$ is a hypertree in $\mathscr H$ then we may realize it with a spanning tree that contains the edge between $e$ and $v$. This can be separated into spanning trees of $\bip\mathscr H_1$ and of $\bip\mathscr H_2$ which induce hypertrees in $\mathscr H_1$ and $\mathscr H_2$, respectively. It is easy to see that this defines an inverse to the correspondence $(f_1,f_2)\mapsto f_1\#f_2$.

To prove the claim on the polynomials, we make the following observation. If the hypertree $f=f_1\#f_2$ is such that valence can be transferred from $a_1\in E_1\setminus\{e_1\}$ to $a_2\in E_2\setminus\{e_2\}$, then 
\begin{enumerate}[(a)]
\item\label{belso} $f_1$ is such that valence can be transferred from $a_1$ to $e_1$ and
\item\label{kulso} $f_2$ is such that valence can be transferred from $e_2$ to $a_2$.
\end{enumerate}
Indeed by the above, the result $f'$ of the assumed transfer of valence has a decomposition $f'=f'_1\#f'_2$ where (since the entries in $f_i$ and $f'_i$ share the same sum) $f'_1$ is a hypertree that shows the truth of \eqref{belso} and the existence of $f'_2$ proves \eqref{kulso}.

Let us now fix an order on $E$ in which $e$ is smallest and use its restrictions to order $E_1$ and $E_2$ (so that $e_i$ becomes the smallest element in $E_i$). Then, we claim that for any pair of hypertrees $f_i$ in $\mathscr H_i$, we have
\begin{equation}\label{eq:addsup}
\bar\iota(f_1\#f_2)=\bar\iota(f_1)+\bar\iota(f_2)\quad\text{and}\quad\bar\varepsilon(f_1\#f_2)=\bar\varepsilon(f_1)+\bar\varepsilon(f_2)
\end{equation}
because, both in the internal and in the external sense, the set of inactive hyperedges for $f_1\#f_2$ is the union of the sets of inactive hyperedges for $f_1$ and for $f_2$. (As $e$, $e_1$, and $e_2$ are smallest, they can not be inactive with respect to any hypertree.) The backward inclusion is obvious because if, say, $f'_1$ results from $f_1$ by a single transfer of valence, then $f'_1\#f_2$ and $f_1\#f_2$ are related by essentially the same transfer. The forward inclusion follows equally easily with the help of \eqref{belso} and \eqref{kulso} above.

Comparing \eqref{eq:addsup} with the definitions of $I$ and $X$ completes the proof.
\end{proof}

The next two statements describe the situation when two bipartite graphs are joined at a single vertex. They follow from Theorem \ref{thm:edgesum} and Lemma \ref{lem:farkinca} by first slightly enlarging a hypergraph as in the Lemma and then joining it to another hypergraph as in the Theorem.

\begin{kov}\label{cor:eleterint}
Let $\mathscr H_1=(V_1,E_1)$ and $\mathscr H_2=(V_2,E_2)$ be hypergraphs with disjoint vertex sets. Choose $e_1\in E_1$ and $e_2\in E_2$ and form the hypergraph $\mathscr H$ with vertex set $V_1\cup V_2$ and hyperedge set $(E_1\setminus\{e_1\})\cup(E_2\setminus\{e_2\})\cup\{\,e_1\cup e_2\,\}$. Then $I_{\mathscr H}=I_{\mathscr H_1}I_{\mathscr H_2}$ and $X_{\mathscr H}=X_{\mathscr H_1}X_{\mathscr H_2}$.
\end{kov}

\begin{kov}\label{cor:pontoterint}
Let $\mathscr H_1=(V_1,E_1)$ and $\mathscr H_2=(V_2,E_2)$ be hypergraphs with disjoint vertex sets. Choose $v_1\in V_1$ and $v_2\in V_2$ and form the hypergraph $\mathscr H$ by identifying them to a single vertex $v$. I.e., the vertex set of $\mathscr H$ is $(V_1\setminus\{v_1\})\cup(V_2\setminus\{v_2\})\cup\{v\}$ and its hyperedge set is identified with $E_1\cup E_2$ so that $v$ is an element of a `new' hyperedge if and only if either $v_1$ or $v_2$ was an element of the `old' hyperedge. Then $I_{\mathscr H}=I_{\mathscr H_1}I_{\mathscr H_2}$ and $X_{\mathscr H}=X_{\mathscr H_1}X_{\mathscr H_2}$.
\end{kov}

\subsection{Deletion and contraction}\label{ssec:delcontr}
These operations are of special importance in the theory of the Tutte polynomial. However in the hypergraph case, so far they play surprisingly small roles. 

\begin{Def} 
Let $\mathscr H=(V,E)$ be a hypergraph and $e\in E$ a hyperedge. \emph{Deleting} $e$ from $\mathscr H$ means passing to the hypergraph $\mathscr H-e=(V,E\setminus\{e\})$. The result of \emph{contracting} $e$ is the hypergraph $\mathscr H/e$ which is obtained from $\mathscr H-e$ by identifying the elements of $e$. I.e., the hyperedge set of $\mathscr H/e$ is essentially $E\setminus\{e\}$ and its vertex set is $(V\setminus e)\cup\{\bar e\}$, where the new vertex $\bar e$ belongs to the hyperedge $e'\in E\setminus\{e\}$ if and only if $e'\cap e\ne\varnothing$.
\end{Def}

It would be more consistent with our previous notation to write $\mathscr H\setminus\{e\}$ instead of $\mathscr H-e$, but for this subsection we decided to conform to the existing literature and use the simpler symbol.

Regarding the hypertree polytope $Q_{\mathscr H}$ and the hyperedge $e$, let us define the \emph{top face} of $Q_{\mathscr H}$ with respect to $e$ as $\{\,f\in Q_{\mathscr H}\mid f(e)=|e|-1\,\}$ and let the \emph{bottom face} be $\{\,f\in Q_{\mathscr H}\mid f(e)=0\,\}$. Because $Q_\mathscr H$, being the set of bases in an integer polymatroid, is itself an integer polytope, the top and bottom faces with respect to any hyperedge coincide with the convex hulls of the hypertrees that they contain.

\begin{all}\label{pro:szendvics}
Let $\mathscr H=(V,E)$ be a hypergraph with hypertree polytope $Q_{\mathscr H}$. For any hyperedge $e\in E$, the hypertree polytopes of $\mathscr H-e$ and $\mathscr H/e$ are naturally isomorphic to the bottom and the top face, respectively, of $Q_{\mathscr H}$ with respect to $e$.
\end{all}

\begin{proof}
It suffices to equate the sets of hypertrees within the respective polytopes. 

Deletion: If we take any hypertree in $\mathscr H-e$ and an arbitrary realization, then by adding to it a single edge of $\bip\mathscr H$ adjacent to $e$, we realize a hypertree that lies along the bottom face. Conversely, any realization of a hypertree from the bottom face of $Q_\mathscr H$ has a single edge adjacent to $e$ and by removing it we obtain a tree which realizes a hypertree in $\mathscr H-e$. Therefore a bijection is defined by extending hypertrees in $\mathscr H-e$ to $e$ with the value zero. (It is possible for he bottom face of $Q_\mathscr H$ to be empty, namely when $\bip\mathscr H-\{e\}=\bip(\mathscr H-e)$ is disconnected. Of course in such a case $Q_{\mathscr H-e}$ is empty as well.) 

Contraction: Let $\bar e$ denote the new vertex in $\mathscr H/e$ that resulted from the contraction of $e$. If $f$ is a hypertree in $\mathscr H/e$, then, using Lemma \ref{lem:anchor}, construct a realization $\Gamma$ for it which contains all edges adjacent to $\bar e$. Based on this tree, we carry out the following construction in $\bip\mathscr H$. Keep all edges of $\Gamma$ that are not adjacent to $\bar e$. If $e'\in E\setminus\{e\}$ is such that $e'\cap e\ne\varnothing$, then connect $e'$ to an arbitrarily chosen element of the intersection. Finally, add all edges that are adjacent to $e$. The result is a spanning tree in $\bip\mathscr H$ so that its induced hypertree is part of the top face of $Q_\mathscr H$ and it agrees with $f$ on $E\setminus\{e\}$.

The inverse of this correspondence is constructed as follows. Let $g$ be a hypertree from the top face of $Q_\mathscr H$. Any of its realizations contains all edges of $\bip\mathscr H$ adjacent to $e$. Take one such tree and (viewing it as a topological space) contract the union of its edges adjacent to $e$ to a single point. The result is another tree which naturally embeds in $\bip(\mathscr H/e)$ as a spanning tree so that it realizes a hypertree in $\mathscr H/e$ that agrees with $g$ over $E\setminus\{e\}$.
\end{proof}

We now turn to generalizations of the classical deletion--contraction formulas \eqref{eq:delcontr}. Lemma \ref{lem:farkinca} can be viewed as a trivial extension of the second one to hypergraphs. (The reader may wish to consult \eqref{eq:regitutte} to see that there is no contradiction.) The first formula also has an easy extension to hypergraphs as follows. 

\begin{all}
Let $\mathscr H=(V,E)$ be a hypergraph with a hyperedge $e$ so that $\bip\mathscr H$ is connected but $\bip\mathscr H-\{e\}$ is the union of $|e|$ connected components. Then $I_\mathscr H=I_{\mathscr H/e}$ and $X_\mathscr H=X_{\mathscr H/e}$ (and both are equal to obvious $|e|$-fold products).
\end{all}

\begin{proof}
Immediate from Lemma \ref{lem:farkinca} and Corollaries \ref{cor:eleterint} and \ref{cor:pontoterint}.
\end{proof}

\begin{megj}
For the abstract duals of the hypergraphs in the previous Proposition, the same argument yields the same conclusion: $I_{\overline{\mathscr H}}=I_{\overline{\mathscr H/e}}$ and $X_{\overline{\mathscr H}}=X_{\overline{\mathscr H/e}}$, because both sides of each equation agree with the obvious $|e|$-fold product.
\end{megj}

Finally, let us generalize the third deletion--contraction formula from \eqref{eq:delcontr}.

\begin{all}\label{pro:delcontr}
Let $\mathscr H=(V,E)$ be a hypergraph that contains a two-element hyperedge $e$ which is so that $\bip\mathscr H-\{e\}$ is connected. Then we have
\begin{equation}\label{eq:elcsusznak}
I_\mathscr H(\xi)=I_{\mathscr H-e}(\xi)+\xi I_{\mathscr H/e}(\xi)\quad\text{and}\quad X_\mathscr H(\eta)=\eta X_{\mathscr H-e}(\eta)+X_{\mathscr H/e}(\eta).
\end{equation}
\end{all}

\begin{proof}
As $|e|=2$, each hypertree in $\mathscr H$ lies either on the top (if its value at $e$ is $1$) or on the bottom (if it is $0$) face of $Q_\mathscr H$ with respect to $e$. Let us choose an order on $E$ in which $e$ is the largest hyperedge. 

With respect to hypertrees along the bottom face, $e$ is internally active. With respect to those on the top face, by Lemma \ref{lem:mohofelett}, $e$ is internally inactive: indeed the Lemma applies because $nj(e)=1$ by our assumption which means that the value of the greedy hypertree at $e$ is $g(e)=|e|-1-nj(e)=0$. With this, our first equation follows easily from Proposition \ref{pro:szendvics}.

The second is equally easy if we notice that with respect to hypertrees of the top face, $e$ is externally active, whereas with respect to those along the bottom face, it is externally inactive. The latter can be argued for example like this: If $f$ is a hypertree in $\mathscr H$ so that $f(e)=0$, then take any of its realizations and add to it the `other' edge adjacent to $e$, too. The cycle thus created goes through at least one more hyperedge $e'$ and if we break the cycle by removing one of its edges adjacent to $e'$, we get a realization of a hypertree which is the result of a transfer of valence to $e$ from the smaller hyperedge $e'$.
\end{proof}

It would be highly desirable to formulate a version of the last proposition for hyperedges of arbitrary size. Such formulas, however, are currently lacking.

\section{Abstract duality}\label{sec:sejtes}

The observations on the interior polynomial contained in Propositions \ref{pro:const} and \ref{pro:rang} show more than just the (already proved) order-independence of the coefficients. The values observed also remain the same if we exchange the roles of hyperedges and vertices, that is they are invariant under abstract duality. We conjecture that this is true for the rest of the interior polynomial as well. In other words, the interior polynomial is in fact an invariant of bipartite graphs.

\begin{sejt}\label{conj:dual}
Let $G=(V_0,V_1,E)$ be a connected bipartite graph which induces the hypergraphs $\mathscr G_0$ and $\mathscr G_1$ as in \eqref{eq:absdual}. Then $I_{\mathscr G_0}=I_{\mathscr G_1}$.
\end{sejt}

So far, supporting evidence for this includes Postnikov's Theorem \ref{thm:post}, which says that $I_{\mathscr G_0}(1)=I_{\mathscr G_1}(1)$, and Theorem \ref{pro:rang}. We present some more below.

Conjecture \ref{conj:dual} was certainly true in Example \ref{ex:harom}, where it was also illustrated that the operation of abstract duality will usually result in a change in the exterior polynomial. (Indeed, the linear terms can already be different.) Let us work out a more general family of examples.

\begin{pelda}
We consider the complete bipartite graph $K_{m,n}=(V_0,V_1,E)$, where $|V_0|=n$ and $|V_1|=m$. For both $\mathscr G_0$ and $\mathscr G_1$ the condition \eqref{eq:hypertree2} in Theorem \ref{thm:politop} becomes vacuous and we see that their hypertree polytopes are the simplices 
\[Q_{\mathscr G_0}=(m-1)\Delta_{V_0}\quad\text{and}\quad Q_{\mathscr G_1}=(n-1)\Delta_{V_1},\] 
both of which contain ${n+m-2\choose n-1}$ hypertrees. (The number of lattice points in a standard $d$-dimensional simplex of sidelength $k$ is ${k+d\choose d}$.)

Fix now some order on $V_0$. If $f\colon V_0\to\N$ is a hypertree in $\mathscr G_0$ and $e\in V_0$, then $e$ is internally inactive with respect to $f$ if and only if $f(e)>0$ and $e$ is not the smallest in the order. Hence to enumerate hypertrees in $\mathscr G_0$ with internal inactivity $k$, we need to 
\begin{enumerate}[(a)]
\item choose the set of the $k$ internally inactive hyperedges and
\item partition the sum of hypertree entries, $m-1$, between them and the smallest hyperedge so that the latter may receive $0$ but the others get a positive amount.
\end{enumerate}
This is an easy exercise 
whose solution is ${n-1\choose k}{m-1\choose k}$. As this is a symmetric expression in $n$ and $m$, we see that Conjecture \ref{conj:dual} holds for complete bipartite graphs.

We leave it for the reader to verify that the number of hypertrees in $\mathscr G_0$ with external inactivity $k$ is ${m+k-2\choose k}$ for all $0\le k\le n-1$. For this of course it does matter if we interchange $n$ and $m$.
\end{pelda}

Let us temporarily call a bipartite graph \emph{good} if it satisfies Conjecture \ref{conj:dual}. Bipartite graphs with an isomorphism that interchanges their color classes are good. We have seen a number of operations that, when applied to good graphs, result in other good graphs. These include adding on a new valence one point (Lemma \ref{lem:farkinca}), joining two bipartite graphs at a vertex (Corollaries \ref{cor:eleterint} and \ref{cor:pontoterint}) or at an edge (Theorem \ref{thm:edgesum}). Adding a new valence two point is also such an operation. To see why, we have to accompany Proposition \ref{pro:delcontr} with the following result.

\begin{all}\label{pro:delcontr'}
Let $\mathscr H=(V,E)$ be a hypergraph with a vertex $v\in V$ belonging to exactly two hyperedges $e_1$ and $e_2$. Define the deletion of $v$ as $\mathscr H'=\overline{\overline{\mathscr H}-v}$ and the contraction of $v$ as $\mathscr H''=\overline{\overline{\mathscr H}/v}$ (which are the obvious ways). Assuming that $\bip\mathscr H-\{v\}$ is connected, we have
\begin{equation}\label{eq:megintcsusznak}
I_\mathscr H(\xi)=I_{\mathscr H'}(\xi)+\xi I_{\mathscr H''}(\xi).
\end{equation}
If we also suppose that $v$ is not the only common element of $e_1$ and $e_2$, we obtain the relation $X_\mathscr H(\eta)=X_{\mathscr H'}(\eta)+\eta X_{\mathscr H''}(\eta)$.
\end{all}

Notice that \eqref{eq:megintcsusznak} matches the first equation in \eqref{eq:elcsusznak} exactly. There is however a discrepancy between the formulas concerning exterior polynomials, not to mention the extra assumption (which, by the way, can not be dropped).

\begin{proof}
Just like in the proof of Proposition \ref{pro:delcontr}, first we are going to split the set of hypertrees in $\mathscr H$ in two so that the two sides are in one-to-one correspondences with hypertrees in $\mathscr H'$ and $\mathscr H''$, respectively, and then study activities. The geometry however will be different this time and therefore the argument more complex.

If $f$ is a hypertree in $\mathscr H'$, define $f'(e_1)=f(e_1)+1$ and $f'(x)=f(x)$ for other hyperedges. This is a hypertree in $\mathscr H$ because a realization of $f$ in $\bip\mathscr H'=\bip\mathscr H-\{v\}$ together with the edge of $\bip\mathscr H$ between $e_1$ and $v$ gives a realization of $f'$. The isometric embedding $f\mapsto f'$ is obviously injective.

Note that $\mathscr H''$ contains a distinguished hyperedge $e$ where $e_1$ and $e_2$ merge. If $g$ is a hypertree in $\mathscr H''$ then $\mathscr H$ has hypertrees which agree with $g$ away from $e_1$ and $e_2$ because the same edges that constitute a realization of $g$ in $\bip\mathscr H''$ form a two-component forest in $\bip\mathscr H$ which can be made a spanning tree by adding to it both edges adjacent to $v$. The set of these,
\[
L_g=\left\{\,h\in Q_\mathscr H\cap\Z^E\biggm|h\big|_{E\setminus\{e_1,e_2\}}=g\big|_{E\setminus\{e_1,e_2\}}\,\right\},\]
is analogous to the line segments that we used in the proof of Theorem \ref{thm:independent}. Let us now associate to $g$ the hypertree $g''\in L_g$ which has the highest value at $e_2$. 

The sets $L_g$, which we will also refer to as \emph{fibers}, lie on parallel lines. Furthermore, any hypertree $h$ in $\mathscr H$ is part of a set $L_g$. To see this, use Lemma \ref{lem:anchor} to construct a realization of $h$ containing both edges adjacent to $v$ and then contract those two edges to a point to obtain a spanning tree in $\bip\mathscr H''$ that induces the right $g$. 

Moreover, if $g_1\ne g_2$, then the lines containing $L_{g_1}$ and $L_{g_2}$ are disjoint (so that the $g$ above is uniquely determined by $h$). This is because $E\setminus\{\,e_1,e_2\,\}$ represents all but one of the hyperedges in $\mathscr H''$ and by \eqref{eq:hypertree3}, the last value (at $e$) is uniquely determined by the rest. In particular, the correspondence $g\mapsto g''$ is also injective.

There do not exist $f$ and $g$ so that $f'=g''$ because
\begin{enumerate}[(a)]
\item\label{atf'} For any $f$, the hypertree $f'$ is such that a transfer of valence is possible from $e_1$ to $e_2$: just repeat the same argument that $f'$ is a hypertree but with the roles of $e_1$ and $e_2$ reversed.
\item\label{atg''} At a hypertree of the form $g''$, the same transfer is never possible. 
\end{enumerate}

On the other hand, elements of $L_g\setminus\{g''\}$, i.e., all hypertrees $h$ in $\mathscr H$ that are such that a transfer of valence is possible from $e_1$ to $e_2$, are of the form $f'$ for some $f$. This is because each such $h$ has a realization that does not contain the edge $\gamma$ of $\bip\mathscr H$ connecting $e_2$ and $v$. To show this, take any realization $\Xi$ of $h$. If it does not contain $\gamma$, we are done. Otherwise, by Lemma \ref{lem:anchor}, we may assume that it contains both edges of $\bip\mathscr H$ adjacent to $v$. \label{xipage} By removing those two edges, we separate $\Xi$ to two smaller trees $\Xi_1$ and $\Xi_2$. Let $E_1\ni e_1$ be the set of hyperedges contained by $\Xi_1$ and let $E_2\ni e_2$ be contained by $\Xi_2$ so that $E$ is the disjoint union of $E_1$ and $E_2$.

There have to be elements of $E_2$ that are adjacent in $\bip\mathscr H$ to vertices of $\Xi_1$. This is because otherwise $\Xi\cap\big(\bip\mathscr H\big|_{E_2}\big)=\Xi_2\cup\{\gamma\}$ would be a spanning tree in $\bip\mathscr H\big|_{E_2}$, which by Lemma \ref{lem:erdo} would mean that $E_2$ is tight at $h$ and that would rule out the possibility of the assumed transfer of valence from $e_1$ to $e_2$. Whenever we locate this kind of element $a$ in $E_2$, we `switch it over' to $E_1$ by our usual method: add to $\Xi$ an edge $\alpha$ connecting $
a$ to a vertex of $\Xi_1$; this creates a unique cycle $C$ in $\Xi\cup\{\alpha\}$ which necessarily contains $\gamma$; remove the edge of $C$ adjacent to $a$ which is not $\alpha$.

By iterating this procedure, we move through realizations of $h$ with smaller and smaller associated sets $E_2$. Therefore after finitely many steps, the hyperedge $a$ that we operate on will be $e_2$ itself. The edge that we remove in that step is necessarily $\gamma$. This finishes the argument that a hypertree in $\mathscr H$ is always either of the form $f'$ or of the form $g''$ but never both.

In order to study activities, let us now order $E$ so that $e_1$ and $e_2$ are the smallest hyperedges. (Later we will use $e_1<e_2$ in the internal case and $e_2<e_1$ in the external case, but for now we leave this unspecified.) We will use this order in $\mathscr H$ to define $\bar\iota$ and $\bar\varepsilon$, 
as well as in $\mathscr H'$ to define $\bar\iota'$ and $\bar\varepsilon'$ 
values. In $\mathscr H''$, let $e$ be the smallest hyperedge and otherwise use the same order in the definitions of $\bar\iota''$ and $\bar\varepsilon''$. 
We start with two easy claims.

\begin{enumerate}[(i)]
\item\label{ftof'} If a hyperedge $a$ is inactive (in either sense) with respect to $f$ in $\mathscr H'$, then $a$ is inactive with respect to $f'$ in $\mathscr H$. This is because if the hypertrees $f_1$ and $f_2$ in $\mathscr H'$ are related by a single transfer of valence then 
$f_1'$ and $f_2'$ are related by the same transfer.
\item\label{g''tog} Regarding $\mathscr H''$, an implication of the opposite kind is easy: If the hyperedge $a\in E\setminus\{\,e_1,e_2\,\}$ is inactive with respect to some element of $L_g$ (in $\mathscr H$), then $a$ is inactive with respect to $g$ (in $\mathscr H''$) as well. Just take the hypertree that results from the transfer of valence which `makes' $a$ inactive and notice that it is part of a set $L_{g_0}$; this $g_0$ differs from $g$ by the single transfer of valence which establishes the claim.
\end{enumerate}

The rest of the proof is concerned with the degree to which the converses of \eqref{ftof'} and \eqref{g''tog} above are true. We start with \eqref{g''tog}.

\begin{lemma}\label{lem:sinpar}
If the hyperedge $a\in E\setminus\{\,e_1,e_2\,\}$ is inactive, in either sense, with respect to the hypertree $g$ in $\mathscr H''$, then $a$ is inactive with respect to any element $h\in L_g$ (e.g.\ $h=g''$) in $\mathscr H$. 
\end{lemma}

To show this, we separate two cases.
\begin{enumerate}[I.]
\item Let $h\in L_g$. If $g$ is such that a transfer of valence is possible from $a$ to $e$ (resp.\ from $e$ to $a$) resulting in the hypertree $g_0$ in $\mathscr H''$, then we claim that $h$ is such that a transfer of valence is possible from $a$ to $e_1$ or $e_2$ (resp.\ from $e_1$ or $e_2$ to $a$). This follows by fixing a $2$-dimensional plane containing $L_g$ and $L_{g_0}$ and constructing the type of elementary argument that we used to prove Lemmas \ref{lem:rombusz} and \ref{lem:nyil}.
\item If $g$ is such that no transfer valence is possible from $a$ to $e$ (resp.\ from $e$ to $a$), then by the argument in \eqref{g''tog} it follows that $h$ 
is such that there is no transfer of valence from $a$ to $e_1$ or to $e_2$ (resp.\ from $e_1$ or from $e_2$ to $a$). As $a$ is inactive, there is another hyperedge $b\in E\setminus\{\,e_1,e_2\,\}$, smaller than $a$, so that $g$ is such that valence can be transferred from $a$ to $b$ (resp.\ from $b$ to $a$), resulting in the hypertree $g_0$. From Lemma \ref{lem:rombusz} and \eqref{g''tog} we obtain that $h$ is such that valence can not be transferred from $b$ to $e_1$ or $e_2$ (resp.\ from $e_1$ or $e_2$ to $b$) either. Thus Lemma \ref{lem:teglalap}, applied to $p=e_1$, $q=e_2$, $r=a$, $s=b$ (resp.\ $p=a$, $q=b$, $r=e_1$, $s=e_2$) implies that $L_g$ and $L_{g_0}$ are located in a rectangular cross-section of $Q_\mathscr H$ so that every hypertree along $L_g$ is such that a transfer of valence is possible from $a$ to $b$ (resp.\ from $b$ to $a$). This completes the proof of the lemma.
\end{enumerate}

Let us assume that $e_1$ and $e_2$ have a second common vertex $v'\ne v$. This is equivalent to saying that for any hypertree $g$ in $\mathscr H''$, the set $L_g$ has at least two elements\footnote{As an alternative to the short argument below, the equivalence can also be deduced from Postnikov's Proposition \ref{pro:felbont}.}, which in turn makes it meaningful for us to let $g'''$ denote the element of $L_g$ adjacent to $g''$. (Note that $g'''$ is of the form $f'$ for some $f$.) Indeed, if $v'$ exists, then a spanning tree of $\bip\mathscr H$ may contain at most three edges of the quadrangle $e_1ve_2v'$ which makes it very easy to realize a transfer of valence between $e_1$ and $e_2$. Conversely, if $e_1\cap e_2=\{v\}$, then there are (greedy) hypertrees in $\mathscr H$ which take the maximal values $|e_1|-1$ and $|e_2|-1$ at $e_1$ and $e_2$, respectively.

In addition to our previous assumptions, we now set $e_2<e_1$. Concerning external inactivities, we prove the following items. 
\begin{itemize}
\item $\bar\varepsilon(g'')=\bar\varepsilon''(g)+1$ for all hypertrees $g$ in $\mathscr H''$. It is clear that $e_1$ is externally inactive with respect to $g''$ whereas of course $e$ is active with respect to $g$. The rest follows from Lemma \ref{lem:sinpar} and \eqref{g''tog}.
\item $\bar\varepsilon(f')=\bar\varepsilon'(f)$ for all hypertrees $f$ in $\mathscr H'$. A converse of \eqref{ftof'} is indeed true in this case. If $f'$ is such that $e_2$ can transfer valence to $e_1$, then the same is obviously true for $f$. If $f'$ is such that $e_1$ or $e_2$ can transfer valence to some other hyperedge $a$ then, even if the result of the transfer happens to be of the form $g''$, it is because the transfer came from $e_1$ and it is easy to see that a transfer of valence from $e_2$ to $a$ turns $f'$ into $g'''$. Finally, if $f'$ and $a$ are such that no transfer of valence is possible from $e_1$ or $e_2$ to $a$ but the hyperedge $b$ can transfer valence to $a$, resulting in the hypertree $h$, then the usual combination of Lemmas \ref{lem:rombusz} and \ref{lem:teglalap} shows that the fibers containing $f'$ and $h$ lie in a rectangular cross-section and therefore $h$ is not of the form $g''$ either.
\end{itemize}
The combination of the last two equations settles the claim in the Proposition concerning exterior polynomials.



Finally, we turn to internal activities and establish \eqref{eq:megintcsusznak}. For the remainder of the proof, we will use the order $e_1<e_2$ among the two smallest hyperedges. Under the simplifying assumption that all fibers contain at least two hypertrees, the result follows relatively easily as above in the external case. Without that assumption, both components of the argument break down: some hypertrees $g''$ are such that $e_2$ can not transfer valence to $e_1$, i.e., $g''$ does not pick up an extra internally inactive hyperedge beyond those of $g$; on the other hand there are also hypertrees $f'$ which do pick up extra internally inactive hyperedges beyond those of $f$ because some transfer of valence can turn them into one of the $g''$ without a convenient $g'''$ nearby. Our task is to show that the two problems cancel each other out.

Let us call a hypertree in $\mathscr H$ \emph{lonely} if it is the unique element of its fiber, i.e., if it is such that no transfer of valence is possible between $e_1$ and $e_2$. We claim that if $h$ is a lonely hypertree, then there is a hyperedge $a\in E\setminus\{\,e_1,e_2\,\}$ so that $h$ is such that both $e_1$ and $e_2$ can transfer valence to $a$. After noting that any realization $\Xi$ of $h$ has to contain both edges of $\bip\mathscr H$ adjacent to $v$, we proceed as earlier (on p.\ \pageref{xipage}) to define the subtrees and hyperedge sets $\Xi_i\supset E_i\ni e_i$, $i=1,2$. Now there has to be a hyperedge $a$ in $E_1$ or $E_2$ that is connected by an edge $\alpha$ of $\bip\mathscr H$ to a vertex in $\Xi_2$ or $\Xi_1$, respectively, for otherwise removing $v$ would disconnect the graph. Adding $\alpha$ to $\Xi$ creates a cycle through $v$ so that by removing either edge of $\bip\mathscr H$ adjacent to $v$, we realize the two desired transfers of valence to $a$.

If $h$ is a lonely hypertree then select the smallest hyperedge $a_h$ with the property above and let the hypertree $t(h)$ be obtained from $h$ by a single transfer of valence from $e_2$ to $a_h$. From the construction it is obvious that there is 
a hypertree $f$ in $\mathscr H'$ so that $f'=t(h)$. It is also clear from Lemma \ref{lem:rombusz} that $t(h)$ is part of a two-element fiber; in particular it is such that $e_2$ can not transfer valence to $e_1$. 

For a hypertree $f$ in $\mathscr H'$, let us introduce
\[T_f=\{\,h\in Q_\mathscr H\cap\Z^E\mid h\text{ is a lonely hypertree with }t(h)=f'\,\}.\]
The proof of \eqref{eq:megintcsusznak} will be immediate from the following items.

\begin{enumerate}[1)]
\item\label{OK} $\bar\iota(g'')=\bar\iota''(g)+1$ for all hypertrees $g$ in $\mathscr H''$ so that $g''$ is not lonely. This follows easily from \eqref{g''tog}, Lemma \ref{lem:sinpar}, and the observation that $e_2$ is internally inactive with respect to such hypertrees.
\item\label{keves} $\bar\iota(g'')=\bar\iota''(g)$ whenever $g''$ is lonely. In these cases $e_2$ is internally active; otherwise see above.
\item\label{tulsok} $\bar\iota(f')=\bar\iota'(f)+|T_f|$ for all hypertrees $f$ in $\mathscr H'$.
\item\label{sorozat} For any hypertree $f$ in $\mathscr H'$, the values of $\bar\iota$ among the elements of $T_f$ are $\bar\iota'(f),\bar\iota'(f)+1,\ldots,\bar\iota'(f)+|T_f|-1$, each occurring exactly once.
\end{enumerate}

We will prove \ref{tulsok}) and \ref{sorozat}) momentarily but let us first take care of \eqref{eq:megintcsusznak}. For each hypertree $g$ in $\mathscr H''$, we need a hypertree in $\mathscr H$ whose $\bar\iota$ is one higher than $\bar\iota''(g)$. If it is not lonely then $g''$ can play this role by \ref{OK}); otherwise \ref{keves}) says that $g''$ is not good but it belongs to a unique set $T_f$ which, by \ref{sorozat}), either contains an element with the desired $\bar\iota$ or if it does not then, according to \ref{tulsok}), $f'$ will do the job. Also, for each hypertree $f$ in $\mathscr H'$ we need a hypertree in $\mathscr H$ with the same $\bar\iota$. If $T_f=\varnothing$ then, by \ref{tulsok}), $f'$ will do; otherwise \ref{sorozat}) says that the element of $T_f$ with the lowest $\bar\iota$ will suffice. Since we used each hypertree in $\mathscr H$ exactly once, \eqref{eq:megintcsusznak} follows.

Let us now fix a hypertree $f$ in $\mathscr H$. For each $h\in T_f$, the hyperedge $a_h$ is internally inactive with respect to $f'$ because a transfer of valence from $a_h$ to $e_2$ turns $f'$ into $h$. We claim that the same hyperedges $a_h$ are internally active with respect to $f$, implying via \eqref{ftof'} that $\bar\iota(f')\ge\bar\iota'(f)+|T_f|$. For this it suffices to prove that any downward transfer of valence from $a_h$ will turn $f'$ into a hypertree of the form $g''$ for some $g$. The hypertree $f'$ is such that valence can not be transferred from $a_h$ to $e_1$ because $h$ was such that valence could not be transferred from $e_2$ to $e_1$. If $f'$ is such that a transfer of valence is possible from $a_h$ to some hyperedge $e_2<x<a_h$ (resulting in the hypertree $j$) then, by the way $a_h$ was chosen in the definition of $t(h)$, the hypertree $h$ is such that valence can not be transferred from $e_1$ to $x$. In other words, $j$ is such that valence can not be transferred from $e_1$ to $e_2$ which, by the characterization \eqref{atg''}, establishes our claim.

To see why $\bar\iota(f')\le\bar\iota'(f)+|T_f|$ holds too, we let $y\in E$ be internally inactive with respect to $f'$ and show that it is either internally inactive with respect to $f$ as well or it is one of the $a_h$ for some $h\in T_f$. If $f'$ is such that a downward transfer of valence from $y$ turns it into $f_0'$ for some $f_0$ then the first option holds. For instance, by Lemma \ref{lem:rombusz} and \eqref{atf'} applied to $f'$, this is the case if $f'$ is such that a transfer of valence is possible from $y$ to $e_1$. Assume now that $f'$ is such that valence can be transferred from $y$ to a hyperedge $e_2<x<y$ resulting in the hypertree $i$. By \eqref{atf'}, $f'$ is such that valence can be transferred from $e_1$ to $e_2$ (resulting in the hypertree $j$) but we may assume that the opposite is the case at $i$. That means that there is a set $U\subset E$ which is tight at $i$, contains $e_2$, but does not contain $e_1$. Because the same set $U$ is not tight at $f'$, it must be the case that $U$ contains $x$ and does not contain $y$. Next, apply Lemma \ref{lem:teglalap} with $p=e_2$, $q=x$, $r=e_1$, $s=y$ to find a rectangular cross-section containing $i$ and notice that it also contains $j$. This implies the existence of two more hypertrees in the cross-section, one of which shows that $f'$ is such that a transfer of valence is possible from $y$ to $e_2$. Let the result of that transfer be the hypertree $h$.

If $h$ is of the form $f_0'$ for some $f_0$ or if it is not but it is not lonely either then we see that $y$ is internally inactive with respect to $f$. If $h$ is lonely, then because $h$ is such that both $e_1$ and $e_2$ can transfer valence to $y$, we have $a_h\le y$; but if $a_h\ne y$ then $f'$ is such that $y$ can transfer valence to $a_h$ (the result being $t(h)$), showing that $y$ is internally inactive with respect to $f$. This finishes the proof of \ref{tulsok}).

What remains is to examine the set $T_f\cup\{f'\}$ of lattice points. Note that $f'=t(h)$ and the hyperedge $a_h$ determine the hypertree $h$. Let us label the elements of $T_f$ so that $a_{h_1}>a_{h_2}>\cdots>a_{h_{|T_f|}}$ holds. Other than the ones belonging to $e_2$ and the $a_{h_i}$, $i=1,2,\ldots,|T_f|$, all components of the elements of $T_f\cup\{f'\}$ are identical. The remaining components are almost the same too as each $h_i$ is derived from $f'$ by a transfer of valence from $a_{h_i}$ to $e_2$. From this description it is clear that $T_f\cup\{f'\}$ is the set of vertices of an inverted $|T_f|$-dimensional unit simplex. Based on this, using the same techniques as above, it is not hard to show that with respect to $h_i$, the hyperedges $e_2,a_{h_{|T_f|}},\ldots, a_{h_{i}}$ are internally active and $a_{h_{i-1}},\ldots,a_{h_1}$ are internally inactive, whereas for all other hyperedges 
we obtain that they are internally active either with respect to all elements of $T_f\cup\{f'\}$ or with respect to none of those. This completes the proof of \ref{sorozat}) and hence that of the Proposition.
\end{proof}

The list of operations before Proposition \ref{pro:delcontr'} reduces Conjecture \ref{conj:dual} to certain prime graphs, one of which is discussed in the next example.

\begin{pelda}
Consider the bipartite graph shown in Figure \ref{fig:peti}. It is in fact a plane bipartite cubic graph, which is a class that will be important for other reasons in Sections \ref{sec:trinity} and \ref{sec:moretrinity}. For now, we note that both of its induced hypergraphs have the interior polynomial
\[1+10\xi+48\xi^2+146\xi^3+302\xi^4+410\xi^5+277\xi^6+49\xi^7+\xi^8.\]
This is a non-trivial example of Conjecture \ref{conj:dual} because the structure is not symmetric: even though both color classes have nine elements with matching valences, there is no isomorphism that interchanges them. This is demonstrated by the slight difference between the corresponding exterior polynomials. If the full dots play the role of hyperedges, we obtain
\[1+8\eta+36\eta^2+110\eta^3+235\eta^4+344\eta^5+318\eta^6+162\eta^7+30\eta^8,\]
whereas if the hollow dots represent hyperedges, we get
\[1+8\eta+36\eta^2+110\eta^3+235\eta^4+348\eta^5+326\eta^6+159\eta^7+21\eta^8.\]
All three polynomials are outputs of a computer code written by P\'eter Juh\'asz.
\begin{figure}[htbp] 
   \centering
   \includegraphics[width=2in]{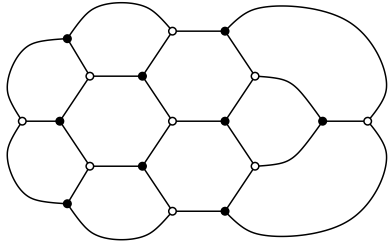} 
   \caption{A plane bipartite cubic graph.}
   \label{fig:peti}
\end{figure}
\end{pelda}

\section{Planar duality}\label{sec:dual}

We call a hypergraph $\mathscr H=(V,E)$ \emph{planar} if the corresponding bipartite graph $\bip\mathscr H$ is planar, i.e., if it admits an embedding into the $2$--sphere $S^2$. We will almost always assume that such an embedding is fixed, but for extra clarity, we will talk about \emph{plane} (hyper)graphs to mean a graph together with a particular embedding. A \emph{region} of a plane (hyper)graph is a connected component of the complement of the image of the embedding. If the graph is connected, each region is homeomorphic to a disk.

For the rest of the paper, it will be convenient to allow multiple edges in (plane) bipartite graphs. Such objects induce the same two hypergraphs (and therefore the same two hypertree polytopes and the same interior and exterior polynomials) as the bipartite graph which results from reducing each set of multiple edges to a single edge. This is all fairly natural given that a spanning tree may contain at most one element from a set of multiple edges.

\begin{Def}\label{def:planedual}
Let $\mathscr H$ be a plane hypergraph so that its associated bipartite graph $\bip\mathscr H$ is connected. We define the \emph{planar dual hypergraph} $\mathscr H^*$ of $\mathscr H$ by keeping the set $E$ of hyperedges and letting the new set of vertices be the set $R$ of regions of $\mathscr H$. We let a region belong to a hyperedge if it is incident with the point representing it on the plane.
\end{Def}

\begin{figure}[htbp]
   \centering
   \includegraphics[width=\linewidth]{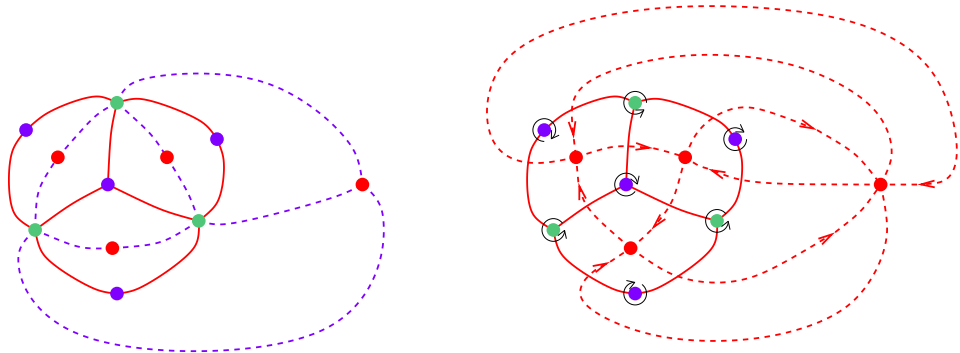} 
   \caption{Left: the planar dual of a hypergraph $\mathscr H$ of three hyperedges. Right: the planar dual of the corresponding bipartite graph $\bip\mathscr H$.}
   \label{fig:dualisok}
\end{figure}

The planar dual hypergraph $\mathscr H^*$ is also planar because the bipartite graph $\bip\mathscr H^*$ can be represented by a planar diagram where we place a new vertex in each element $r\in R$ and connect it to all points around the boundary $\partial r$ that represent elements of $E$. See Figure \ref{fig:dualisok} for an example. In fact we will identify $R$ and the set of points just introduced and think of the planar dual hypergraph as given with the embedding that we have just outlined.

It is easy to see that our notion generalizes the usual duality of plane graphs. Indeed, a plane graph $G=(V,E)$ is also a hypergraph. Passing to $\bip G$ means placing a new vertex to the midpoint of each edge. Now in $G^*$, these same midpoints are connected to the points which represent the two regions on either side of the original edge. Finally, if we forget the midpoint, these two connections merge to form the usual dual edge.

The following is also almost immediate.

\begin{all}\label{pro:dualdual}
For a plane hypergraph $\mathscr H$, the planar dual of the planar dual of $\mathscr H$ is naturally isomorphic to $\mathscr H$.
\end{all}

\begin{proof}
We explained after Definition \ref{def:planedual} how $\bip\mathscr H^*$ is embedded in $S^2$. Each point representing a vertex $v\in V$ is of course in some connected component of the complement of $\bip\mathscr H^*$. This correspondence is onto and one-to-one.
\end{proof}

The hypertree polytopes $Q_\mathscr H$ and $Q_{\mathscr H^*}$ of a planar dual pair are subsets of the same Euclidean space $\R^E$. Next, let us address the effect of planar duality on the hypertree polytope, as well as the interior and exterior polynomials.

\begin{tetel}\label{thm:sikdualis}
Let the hypergraphs $\mathscr H$ and $\mathscr H^*$ be planar duals. Then $Q_\mathscr H$ and $Q_{\mathscr H^*}$ are centrally symmetric. Furthermore, up to interchanging the indeterminates $\xi$ and $\eta$, we have $I_{\mathscr H^*}=X_{\mathscr H}$ and $X_{\mathscr H^*}=I_{\mathscr H}$.
\end{tetel}

In Figure \ref{fig:dualisok} we see that the hypergraph $\mathscr G_0$ of our earlier examples is planar and in fact planar self-dual. This explains why \eqref{eq:egyik} showed its interior and exterior polynomials to coincide.

Parts of the following proof were simplified from an earlier version by A.\ Bene.

\begin{proof}[Proof of Theorem \ref{thm:sikdualis}]
Proposition \ref{pro:dualdual} implies that for the claim on the polytopes, it suffices to show that if $f$ is a hypertree in $\mathscr H$ then $f^*(e)=|e|-1-f(e)$ defines a hypertree in $\mathscr H^*$. (So the center in which $Q_\mathscr H$ and $Q_{\mathscr H^*}$ are symmetric becomes the point whose $e$-coordinate is $\frac12(|e|-1)$ for all $e\in E$.)

To this end, choose a spanning tree $\Gamma\subset\bip\mathscr H$ that realizes $f$. Form its planar dual $\Gamma^d$ in the classical sense, i.e., put a vertex in each element of $R$ and connect these by edges that bisect each edge of $\bip\mathscr H$ that is not an edge in $\Gamma$. This is well known to be a spanning tree in the dual graph, in particular it is cycle-free, connected, and contains all points of $R$. We are going to use a continuous deformation followed by some discrete steps to turn $\Gamma^d$ into a spanning tree $\Gamma^*\subset\bip\mathscr H^*$ that realizes $f^*$.

Each edge $\gamma^d$ of $\Gamma^d$ bisects an edge $\gamma$ of $\bip\mathscr H$ and $\gamma$ has exactly one end $e$ in $E$. Let us push the midpoint of $\gamma^d$ along $\gamma$ so that it gets very close to $e$. If $r\in R$ and $e\in E$ are adjacent with two edges $\gamma_1,\gamma_2$ of $\bip\mathscr H$, then we merge the halves of $\gamma_1^d$ and $\gamma_2^d$ adjacent to $r$, as shown in the right side of Figure \ref{fig:csillag}. The result of the modifications so far is still a plane tree which we will continue to denote with $\Gamma^d$.

\begin{figure}[htbp] 
\labellist
\footnotesize
\pinlabel $e$ at 275 273
\pinlabel $e$ at 1050 273
\pinlabel $r$ at 410 490
\pinlabel $\gamma_2$ at 375 550
\pinlabel $\gamma_2^d$ at 300 550
\pinlabel $\downarrow$ at 300 505
\pinlabel $\gamma_1$ at 500 475
\pinlabel $\gamma_1^d$ at 515 400
\pinlabel $\leftarrow$ at 470 400
\endlabellist
   \centering
   \includegraphics[width=4in]{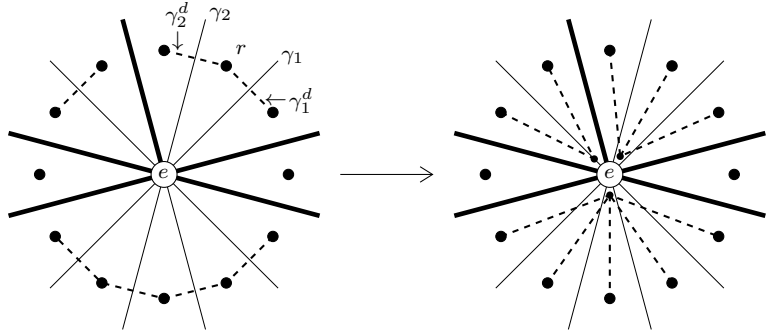} 
   \caption{Modifying a tree $\Gamma^d$ in $(\bip\mathscr H)^*$ toward a tree $\Gamma^*$ in $\bip\mathscr H^*$. Solid edges represent $\bip\mathscr H$; the thick ones are in $\Gamma$ and the thin ones are not.}
   \label{fig:csillag}
\end{figure}

Now let $e\in E$ be arbitrary. Write the edges of $\bip\mathscr H$ emanating from $e$, in a cyclic order determined by the embedding, as $\gamma_1,\ldots,\gamma_{|e|}$. Let the edges among these that do not belong to $\Gamma$ form $k$ maximal subsequences of consecutive elements. Let the lengths of these subsequences be $s_1,\ldots,s_k$. In our construction so far, each of these $k$ groups resulted in a point $e_i$ near $e$ with $s_i+1$ edges connecting it to elements of $R$. See Figure \ref{fig:csillag}, which shows $k=3$ with $s_1=2$, $s_2=1$, and $s_3=4$. We wish to push these points $e_1,\ldots,e_k$ all the way to $e$ so that the edges adjacent to them become edges of $\bip\mathscr H^*$.

The simplest case is that of $k=1$ when we can just identify $e_1$ with $e$. Note that after this, the number of edges of $\Gamma^d$ adjacent to $e$ is $s_1+1=|e|-(f(e)+1)+1=f^*(e)+1$. When $k\ge2$, we still start with pushing $e_1$ to $e$. When we do the same to $e_2$, however, a unique cycle is created in $\Gamma^d$ and we remove one of its two edges adjacent to $e$ in order to have a tree again. Then we push $e_3$ to $e$ and remove another edge, and so forth until the last point $e_k$. In the tree that results at the end, the valence of $e$ is $(s_1+1)+s_2+\cdots+s_k=|e|-(f(e)+1)+1=f^*(e)+1$. After applying this procedure to all elements of $e\in E$ with $k\ge1$, we arrive at a tree that is adjacent to all such hyperedges as well as to all elements of $R$. 

We intentionally left out the case $k=0$ which of course corresponds to $f^*(e)=0$. To finish the construction of the spanning tree $\Gamma^*\subset\bip\mathscr H^*$ realizing $f^*$, we add one more edge for each such hyperedge so that it is connected to an arbitrarily chosen adjacent element of $R$.

The claim on the polynomials is now easily obtained 
by the following observation. The correspondence $f\leftrightarrow f^*$ is a bijection between hypertrees in $\mathscr H$ and $\mathscr H^*$. The hypertree $f$ can be transformed into $g$ by a single transfer of valence if and only if $f^*$ can be transformed into $g^*$ by the opposite transfer. Therefore with respect to any order of the hyperedges, we have $\bar\iota(f)=\bar\varepsilon(f^*)$.
\end{proof}

\begin{Def}\label{def:hyperdual}
If the hypertree $f$ in the plane hypergraph $\mathscr H$ and the hypertree $f^*$ in $\mathscr H^*$ are related as above (namely, $f(e)+f^*(e)=|e|-1$ for all hyperedges $e$), then we will call them \emph{planar dual hypertrees}.
\end{Def}

\section{Trinities}\label{sec:trinity}

\subsection{Basic observations}\label{ssec:basic}
\begin{figure}[htbp]
\labellist
\footnotesize
\pinlabel $r_0$ at 706 334
\pinlabel $r_1$ at 318 388
\pinlabel $r_2$ at 102 388
\pinlabel $r_3$ at 210 208
\pinlabel $e_0$ at 56 244
\pinlabel $e_1$ at 382 262
\pinlabel $e_2$ at 219 496
\pinlabel $v_0$ at 38 441
\pinlabel $v_1$ at 201 334
\pinlabel $v_2$ at 218 117
\pinlabel $v_3$ at 417 424
\pinlabel $t_1$ at 110 450
\pinlabel $t_2$ at 110 330
\pinlabel $t_3$ at 160 175
\pinlabel $t_4$ at 260 260
\pinlabel $t_5$ at 250 410
\pinlabel $t_6$ at 490 480
\pinlabel $t_7$ at 390 360
\pinlabel $t_8$ at 470 170
\endlabellist
   \centering
   \includegraphics[width=4in]{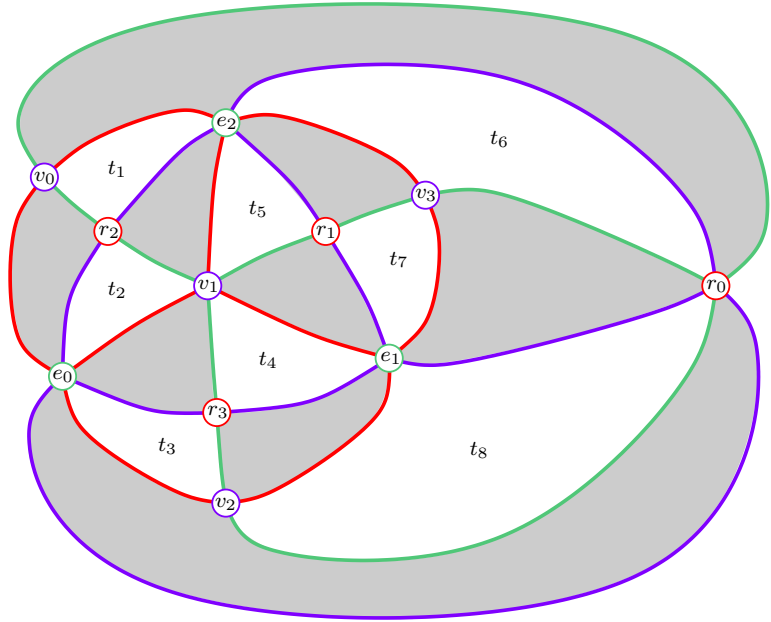} 
   \caption{A trinity of plane bipartite graphs.}
   \label{fig:trinity}
\end{figure}

Starting from a plane hypergraph and iterating the constructions of planar and abstract duality, a total of six plane hypergraphs can be built. The corresponding bipartite graphs form a triple; see Figure \ref{fig:trinity} for an example. The picture thus obtained has a rich combinatorics and admits several equivalent definitions (cf.\ Remark \ref{rem:pbcg}). To highlight the perfect symmetry between the constituent parts, we choose to build our discussion around the following notion.

\begin{Def}
A \emph{trinity} is a triangulation of the sphere $S^2$ together with a three-coloring of the $0$-simplices. (I.e., $0$-simplices joined by a $1$-simplex have different colors.) According to dimension, we will refer to the simplices as \emph{points, edges}, and \emph{triangles}.
\end{Def}

Most claims of this subsection can be repeated for three-colored triangulations of other closed, orientable surfaces. The material in the rest of the section, however, is quite specific to the sphere. We hope to return to this point in our forthcoming joint paper with Bene.

We will use the names red, emerald, and violet for the colors in the trinity and denote the respective sets of points with $R$, $E$, and $V$. Let us color each edge in the triangulation with the color that does not occur among its ends. Then $E$ and $V$ together with the red edges form a bipartite graph that we will call the \emph{red graph} and denote with $G_R$. Each region of the red graph contains a unique red point. Likewise, the \emph{emerald graph} $G_E$ has red and violet points, emerald edges, and regions marked with emerald points. Finally, the \emph{violet graph} contains $R$ and $E$ as vertices, violet edges, and a violet point in each of its regions.

\begin{Def}
We will refer to the plane bipartite graphs $G_R$, $G_E$, and $G_V$ above as the \emph{constituent bipartite graphs} of the trinity. The hypergraphs induced by the constituent bipartite graphs are said to be \emph{contained} in the trinity.
\end{Def}

A trinity contains six plane hypergraphs. It can be uniquely reconstructed from each of the six as follows.

\begin{all}
Let $T$ be a trinity with points set $R\cup E\cup V$ as above. Let $\mathscr H=(V,E)$ be the hypergraph with emerald hyperedges and violet vertices so that $\bip\mathscr H=G_R$. Then we also have $(R,E)=\mathscr H^*$, $(E,R)=\overline{\mathscr H^*}$, $(V,R)=\left(\overline{\mathscr H^*}\right)^*=\overline{\left(\overline{\mathscr H}\right)^*}$, $(R,V)=\left(\overline{\mathscr H}\right)^*$, and $(E,V)=\overline{\mathscr H}$.
\end{all}

\begin{proof}
Obvious from an examination of Figure \ref{fig:trinity}.
\end{proof}



The triangles of a trinity can also be colored but not with the original three colors. Notice that each triangle is adjacent with exactly one edge and one point of each color. Compared to the orientation of the sphere, the cyclic order of the colors around each triangle may be positive or negative. If two triangles share an edge, these orientations are opposite. 
Hence the triangles have a black and white checkerboard coloring according to orientation, cf.\ Figure \ref{fig:trinity}. In particular, the dual graph of a trinity is a plane bipartite cubic (i.e., uniform trivalent) graph. It turns out that the converse is also true.

\begin{all}\label{pro:pbcg}
Planar duals of plane bipartite cubic graphs are three-colorable.
\end{all}

\begin{proof}
Let $\theta$ be a plane bipartite cubic graph and let us call its color classes black and white. Direct each edge in $\theta$ from its black to its white endpoint. Now if we assign the modulo $3$ coefficient $1$ to the thus directed edges, the result is a cycle $C$ and hence a homology class in $H_1(S^2,\Z_3)=0$. Therefore the modulo $3$ intersection number of a closed loop with $C$ is zero.

Choose a base region of $\theta$ and assign the remainder class $0\in\Z_3$ to it. Given any other region, connect it to the base and assign to the region the modulo $3$ intersection number of the path and $C$. By the above remark, this is well defined and it is obviously a three-coloring with the colors $0$, $1$, and $2$.
\end{proof}

\begin{megj}\label{rem:pbcg}
According to Proposition \ref{pro:pbcg} and the paragraph above it, the notion of a trinity is equivalent to that of a plane bipartite cubic graph.
\end{megj}

Finally, notice that the sets of red edges, emerald edges, violet edges, white triangles, and black triangles all have the same cardinality $n$. In particular, adjacency defines natural bijections between white triangles and edges of each color. Now if we apply Euler's formula to the trinity $T$, we get $|R|+|E|+|V|-3n+2n=2$, that is
\begin{equation}\label{eq:euler}
\text{ the total number of points exceeds that of the white triangles by }2.
\end{equation}

\subsection{The Tree Trinity Theorem}
Next, form the planar dual $G^*_{R}$ of the red graph $G_{R}$. This is now in the classical sense, i.e., the vertex set of $G^*_{R}$ is $R$ and its edges are in a one-to-one correspondence with the red edges of the trinity. This graph has a natural orientation (more precisely, two natural orientations that are opposites of each other), defined as follows. Put a positive spin (as in a small spinning top) to the emerald points and a negative spin to the violet points. If $e\in E$ and $v\in V$ are connected by a (red) edge, then the two spins induce the same orientation on the dual edge. See Figure \ref{fig:dualisok} for an illustration. Notice that at each red point $r$, the edges of  $G^*_{R}$ that are adjacent to it are oriented toward and away from $r$ in an alternating fashion. In particular, each vertex of  $G^*_{R}$ has the same number of incoming and outgoing edges.

\begin{figure}[htbp]
   \centering
   \includegraphics[width=3in]{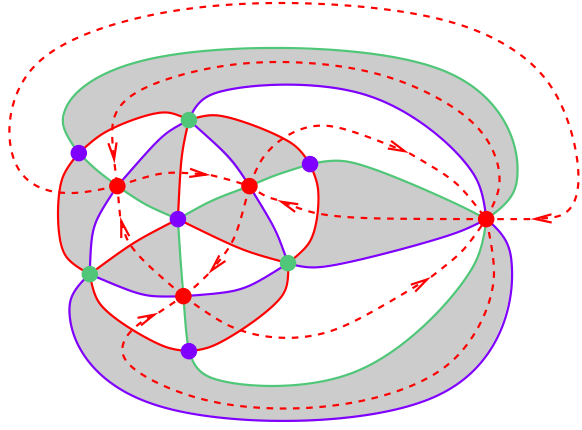} 
   \caption{The trinity of Figure \ref{fig:trinity} and the directed graph $G^*_{R}$.}
   \label{fig:dualitas}
\end{figure}

Similar properties hold of course for the graphs  $G^*_{E}$ and  $G^*_{V}$. It is in fact this triple of directed graphs that Tutte called a trinity. In \cite{tutte1}, he proved a general property of directed graphs which will be stated after the next definition.

\begin{Def}
Let $D$ be a finite
directed graph (possibly with multiple edges) and fix a vertex $r\in D$, called the \emph{root}. A spanning tree in $D$ is called a \emph{spanning arborescence} with respect to $r$ if the unique path in the tree from $r$ to any other vertex is oriented toward that vertex.
\end{Def}
 
Directed graphs with the same in-degree and out-degree everywhere are called \emph{balanced}.

\begin{tetel}[Tutte]
Let $D$ be a finite, balanced directed graph. The number of its spanning arborescences does not depend on the choice of root. The same holds for trees that are oriented toward the root and the counts of the two kinds of tree coincide.
\end{tetel}

The invariant of balanced directed graphs introduced in the Theorem above is denoted with $\rho(D)$ and called the \emph{arborescence number}. The following will be a useful characterization.

\begin{lemma}\label{lem:tree}
Spanning arborescences of a finite directed graph with root $r$ are exactly those cycle-free subgraphs which also have the property that for any vertex other than $r$, exactly one edge of the subgraph is directed toward that vertex. 
\end{lemma}

\begin{proof}
Let $D=(V,E)$ be a finite directed graph. If it is disconnected then both properties above describe the empty set. If $D$ is connected and $A$ is a spanning arborescence with respect to $r$, then of course it is cycle-free and it has at least one edge directed into any vertex other than $r$, namely the last edge on the path from $r$ to that vertex. But because $A$ contains exactly $|V|-1$ edges, these already exhaust all of them.

Conversely, if a subgraph is cycle-free than it may have at most $|V|-1$ edges. If it has one incoming edge for all points other than $r$, then it has to have exactly $|V|-1$ of them (so it is a spanning tree) and none can be directed into $r$. For any vertex $v\ne r$, the first edge of the unique path from $r$ to $v$ has to be directed toward $v$. The same holds for the second edge lest they share terminal points with the first. Etc.
\end{proof}

Let us now return to trinities and recall a beautiful result by Tutte.

\begin{tetel}[Tree Trinity Theorem]\label{thm:ttt}
The arborescence numbers $\rho(G^*_{R})$, $\rho(G^*_{E})$, and $\rho(G^*_{V})$ of the directed graphs associated to a trinity $T$ are the same.
\end{tetel}

We will call this common value the \emph{arborescence number of the trinity} and denote it with $\rho(T)$. We quickly sketch a proof, based on \cite{tutte3}, as elements of it will be relevant later. 

\begin{proof}
Recall that the triangles of a trinity have a black-and-white color scheme. We may concretize our choice of orientation in the three directed graphs by requiring that the head of each edge be in a white triangle while the tail end of each is in a black triangle. Compare Figures \ref{fig:dualisok} and \ref{fig:dualitas}. Then by Lemma \ref{lem:tree}, a spanning arborescence of $G^*_{R}$, say, is an assignment of an adjacent white triangle to each non-root red point (the unique incoming edge will reach the point through the chosen triangle). We adapt this point of view for the rest of the proof.

Let us now distinguish a white triangle $t_0$ (called the outer triangle) and choose its adjacent points $r_0\in R$, $e_0\in E$, and $v_0\in V$ as the roots in the three directed graphs. Then $t_0$ can never be selected into an arborescence of any color. Notice that by \eqref{eq:euler}, the rest of the triangles and the rest of the points are in a one-to-one correspondence. 

Fix a spanning arborescence $A$ in $G^*_{R}$. We are going to associate to it two other spanning arborescences, in $G^*_{E}$ and $G^*_{V}$ respectively, so that 
\begin{equation}\label{eq:bijection}
\parbox{3.5in}{the union of the three mappings is a bijection from the set of non-root points to the set of non-outer white triangles.}
\end{equation}
In other words, no white triangle will get assigned to two different points. We also claim that, given $A$, there is a unique way to do this.

The dual $A^*$ of $A$ is a spanning tree in $G_{R}$. This tree $A^*$ consists of those edges whose duals are not in $A$, in other words of those red edges whose adjacent white triangles have not been assigned to a red point. Note that the edge connecting $e_0$ and $v_0$ is part of $A^*$. 

As $A^*$ is a tree, it has at least two leaves, i.e., valence one points. If those are $e_0$ and $v_0$, then $A^*$ consists of a single edge, hence there can not be any other emerald or violet points and we have no assigning to do. Otherwise, as our leaf is adjacent to a unique white triangle that has not yet been assigned, we make the obvious assignment and delete the leaf (with its single adjacent edge) from $A^*$. We continue this procedure until $A^*$ is reduced to the single edge of $t_0$ between $e_0$ and $v_0$. By that time, all other emerald and violet points have been associated to an adjacent white triangle so that no white triangle was used twice.

We have to check (cf.\ Lemma \ref{lem:tree}) that the subgraphs we constructed in $G^*_{E}$ and $G^*_{V}$ are cycle-free. This is shown using proof by contradiction, arguing that a violet or emerald cycle would prevent some red point from being connected to the red root $r_0$ in $A$. See \cite{tutte3} for details.

As the color red played no special role above, any triple of spanning arborescences satisfying \eqref{eq:bijection} can be uniquely reconstructed from any of its three members. We conclude that spanning arborescences of our three directed graphs are arranged in disjoint triples and so their numbers agree. 
\end{proof}

A mapping satisfying \eqref{eq:bijection} will be called a \emph{Tutte matching} (with respect to some fixed outer white triangle $t_0$). Tutte matchings are exactly the nonzero terms in the expansion of a determinant that first appeared in the work of Berman.

\begin{tetel}[Berman]\label{thm:berman}
The common arborescence number of the three directed graphs of a trinity can be expressed as the absolute value of the determinant of the adjacency matrix of non-outer white triangles and non-root points. 
\end{tetel}

This is of course a square matrix by the observation \eqref{eq:euler}. To prove the theorem, Berman \cite{berman} checks that any Tutte matching is in fact the union of three spanning arborescences (which boils down to cycle-freeness and another argument by contradiction) and that all non-zero expansion terms come with the same sign. 

\begin{pelda}
The adjacency matrix associated with the trinity of Figure \ref{fig:trinity} is
\[
M=
\bordermatrix
{~&\text{\footnotesize $t_1$}&\text{\footnotesize $t_2$}&\text{\footnotesize $t_3$}&\text{\footnotesize $t_4$}&\text{\footnotesize $t_5$}&\text{\footnotesize $t_6$}&\text{\footnotesize $t_7$}&\text{\footnotesize $t_8$}\cr
\text{\footnotesize $r_1$}&0&0&0&0&1&0&1&0\cr
\text{\footnotesize $r_2$}&1&1&0&0&0&0&0&0\cr
\text{\footnotesize $r_3$}&0&0&1&1&0&0&0&0\cr
\text{\footnotesize $e_1$}&0&0&0&1&0&0&1&1\cr
\text{\footnotesize $e_2$}&1&0&0&0&1&1&0&0\cr
\text{\footnotesize $v_1$}&0&1&0&1&1&0&0&0\cr
\text{\footnotesize $v_2$}&0&0&1&0&0&0&0&1\cr
\text{\footnotesize $v_3$}&0&0&0&0&0&1&1&0\cr}.
\]
Its determinant is $7$, which equals the number of spanning arborescences in each of the three dual directed graphs, including the one shown in Figure \ref{fig:dualitas}.
\end{pelda}

\section{Hypertrees in a trinity}\label{sec:moretrinity}

\subsection{Hypertrees and arborescences}
The main result of this last section is the following.

\begin{tetel}\label{thm:alexander}
Let $G=(V_0,V_1,E)$ be a plane bipartite graph. The arborescence number of its planar dual $G^*$ agrees with the number of hypertrees in the hypergraphs $\mathscr G_0=(V_1,V_0)$ and $\mathscr G_1=(V_0,V_1)$ induced by $G$. I.e., 
\[\rho(G^*)=|Q_{\mathscr G_0}\cap\Z^{V_0}|=|Q_{\mathscr G_1}\cap\Z^{V_1}|.\]
\end{tetel}

\begin{proof}
Let us denote the vertex set of $G^*$ with $R$ and fix a root $r_0\in R$. Any spanning arborescence $A$ rooted at $r_0$ has a planar dual spanning tree $A^*$ in $G$ and that in turn has a valence distribution which gives rise to a hypertree $f_A\colon V_0\to\N$. We will show that the mapping $A\mapsto f_A$ is one-to-one and onto. (A similar statement holds of course for hypertrees in $\mathscr G_1$ instead of $\mathscr G_0$. Thus what we find as a byproduct is a set of spanning trees in $G$ that simultaneously realize all hypertrees in $\mathscr G_0$ and $\mathscr G_1$.)

For any $v\in V_0$, the edges of $G$ adjacent to it give rise to a directed cycle $C_v$ in $G^*$ and the edge set of $G^*$ is the disjoint union of these cycles. Note that each cycle $C_v$ travels around the corresponding point $v$ following the same orientation.

First we show that if two arborescences $A$ and $B$ contain the same number of edges from each $C_v$ (which is equivalent to $f_{A}=f_{B}$), then $A=B$. Assume this is not so. Then at least one cycle $C_{v_1}$ contains an edge $\alpha_1$ of $A$ that is not an edge in $B$. Let the head of $\alpha_1$ be $r_1$. Obviously $r_1\ne r_0$ because $A$ has an edge pointing to it. The other arborescence $B$ also has to have such an edge $\beta_1$ but that has to be part of a different cycle $C_{v_2}$. The edge $\beta_1$ cannot belong to $A$ because then $A$ would have two edges pointing to $r_1$. Then, because $A$ and $B$ intersect $C_{v_2}$ in the same number of edges, $C_{v_2}$ also contains an edge $\alpha_2$ which belongs to $A$ but does not belong to $B$.

Iterating our argument, we find cycles $C_{v_1},C_{v_2},C_{v_3},\ldots$ in $G^*$ until the first repetition occurs. After discarding some cycles at the beginning and relabeling the rest, we may assume that $C_{v_1}=C_{v_{k+1}}$ was the first such coincidence. Then we find ourselves in the situation depicted in Figure \ref{fig:korbeer}. The interiors of the cycles $C_{v_i}$ are disjoint from $G^*$, hence the non-root points $r_1,\ldots,r_k$ separate $G^*$ into the smaller directed graphs $D_1$ and $D_2$, see Figure \ref{fig:korbeer}, which have only those points in common. Both $D_1$ and $D_2$ have to contain vertices other than $r_1,\ldots,r_k$ because otherwise $A$ or $B$ would not be cycle-free. Only one of them, however, can contain the root $r_0$. Now the contradiction is apparent from the fact that no directed path in $A$ can pass from $D_1$ to $D_2$ and no directed path in $B$ can pass from $D_2$ to $D_1$.

\begin{figure}[htbp]
\labellist
\footnotesize
\pinlabel $D_1$ at -100 300
\pinlabel $D_2$ at 320 260
\pinlabel $C_{v_1}$ at 100 260
\pinlabel $r_1$ at 110 350
\pinlabel $\alpha_1$ at 190 340
\pinlabel $\beta_1$ at 100 430
\pinlabel $C_{v_2}$ at 210 430
\pinlabel $r_2$ at 295 430
\pinlabel $\alpha_2$ at 310 360
\pinlabel $\beta_2$ at 330 500
\pinlabel $C_{v_3}$ at 420 410
\pinlabel $r_3$ at 510 400
\pinlabel $\alpha_3$ at 520 320
\pinlabel $\beta_3$ at 570 450
\pinlabel $\alpha_{k-2}$ at 540 110
\pinlabel $\beta_{k-2}$ at 500 10
\pinlabel $C_{v_{k-1}}$ at 400 100
\pinlabel $r_{k-1}$ at 305 75
\pinlabel $\alpha_{k-1}$ at 270 160
\pinlabel $\beta_{k-1}$ at 255 15
\pinlabel $C_{v_k}$ at 180 100
\pinlabel $r_k$ at 130 160
\pinlabel $\alpha_k$ at 170 200
\pinlabel $\beta_k$ at 70 170
\endlabellist
   \centering
   \includegraphics[width=2in]{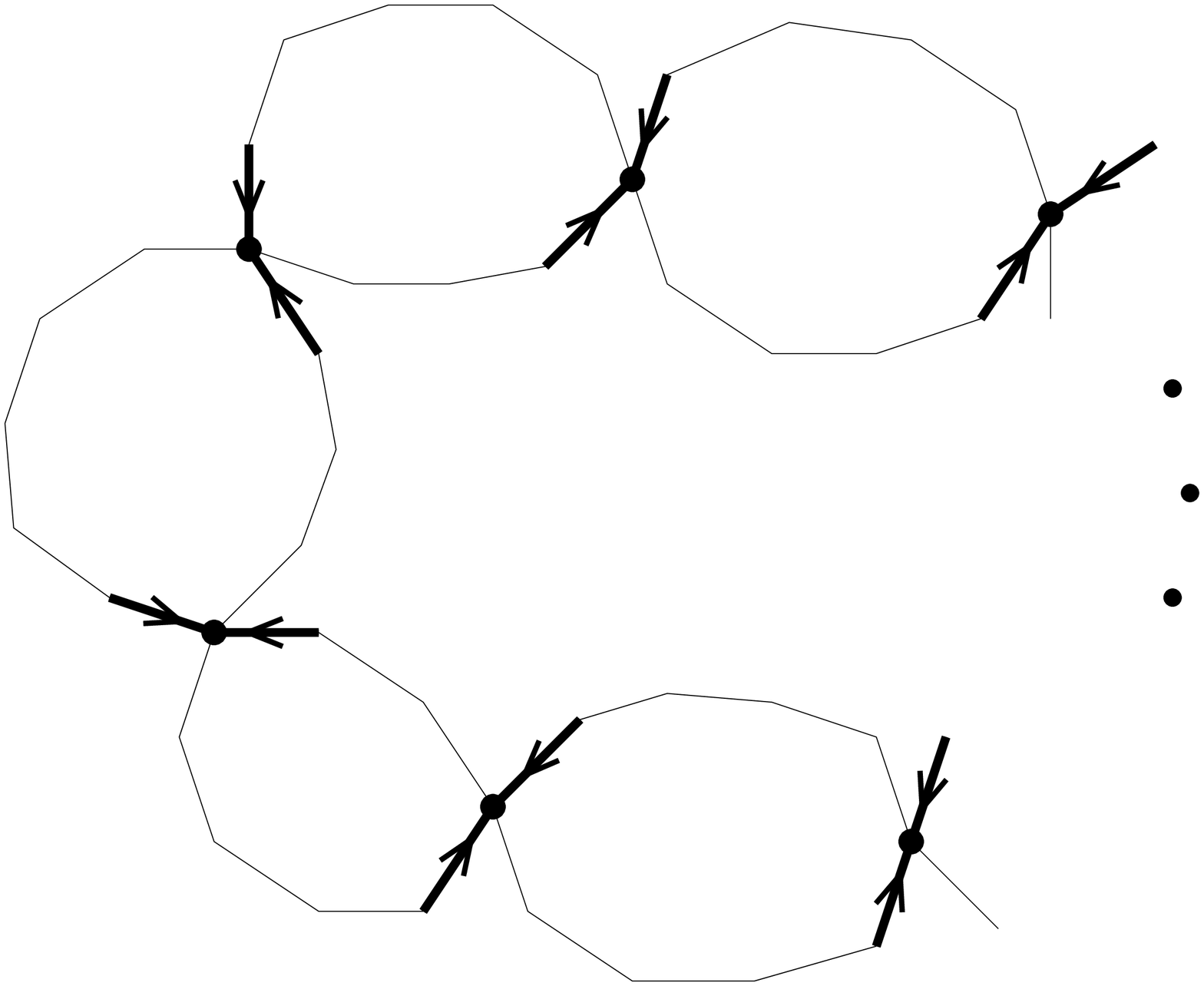} 
   \caption{A cycle of cycles in the directed graph $G^*$.}
   \label{fig:korbeer}
\end{figure}

To prove that any hypertree $f\colon V_0\to\N$ can be obtained in the form $f=f_A$, we employ the following strategy. Start from an arbitrary spanning tree $\Gamma$ in $G$ that realizes $f$ and take its dual tree $\Gamma^*\subset G^*$. We are going to change $\Gamma$ step-by-step through other realizations of the same hypertree $f$ so that the dual moves closer and closer to being an arborescence. It is sufficient to keep track of the changes made to the dual tree $\Gamma^*$. Then, the condition that we preserve the hypertree translates to requiring that for each $v\in V_0$, the number of edges in $C_v\cap\Gamma^*$ stays invariant throughout the process.

It is of course crucial to describe carefully what we mean by getting closer to an arborescence. In a rooted tree, edges $\gamma$ have a well defined distance $d(\gamma)\in\N$ to the root (those adjacent to the root have $d=0$ etc.). If the tree is directed, then there is also a clear sense for each edge to be directed toward the root or away from the root. The first kind of edge will be called \emph{bad} and the second kind \emph{good}. The tree is an arborescence if and only if all of its edges are good. Let us now associate the following quantities to a finite, directed, rooted tree $(T,r)$.
\begin{itemize}
\item Let $n(T,r)$ denote the smallest value of $d$ among bad edges. (So that within a radius of $n$ from the root, the tree is an arborescence.)
\item For $1\le m\le n(T,r)$, let $\lambda_{T,r}(m)$ be the number of edges $\gamma$ with $d(\gamma)=m-1$. These values are positive. For $m>n(T,r)$, we define $\lambda_{T,r}(m)=0$. 
\end{itemize}
Then, for a pair of rooted trees we write $(T_1,r_1)\prec(T_2,r_2)$ if either
\begin{enumerate}
\item\label{lexi} the sequence $\lambda_{T_1,r_1}(1),\lambda_{T_1,r_1}(2),\lambda_{T_1,r_1}(3),\ldots$ is smaller than the sequence $\lambda_{T_2,r_2}(1),\lambda_{T_2,r_2}(2),\lambda_{T_2,r_2}(3),\ldots$ in lexicographic order, or
\item\label{tuske} the two sequences coincide (implying $n(T_1,r_1)=n(T_2,r_2)=n$) but the number of bad edges in $T_1$ with $d=n$ is higher than the same count in $T_2$.
\end{enumerate}

This may sound complicated but the actual idea of the proof is very simple. If $\Gamma^*$ is not a spanning arborescence already, then it contains a bad edge $\gamma$ with $d(\gamma)=n(\Gamma^*,r_0)$. If we remove $\gamma$ from $\Gamma^*$, then the tree falls apart into a root component and a non-root component. Now $\gamma$ is part of a cycle $C_v$ which contains points from both components. Therefore $C_v$ also has an edge $\gamma'$ that goes from a point of the root component to a point of the non-root component. Let ${\Gamma^*}'$ denote the tree obtained from $\Gamma^*$ by replacing $\gamma$ with $\gamma'$. 

Note that $\gamma'$ is a good edge of the new tree. The difficulty is that in what used to be the non-root component, distances to the root and relative orientations may have changed. What we do know though is that each of those edges is farther away from the root than $\gamma'$. Let us separate two cases according to whether the non-root component is reattached `close to' or `far from' the root.
\begin{enumerate}[(i)]
\item If $d(\gamma')<n(\Gamma^*,r_0)$ (where $d(\gamma')$ is of course measured in ${\Gamma^*}'$), then $(\Gamma^*,r_0)\prec({\Gamma^*}',r_0)$ by \eqref{lexi}. Indeed, with the addition of $\gamma'$, $\lambda(d(\gamma')+1)$ went up by $1$ whereas the $\lambda$ values before it stayed the same.
\item If $d(\gamma')\ge n(\Gamma^*,r_0)$, then $(\Gamma^*,r_0)\prec({\Gamma^*}',r_0)$ either by \eqref{lexi} (if $\gamma$ was the unique bad edge of $\Gamma^*$ with $d(\gamma)=n(\Gamma^*,r_0)$, in which case $\lambda(n(\Gamma^*,r_0)+1)$ moves from $0$ to a positive value) or by \eqref{tuske}. 
\end{enumerate}
We conclude that if the dual of a realization of $f$ is not a spanning arborescence, then it is possible to change it to another tree, dual to another realization, that is larger in the sense of the linear order $\prec$. Because there are altogether finitely many spanning trees in $G^*$, it is now obvious that a spanning arborescence will be reached after finitely many improvements. Thus we find that our map $A\mapsto f_A$ is onto.
\end{proof}

Theorem \ref{thm:alexander} of course verifies Postnikov's Theorem \ref{thm:post} in the case of a planar bipartite graph $G$. Postnikov  \cite[Lemma 12.6]{post} describes a sufficient condition on a set of spanning trees of $G$ in order for them to define a triangulation of the root polytope and therefore to simultaneously realize all hypertrees of the two induced hypergraphs. It is interesting to note that our set of spanning trees satisfies his condition, as follows.

\begin{all}
Let $G$ be a plane bipartite graph with (classical) dual $G^*$ as before. Let $A$ and $B$ denote two spanning arborescences of $G^*$ with respect to the same root $r$. Then there is no cycle $\alpha_1,\beta_1,\alpha_2,\beta_2\ldots,\alpha_k,\beta_k$ in $G$ composed of edges $\alpha_1,\alpha_2,\ldots,\alpha_k$ of $A^*$ and $\beta_1,\beta_2,\ldots.\beta_k$ of $B^*$.
\end{all}

\begin{proof}
It is easy to see that edges of $A$ could only cross such a cycle in one direction and edges of $B$ could only cross it in the opposite direction. That leads to the usual contradiction of some points becoming inaccessible from the root for one of the two arborescences.
\end{proof}

Tutte's Theorem \ref{thm:ttt} says that directed graphs in a trinity have the same number of spanning arborescences, just like the two graphs in a planar dual pair have the same number of spanning trees. We may now add that the six hypergraphs contained in the trinity have the same number of hypertrees. Moreover, the following is immediate from Theorem \ref{thm:alexander}.

\begin{kov}\label{cor:rho}
The arborescence number $\rho$ associated to a trinity is also the number of hypertrees in any of the six hypergraphs contained in the trinity.
\end{kov}

Let us make one more observation. We saw how each Tutte matching contains three spanning arborescences (one for each of $G_R^*$, $G_E^*$, and $G_V^*$) and how the (classical) planar duals of these realize a hypertree in each of the six hypergraphs. This way, each hypertree turns up in relation to exactly one Tutte matching. In fact we also have the following.

\begin{all}
In any trinity, the set of the six hypertrees induced by a Tutte matching consists of three pairs of planar duals (in the sense of Definition \ref{def:hyperdual}).
\end{all}

\begin{proof}
Pick an arbitrary, say emerald point from the trinity and denote it with $e$. It is matched to one adjacent white triangle (or, if the point is the root, substitute the outer white triangle here). Out of its other $|e|-1$ adjacent white triangles, $m$ are matched to red points and $|e|-1-m$ to violet points. Therefore there are $|e|-m$ red edges and $m+1$ violet edges adjacent to $e$ that appear in the spanning trees of $G_R$ and $G_V$, respectively, induced by the Tutte matching. This causes the $e$-coordinates of the corresponding hypertrees to be $|e|-m-1$ and $m$, respectively, the sum of which is indeed $|e|-1$.
\end{proof}

\subsection{Determinant formulas}
We are going to extend Berman's Theorem \ref{thm:berman} to write hypertree polytopes associated with plane hypergraphs in determinant form. Let $M$ denote the adjacency matrix with rows indexed by the non-root points and columns indexed by the non-outer white triangles. If a triangle $t_i$ is adjacent with the emerald point $e_j$ and the non-root violet point $v_k$, then at the intersection of row $v_k$ and column $t_i$, change the entry $1$ to $e_j$. (After it becomes a matrix entry, we will think of $e_j$ as an indeterminate associated with the original point.) 

Call this matrix the \emph{enhanced adjacency matrix} and denote it with $M_{e\to v}$. Notice that even though $M$ contains no row indexed by a root, indeterminates belonging to roots do appear in enhanced adjacency matrices.

It would be rather pointless to write the symbol $e_j$ in the row indexed by that point because that would simply multiply the determinant by $e_j$. The little twist that we introduced, on the other hand, turns out to be quite useful.

\begin{tetel}\label{thm:det}
Let the plane bipartite graphs $G_{R}$, $G_{E}$, and $G_{V}$ be the constituents of a trinity. Fix an outer white triangle $t_0$ which is adjacent to the roots $r_0\in R$, $e_0\in E$, and $v_0\in V$ and form the enhanced adjacency matrix $M_{e\to v}$ as above. Then for the hypergraph $\mathscr H=(V,E)$ and its hypertree polytope, we have
\[Q_{\mathscr H}=\det M_{e\to v}\]
in the following sense. The determinant on the right hand side is a sum of monomials in the indeterminates $e\in E$. Either each monomial has coefficient $+1$ or each has $-1$. If we write the exponents in each monomial as a vector, the set we obtain is exactly the left hand side.
\end{tetel}

As a special case, we may use this theorem to write the spanning tree polytope of a plane graph as a determinant, too.

\begin{proof}
The claim on the uniform signs of course follows from Berman's Theorem \ref{thm:berman} on the adjacency matrix $M$. There and in the proof of Theorem \ref{thm:ttt} we saw that the nonzero terms in the expansion of $\det M$, which are the Tutte matchings (of non-root points in the trinity to adjacent non-outer white triangles), are also triples of spanning arborescences in the directed graphs $G^*_{R}$, $G^*_{E}$, and $G^*_{V}$. Furthermore, all spanning arborescences from each of the three directed graphs occur as part of exactly one triple.

Then in the proof of Theorem \ref{thm:alexander} we found that the planar duals $A^*$ of the spanning arborescences $A$ of $G^*_{R}$ realize each hypertree in $Q_{\mathscr H}$ exactly once. Recall that the value of the hypertree $f_A$ at an emerald point $e_j$ is the number of adjacent (to $e_j$) edges of $A^*$ minus one. A red edge is in $A^*$ if and only if its adjacent white triangle is not assigned to its red point, i.e., if the triangle is the outer one or if it is assigned to an emerald or violet point. 

If $e_j\ne e_0$ then it is not adjacent to $t_0$ and out of all white triangles adjacent to $e_j$, exactly one is assigned to an emerald point (namely, $e_j$). Thus $f_A(e_j)$ is the number of white triangles that are adjacent to $e_j$ and which are assigned to their violet points. The same holds for $e_0$ where now $t_0$ plays the role of the correction term $-1$.

Finally, we just have to check the effect of enhancing $M$ into $M_{e\to v}$ on the expansion term that corresponds to $A$. In that monomial, by definition, the indeterminate $e_j$ will appear once for each time that a (non-outer) white triangle $t_i$, which is adjacent to the point $e_j$, gets assigned to its adjacent violet point $v_k$. This completes the proof.
\end{proof}

\begin{pelda}
The hypertree polytopes of Example \ref{ex:ketto} can be recomputed as the determinants of the following two matrices $M_{e\to v}$ and $M_{v\to e}$ (warning: the labels of the points have changed, cf.\ Figure \ref{fig:trinity}):
\[
\bordermatrix
{~&\text{\footnotesize $t_1$}&\text{\footnotesize $t_2$}&\text{\footnotesize $t_3$}&\text{\footnotesize $t_4$}&\text{\footnotesize $t_5$}&\text{\footnotesize $t_6$}&\text{\footnotesize $t_7$}&\text{\footnotesize $t_8$}\cr
\text{\footnotesize $r_1$}&0&0&0&0&1&0&1&0\cr
\text{\footnotesize $r_2$}&1&1&0&0&0&0&0&0\cr
\text{\footnotesize $r_3$}&0&0&1&1&0&0&0&0\cr
\text{\footnotesize $e_1$}&0&0&0&1&0&0&1&1\cr
\text{\footnotesize $e_2$}&1&0&0&0&1&1&0&0\cr
\text{\footnotesize $v_1$}&0&e_0&0&e_1&e_2&0&0&0\cr
\text{\footnotesize $v_2$}&0&0&e_0&0&0&0&0&e_1\cr
\text{\footnotesize $v_3$}&0&0&0&0&0&e_2&e_1&0\cr};
\bordermatrix
{~&\text{\footnotesize $t_1$}&\text{\footnotesize $t_2$}&\text{\footnotesize $t_3$}&\text{\footnotesize $t_4$}&\text{\footnotesize $t_5$}&\text{\footnotesize $t_6$}&\text{\footnotesize $t_7$}&\text{\footnotesize $t_8$}\cr
\text{\footnotesize $r_1$}&0&0&0&0&1&0&1&0\cr
\text{\footnotesize $r_2$}&1&1&0&0&0&0&0&0\cr
\text{\footnotesize $r_3$}&0&0&1&1&0&0&0&0\cr
\text{\footnotesize $e_1$}&0&0&0&v_1&0&0&v_3&v_2\cr
\text{\footnotesize $e_2$}&v_0&0&0&0&v_1&v_3&0&0\cr
\text{\footnotesize $v_1$}&0&1&0&1&1&0&0&0\cr
\text{\footnotesize $v_2$}&0&0&1&0&0&0&0&1\cr
\text{\footnotesize $v_3$}&0&0&0&0&0&1&1&0\cr}.
\]
The first is equal to $e_0^2 e_1 + e_0 e_1^2 + e_0^2 e_2 + e_0 e_1 e_2 + e_1^2 e_2 + e_0 e_2^2 + e_1 e_2^2$, and the second one is $v_0 v_1 + v_1^2 + v_0 v_2 + v_1 v_2 + v_0 v_3 + v_1 v_3 + v_2 v_3$. Observe that the sequences of exponents ($(2,1,0)$ etc.\ and $(1,1,0,0)$ etc., respectively) are exactly the hypertrees that we saw earlier.
\end{pelda}

Of course, as we have already done in the example above, we may use Theorem \ref{thm:det} to compute any of the six hypertree polytopes associated to a trinity. This leads to another proof of Theorem \ref{thm:sikdualis} in terms of manipulating determinants. It is quick but it is probably much less instructive than the proof given in Section \ref{sec:dual}. 

Indeed, starting from $M_{e\to v}$, pull out $e_j$ from any column that belongs to a white triangle adjacent with $e_j$. Do this for all emerald points $e_j$. Then for the non-root emerald points, pull out $e_j^{-1}$ from the row indexed by $e_j$. The result is a monomial times the enhanced adjacency matrix $M_{e\to r}$ but with $e_j^{-1}$ written everywhere where $e_j$ should appear. Therefore the determinant of the matrix is $-Q_{\mathscr H^*}$. 

If we examine the monomial factor, we see that the exponent of $e_j$ in it is exactly $|e_j|-1$. 
Here $|e_j|$ means the size of $e_j$ as a hyperedge, which is the same in $\mathscr H$ as in $\mathscr H^*$. It can also be described as the number of violet and red points connected to $e_j$ by red, respectively violet edges of the trinity, as well as the number of white triangles adjacent to $e_j$. The claim is obvious for the non-root emerald points and the exponent of $e_0$ is $|e_0|-1$ because $e_0$ is adjacent to the outer white triangle $t_0$ which had no column in the matrix.

\subsection{Summary and two final notes}
\begin{figure}[htbp]
\labellist
\footnotesize
\pinlabel hypergraph at 370 580
\pinlabel {interior polynomial} at 370 636
\pinlabel {exterior polynomial} at 370 692
\pinlabel polytope at 370 524
\pinlabel $V,E$ at 526 497
\pinlabel $I$ at 574 533
\pinlabel $X$ at 622 569
\pinlabel $R,E$ at 565 335
\pinlabel $X$ at 617 323
\pinlabel $I$ at 669 311
\pinlabel $E,R$ at 459 202
\pinlabel $X'$ at 484 152
\pinlabel $Y$ at 509 102
\pinlabel $V,R$ at 281 202
\pinlabel $Y$ at 256 152
\pinlabel $X'$ at 231 102
\pinlabel $R,V$ at 175 335
\pinlabel $Y'$ at 123 323
\pinlabel $I'$ at 71 311
\pinlabel $E,V$ at 214 497
\pinlabel $I'$ at 166 533
\pinlabel $Y'$ at 118 569
\pinlabel $Q$ at 370 380
\pinlabel $Q_E$ at 505 410
\pinlabel $Q_V$ at 235 410
\pinlabel $Q_R$ at 370 237
\endlabellist
   \centering
   \includegraphics[width=3.6in]{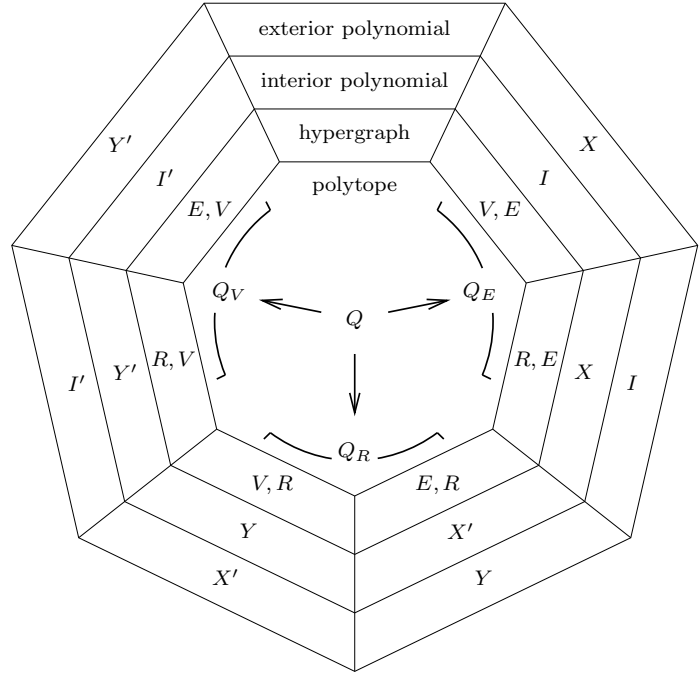} 
   \caption{Relations between the hypertree polytopes, interior, and exterior polynomials of the six hypergraphs contained in a trinity. If Conjecture \ref{conj:dual} is true, then $I=I'$, $X=X'$, and $Y=Y'$.}
   \label{fig:hetszog}
\end{figure}

Let us summarize our findings on hypertree polytopes and interior and exterior polynomials associated to the six hypergraphs contained in a trinity. See Figure \ref{fig:hetszog}. The six polytopes form three centrally symmetric pairs by Theorem \ref{thm:sikdualis}, i.e., there are only three `shapes' associated to the trinity. Let us denote these, indicating the respective sets of hyperedges, with $Q_R$, $Q_E$, and $Q_V$. These polytopes have different dimensions but by Corollary \ref{cor:rho} each contains the same number of integer lattice points, namely $\rho$, where $\rho$ is the arborescence number of the trinity.

A priori, there are twelve interior and exterior polynomials associated to the six hypergraphs. But again by Theorem \ref{thm:sikdualis}, that number is reduced to six, as indicated in Figure \ref{fig:hetszog}. These polynomials all have non-negative integer coefficients with their sum equal to $\rho$ in each case. Furthermore, if Conjecture \ref{conj:dual} holds, then $I=I'$, $X=X'$, and $Y=Y'$, so there are in fact only three different polynomials.

Our first `final note' is that the polytopes $Q_{R}$, $Q_{E}$, and $Q_{V}$ are all projections of a single set $Q$ of lattice points. This can be defined by superimposing $M_{e\to v}$, $M_{v\to r}$, and $M_{r\to e}$ to obtain the matrix $M_{e\to v\to r}$ and taking its determinant. In our running example, the result is
\[
M_{e\to v\to r}=
\bordermatrix
{~&\text{\footnotesize $t_1$}&\text{\footnotesize $t_2$}&\text{\footnotesize $t_3$}&\text{\footnotesize $t_4$}&\text{\footnotesize $t_5$}&\text{\footnotesize $t_6$}&\text{\footnotesize $t_7$}&\text{\footnotesize $t_8$}\cr
\text{\footnotesize $r_1$}&0&0&0&0&v_1&0&v_3&0\cr
\text{\footnotesize $r_2$}&v_0&v_1&0&0&0&0&0&0\cr
\text{\footnotesize $r_3$}&0&0&v_2&v_1&0&0&0&0\cr
\text{\footnotesize $e_1$}&0&0&0&r_3&0&0&r_1&r_0\cr
\text{\footnotesize $e_2$}&r_2&0&0&0&r_1&r_0&0&0\cr
\text{\footnotesize $v_1$}&0&e_0&0&e_1&e_2&0&0&0\cr
\text{\footnotesize $v_2$}&0&0&e_0&0&0&0&0&e_1\cr
\text{\footnotesize $v_3$}&0&0&0&0&0&e_2&e_1&0\cr},
\]
with a determinant of $e_0^2 e_1 r_0^2 v_0 v_1^2 + e_0 e_1^2 r_0 r_3 v_0 v_1 v_2 + 
 e_1^2 e_2 r_1 r_2 v_1^2 v_2 + e_0^2 e_2 r_0 r_1 v_0 v_1 v_3 + 
 e_0 e_2^2 r_0 r_2 v_1^2 v_3 + e_0 e_1 e_2 r_1 r_3 v_0 v_2 v_3 + 
 e_1 e_2^2 r_2 r_3 v_1 v_2 v_3$. Obviously, $Q_E=\det M_{e\to v}$, $Q_V=\det M_{v\to r}$, and $Q_R=\det M_{r\to e}$ are all obtained from $\det M_{e\to v\to r}$ by substituting $1$ for the unneeded variables. This corresponds to projections of the associated sets of lattice points which we indicated with the three arrows in Figure \ref{fig:hetszog}. In particular, all four determinants contain the same number of lattice points, which is the arborescence number of the trinity. It is not yet clear, however, whether $Q=\det M_{e\to v\to r}$ is itself a lattice polytope.

Finally, there exists a curious relationship between our results and an invariant introduced by Jaeger \cite{jaegerpoly}. We noted in subsection \ref{ssec:basic} that the planar dual of a trinity is a plane bipartite cubic graph. Jaeger associated a one-variable polynomial $S(u)$ (with non-negative integer coefficients) to such objects. His definition does not make use of the three-coloring of Proposition \ref{pro:pbcg} and indeed, his polynomial is different from any of the six (three) that we introduced. (For example, the trinity of Figure \ref{fig:trinity} has $S(u)=6u^3+u^5$.) The sum of the coefficients is, however, the same.

\begin{all}
Let $T$ be a trinity. Then $S(T^*,1)=\rho(T)$, the arborescence number of $T$.
\end{all}

\begin{proof}
As the definition of $S$ is based on a recursion, it is natural to prove our claim by induction. Jaeger sets his initial condition $S(u)=1$ at a `free loop,' but an equivalent notion results if we start by setting $S(u)=u$ for the theta graph (which has two vertices and three edges connecting them). The theta graph is dual to the trinity formed by one black and one white triangle with a common boundary. 
Each of the three directed graphs associated to the latter has a single vertex and one loop edge and, therefore, only one arborescence.

\begin{figure}[htbp]
\labellist
\small
\pinlabel $T$ at -40 400
\pinlabel $T_0$ at -40 30
\pinlabel $T$ at 720 400
\pinlabel $T_1$ at 670 10
\pinlabel $T_2$ at 1110 10
\pinlabel $e$ at 880 300
\pinlabel $v'$ at 755 260
\pinlabel $v'$ at 1000 30
\pinlabel $v''$ at 1020 440
\pinlabel $v''$ at 1255 205
\pinlabel $v$ at 610 100
\endlabellist
   \centering
   \includegraphics[width=3in]{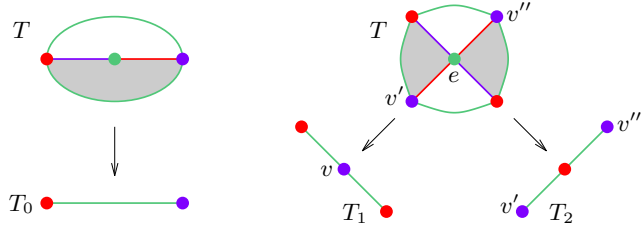} 
   \caption{Collapsing moves used in the definition of Jaeger's polynomial.}
   \label{fig:collapse}
\end{figure}

The invariant $S$ obeys two kinds of recursive `collapsing' rules. Figure \ref{fig:collapse} shows the associated pictures, drawn not for the plane bipartite cubic graphs but for their dual trinities (with an arbitrarily chosen local coloring). The part of the trinity outside of the emerald contour does not change under the collapsing operations. (The red and violet points in the pictures have other neighbors that are not shown.) It is easy to show that if a trinity contains at least four triangles than at least one of the moves can always be carried out.

The rules themselves are $S(T^*,u)=uS(T_0^*,u)$ for the collapsing move on the left and $S(T^*,u)=S(T_1^*,u)+S(T_2^*,u)$ for the one on the right. In the first case, $\rho(T)=\rho(T_0)$ by an application of Corollary \ref{cor:rho} and either Lemma \ref{lem:farkinca} (if we choose to work with the red or the violet graph) or the observation that a double edge has been removed from the emerald graph.

In the second case, separate the points of $T$ by color into the usual sets $R$, $E$, and $V$. Notice that the set $V_1$ of violet points in $T_1$ is obtained from $V$ by identifying $v'$ and $v''$ with a single point $v$. The set of emerald points in both $T_1$ and $T_2$ is $E\setminus\{e\}$. Furthermore, the hypergraph $(V_1,E\setminus\{e\})$ contained in $T_1$ and the hypergraph $(V,E\setminus\{e\})$ contained in $T_2$ are the contraction and deletion, respectively, of the hyperedge $e$ from the hypergraph $(V,E)$ contained in $T$. Therefore we get the desired conclusion from Corollary \ref{cor:rho} and Proposition \ref{pro:delcontr}.
\end{proof}

\end{document}